\documentclass[12pt,a4paper]{article}%
\usepackage{setspace}

\usepackage[margin=1in,a4paper]{geometry}
\usepackage{commath}
\usepackage{amsfonts}
\usepackage{amsthm}
\usepackage{amsmath}
\usepackage{amssymb}
\usepackage{enumerate}
\usepackage[longnamesfirst]{natbib}
\usepackage[bookmarks,colorlinks = true,linkcolor = blue,urlcolor  =
blue,citecolor = blue,anchorcolor = blue]{hyperref}
\usepackage{booktabs}
\usepackage{xr}
\usepackage{graphicx}
\usepackage{multirow}
\usepackage{float}
\usepackage{color}
\usepackage{blkarray}
\usepackage[autostyle]{csquotes}
\usepackage{setspace}%
\usepackage{tikz}
\usetikzlibrary{positioning} 
\usepackage{xcolor}
\providecommand{\U}[1]{\protect \rule{.1in}{.1in}}
\newtheorem{lemma}{Lemma}%[section]
\newtheorem{assumption}{Assumption}%[section]
\newtheorem{proposition}{Proposition}%[section]
\newtheorem{thm}{Theorem}%[section]
\newtheorem{corollary}{Corollary}%[section]
\numberwithin{equation}{section}
\DeclareMathOperator*{\ve}{vec}

\DeclareMathOperator*{\vech}{vech}
\DeclareMathOperator*{\tr}{tr}
\DeclareMathOperator*{\var}{var}

\DeclareMathOperator*{\diag}{diag}

\begin{document}

\title{A Large Confirmatory Dynamic Factor Model for Stock Market Returns in Different Time Zones\thanks{The \texttt{Matlab} code for this paper is available at: \href{https://doi.org/10.13140/RG.2.2.29700.95365}{https://doi.org/10.13140/RG.2.2.29700.95365}. Any remaining errors are our own.} }
\author{Oliver B. Linton\thanks{Faculty of Economics, University of Cambridge. Austin Robinson Building, Sidgwick Avenue, Cambridge, CB3 9DD, UK. Email: \texttt{obl20@cam.ac.uk.}}\\University of Cambridge
\and Haihan Tang\thanks{Corresponding author. International School of Finance, Fudan University. 220 Handan Road, Yangpu District, Shanghai, 200433, China. Email:
\texttt{hhtang@fudan.edu.cn.}}\\Fudan University
\and Jianbin Wu\thanks{School of Business, Nanjing University. 22 Hankou Road, Gulou District, Nanjing, 210093, China. Email: \texttt{wujianbin@nju.edu.cn.} }\\Nanjing University }
\date{\today }
\maketitle

\begin{abstract}
\noindent We propose a confirmatory dynamic factor model for a large number of stocks whose returns are observed daily  across multiple time zones. The model has a global factor and a continental factor that both drive the individual stock return series. We propose two estimators of the model: a quasi-maximum likelihood estimator (QML-just-identified), and an improved estimator based on an Expectation Maximization (EM) algorithm (QML-all-res). Our estimators are consistent and asymptotically normal under the large approximate factor model setting. In particular, the asymptotic distributions of QML-all-res are the same as those of the infeasible OLS estimators that treat factors as known and utilize all the restrictions on the parameters of the model. We apply the model to MSCI equity indices of 42 developed and emerging markets, and find that most markets are more integrated when the CBOE Volatility Index (VIX) is high.

\begin{description}
\item[Keywords:] Daily Global Stock Market Returns; Time-Zone Differences; Confirmatory Dynamic Factor Models; Quasi Maximum Likelihood; EM Algorithm.

\item[JEL classification] C32; C55; C58; G15.

\end{description}
\end{abstract}

\section{Introduction}

Correlations among stock returns in different markets are an important ingredient of the strategic portfolio allocation including international stocks and for understanding the degree of segmentation of international capital markets. The last three decades have witnessed a heightening interest in measuring and modelling such correlations, whether under the rubric of stock market integration,  
international stock co-movement, interdependence or otherwise (\cite{gagnonkarolyi2006}, \cite{sharmaseth2012}). %\cite{gagnonkarolyi2006} and \cite{sharmaseth2012} have carefully reviewed the early literature and categorized these studies according to methodologies, data sets, and findings. 

The world's stock markets are separated in time by substantial time-zone differences, to the extent that for example the American and main Asian markets do not overlap at all (although the European markets overlap a little with both the American and Asian markets). Table \ref{table top ten stock exchange} lists the trading hours of the world's top ten stock exchanges in terms of market capitalization as of September 2023. This shows the overlaps and lack of overlaps between them. Whenever using daily data to measure correlations, researchers have to address the issue of \textit{non-synchronous trading}. This is because the closing prices of the American markets, say, contain a lot of information that could never have been in the closing prices of the Asian markets for the same calendar days (\cite{SCHOTMAN2006462}) and vice versa. Ignoring non-synchronous trading issues can bring difficulties and biases to research and to stock trading. Specifically, the correlation between stocks may be underestimated (\cite{burns1998correlations}).

\begin{table}[ptb]
\centering
\begin{tabular}[c]{clllc}%
\toprule[0.3mm] \toprule[0.3mm] & Stock Exchange (S.E.) & Market Cap & Trading Hours & Country\\
&  & (US Trillions) & (Beijing Time) & \\
\midrule[0.3mm] 
1 & New York S.E. & \$25.24 & 21:30-04:00($+1$)$^{\dagger}$ & US\\
2 & NASDAQ S.E. & \$20.58 & 21:30-04:00($+1$)$^{\dagger}$ & US\\
3 & Shanghai S.E. & \$6.6 & 09:30-11:30, 13:00-15:00 & China\\
4 & Euronext S.E. & \$6.26 & 15:00-23:30$^{\dagger}$   & Euro. Countries\\
5 & Tokyo S.E. & \$5.75 & 08:00-10:30, 11:30-14:00 & Japan\\
6 & Shenzhen S.E. & \$4.38 & 09:30-11:30, 13:00-15:00 & China\\
7 & Hong Kong S.E. & \$4.1 & 09:30-12:00, 13:00-16:00 & China\\
8 & India National S.E. & \$3.59 & 11:45-18:00 & India\\
9 & London S.E. & \$3.42 & 15:00-23:30$^{\dagger}$ & UK\\
10 & Saudi S.E. & \$3.06 & 15:00-20:00 & Saudi Arabia\\
\bottomrule[0.3mm] \bottomrule[0.3mm] 
\end{tabular}
\caption{{\protect\small As of September 2023. ($+1$) means $+$one day. $^{\dagger}$ means summer time. The trading hours of the markets of Euronext S.E. vary and the table reports those of the Paris S.E.}}%
\label{table top ten stock exchange}%
\end{table}

The objective of this paper is to develop a framework to model correlations of daily stock returns in different markets across multiple time zones. This allows one to analyze international stock returns as if they were measured at the same time, and so allows empirical researchers to eliminate the bias that would otherwise be manifested in whatever question they are addressing. The machinery will be a dynamic factor model, a special case of the generalized dynamic factor model first proposed by \cite{fornietal2000}. To make the model tractable, we make the following simplifying assumption: All the markets belong to either Asia, Europe, or America.

We suppose that each logarithmic 24-hour close-to-close return follows a dynamic factor model driven by a global factor and a continental factor. This model reflects a situation in which international information represented by the global factor affects all the three continents simultaneously, but is only revealed in the returns of the three continents sequentially as their markets open in turn and trade on the new information (\citet[][p.235]{kochkoch1991}). Local information represented by a continental factor accumulated since the last closure of that continent will also have an impact on the logarithmic 24-hour close-to-close returns of those markets in our framework. The modelling idea is that the global and continental factors can only impact the returns contemporaneously, consistent with the weak-form efficient market hypothesis. The dynamics of the factors can be consistent with a weak-form efficient market (if the dynamic parameter of the global factor is zero); it can also be consistent with market microstructure or sluggish responses of prices introducing some short-term auto-correlation that is manifested in the dynamic evolution of the global factor. 

Our model is more like a \textit{structural} or \textit{confirmatory} factor model,\footnote{See \cite{joreskog1969} and \cite{rubinthayer1982} for early discussions of confirmatory factor models.} where all the factors and loadings have clear financial interpretations. On the other hand, the identification approaches commonly used in the literature are more like reduced forms and do not have such clear financial interpretations. Although our identification restrictions are not frequently utilized in the factor analysis literature, they are commonly employed in the field of macroeconomic nowcasting. In macroeconomic nowcasting, certain macroeconomic indices often exhibit zero loadings on specific factors (cf. \cite{giannonereichlinsmall2008}, \cite{aruobadieboldscotti2009}, %\cite{2010Nowcasting},
\cite{banburaetal2013}, and \cite{BokBrandyn2018}). Notably, the methodology outlined by \cite{BokBrandyn2018} forms the basis for the New York Fed staff nowcast technical paper (\cite{almuzara2023new}).

We propose two estimators of the parameter matrices, $\Lambda, M, \Sigma_{ee}$, of the model: a quasi-maximum likelihood estimator (QML-just-identified), and an improved estimator (QML-all-res). The QML-all-res estimator is built upon the QML-just-identified estimator and takes all the restrictions implied by our model into account. The idea is that all these restrictions could be incorporated in an Expectation-Maximization (EM) algorithm. The numbers of parameters in $\Lambda$, $M$ (symmetric) and $\Sigma_{ee}$ (diagonal), without any imposed restriction of our model, are $6N\times 14$, $14\times 15/2$, and $6N$, respectively (see Section \ref{sec Model} for details). We do not incur the curse of dimensionality in estimation because we ignore the auto-correlation of the stacked factors, as well as the auto- and cross-correlation of the stacked idiosyncractic components, when setting up the quasi-log-likelihood. Because of this, the EM algorithm of QML-all-res consists of analytic formulas for factors and most parameters, and is viable in large dimensions.

A key issue we address in our paper is identification. We show explicitly how to impose some of the restrictions implied by our model to guarantee identification of the QML-just-identified estimator. In particular, we impose a subset of restrictions, specifically $14^{2}$ of them, because there are 14 factors in the two-day representation (kind of a static form) of our model. Since these $14^{2}$ restrictions cannot be written compactly in a matrix form, it took us a considerable amount of work to derive the large sample results of the QML-just-identified estimator. Let $N$ and $T$ denote the cross-sectional and temporal dimensions, respectively (see Section \ref{sec Model} for details). We prove that as $N,T\rightarrow \infty$, QML-just-identified is consistent, and asymptotically normal if $\sqrt{T}/N\to 0$.  As $N,T\rightarrow \infty$, QML-all-res is asymptotically normal if $\log N/T\to 0$ and $\sqrt{T}/N\to 0$. In particular, the asymptotic distributions of QML-all-res are the same as those of the infeasible OLS estimators that treat factors as known and utilize all the restrictions of the parameters of the model. We provide consistent standard errors for the QML-all-res estimator under a setting justifiable in finance.  We next show that the estimated factors formed using QML-all-res are asymptotically normal as $N,T\rightarrow \infty$ if $\sqrt{N}/T\to 0$, up to logarithmic terms. Uniform consistency for the estimated factors is also obtained. %Last, 
All the large sample results are obtained under a large approximate factor model setting.

We apply our model to MSCI equity indices of 42 developed and emerging markets. %(the cross-sectional dimension of the two-day representation is 424, while the time series dimension $T=1179$). 
Taking Asia-Pacific as an example, we find that the global factor has the largest loadings during the Asian trading period (a sub-period of a day). China's mainland and Hong Kong of China have large loadings on the continental factor, whereas Japan has minimal loadings on the continental factor. In addition, we adopt a methodology similar to \cite{2009Global} to assess market integration. The US demonstrates the highest level of global integration but minimal regional integration. We also find evidence that the weak-form efficient market hypothesis does not hold across the globe and it takes more than a $1/3$ day for some international news to fully unfold or dissipate.

Furthermore, we investigate whether markets become more integrated during periods of high CBOE Volatility Index (VIX). In particular, we divide the sample into four sub-samples: Before-HighVix, Before-LowVix, After-HighVix, and After-LowVix. Before-HighVix is the intersection of the before and high-VIX sub-samples, and so forth. By comparing the results between Before-HighVix and Before-LowVix (as well as between After-HighVix and After-LowVix), we find most markets exhibit higher levels of integration during periods of higher VIX. This heightened integration is primarily driven by the increase of magnitudes of loadings of the global factor during the Asian sub-period. This also means that, during periods of higher VIX, the returns of the Asian markets impart trading signals to subsequent markets. %contain stronger trading signals \textcolor{red}{for other markets (Jianbin: should be ``than other markets?")}.

\subsection{Literature Review}

There are several general approaches to address the non-synchronous trading across the international markets. The first approach is to use returns of relatively low frequencies: two-day returns (\cite{KristinRoberto2002}, %\cite{corsetti2005some}, 
\cite{DUNGEY2015271}), weekly returns (%\cite{cappiello2006asymmetric}, 
\cite{BekaertHodrickZhang2009}, \cite{bekaert2014global}%\cite{nictoi2019drives}
), and monthly returns (\cite{bekaert2019}, \cite{RapachStraussZhou2013JF}).

The second approach involves use of intraday prices as the ``pseudo" closing prices of one market to better match the closing prices of another market closing at an earlier time (\cite{SCHOTMAN2006462}, \cite{savva2009spillovers}%\cite{gagnon2009information}
). For example, \cite{SCHOTMAN2006462} used an intraday price of the German market as its pseudo closing price when studying the German and Hungarian markets. The timing of the intraday price is chosen to match the closing time of the Hungarian market. %\cite{savva2009spillovers} and \cite{gagnon2009information} also used this approach.

The third approach is attributed to \cite{burns1998correlations}, who proposed a Vector Moving Average (VMA) structure for handling the non-synchronous trading. This method has been widely utilized in various studies (%\cite{acharya2012capital}, 
\cite{BENSAIDA2019}, \cite{OPIE2020}).%, \cite{engle2015systemic}, \cite{engle2016dynamic} etc.).

The fourth approach involves using contemporaneous or lagged variables accordingly to reflect the time zone differences. For instance, \citet[][p.11]{bekaert2023risk} wrote, ``For the euro area, JP (Japanese) and EA (European) shocks that materialize before or during the European opening hours enter contemporaneously while the other shocks as well as the US shocks enter the information set on the next trading day." %\cite{2009Global}, 
Also see \cite{CAI2009} and \cite{lehkonen2015stock}. 

We have a few remarks about these approaches. Using returns of relatively low frequencies as in the first approach prevents researchers from studying stock co-movement at higher frequencies. The second approach utilizes intraday data which might be unavailable or costly. The third approach is sensitive to the ordering of returns. For instance, using a vector of returns of Asia, Europe, and the US of the same calendar day might produce a different result from using a vector of returns of Europe and the US of a calendar day, together with the return of Asia of the next day. The fourth approach is a very flexible way to address the non-synchronous trading issue, but does not allow one to identify the source of variation or its relative impacts as is possible in our framework.

The methodological approach of having global and continental factors, in some respects, is inspired by \cite{lin1994bulls} who investigated the correlation between the Nikkei 225 and S\&P500 indices. They proposed a signal-extraction model containing a global factor and a local factor for the returns, but the model cannot generalize to a case of many stocks in multiple time zones easily. Our approach resembles that of \cite{koseotrokwhiteman2003} who modeled 60 markets' yearly macroeconomic aggregates (output, consumption and investment) using a dynamic factor model consisting of global, regional, and market-specific factors. The difference is that we work with daily stock returns and have the feature of sequential revelation of the international information, measured at 1/3-day frequency, in the returns of the three continents. Likewise, our approach is related to \cite{hallinliska2011} and
\cite{barigozzihallinsoccorsi2019}, who considered generalized dynamic factor models in the presence of blocks.\footnote{\cite{barigozzihallinsoccorsi2019} analyzed the US, European and Japanese stock returns. Their model allows for not only global or continental factors, but also factors affecting any two of the three markets. Nevertheless they found no global factors in their sample.}
%and to nowcasting\footnote{In the nowcasting literature, researchers use factor models to extract the information contained in the data at higher frequencies than the target variable in order to forecast the target variable.} (\cite{giannonereichlinsmall2008}, \cite{aruobadieboldscotti2009}, \cite{2010Nowcasting}, \cite{banburaetal2013}, \cite{BokBrandyn2018}, \cite{almuzara2023new} etc.).

%In contrast, in the literature of block-specific factors and nowcasting, factor models are dimension-reducing tools with no economic interpretations (i.e., \textit{exploratory} factor analysis). %In some sense, our model belongs to the class of structural dynamic factor models (\cite{stockwatson2016}).

%Likewise, our approach is related  In the nowcasting literature, researchers use factor models to extract the information contained in the data at higher frequencies than the target variable in order to forecast the target variable. Here, if we make some additional, mild assumptions on the data generating process, we could also obtain the nowcast of a stock's price level after the market closure; this is the similarity.

On the theoretical side, research about estimation of large factor models via
the likelihood approach has matured over the last decade; \cite{barigozzi2023}
provided the most recent and critical review on this topic. The likelihood
approach enjoys several advantages compared to the principal components method (\citet[p.204]{banburaetal2013}, \citet[p.37]{barigozzi2023}). The most prominent advantage in our context (i.e., confirmatory factor models) is its flexibility to incorporate (nonlinear) restrictions on the factor loadings and covariance matrices, something cannot be easily achieved in the principal components method (\citet[p.431]{stockwatson2016}).

\cite{dozgiannonereichlin2012} established an average rate of convergence of the estimated factors using a quasi maximum likelihood (QML) estimator via the Kalman smoother. %\footnote{In fact, \cite{dozgiannonereichlin2012} called their estimator the quasi-maximum likelihood estimator (QMLE) instead of the MLE. We re-label it as the MLE since we shall reserve the phrase QMLE for another purpose to be made specific shortly. Similarly for \cite{baili2012}'s QMLE in the next paragraph.} 
There is a rotation matrix attached to the estimated
factors as the authors did not address identification of factor models.
%Nor did they derive consistency for the estimated factor loadings, or the limiting distributions of any estimate.
Moreover, they did not provide consistency or the limiting distributions of the loadings. In an important paper, \cite{baili2012} took a different approach to study the large exact factor models. \cite{baili2012} obtained consistency, the rates of convergence, and the limiting distributions of the QML estimator of the factor loadings, idiosyncratic variances, and sample covariance matrix of factors. Factors are then estimated via a generalized least squares (GLS) method. \cite{baili2016} generalized the results of \cite{baili2012} to the large approximate factor models.

Instead of maximizing a likelihood and finding the QML estimator (the global maximum), people usually use EM algorithms to estimate models in practice
(\cite{watsonengle1983}, \cite{quahsargent1992},
\cite{dozgiannonereichlin2012}, \cite{rubinthayer1982}, \cite{baili2012},
\cite{baili2016} etc.). Since an EM algorithm runs only for a finite number of
iterations and converges to a local maximum under some conditions when the number of iterations goes to infinity, strictly speaking the estimate obtained by an EM algorithm is only an approximation to the QML estimator. To the best of our knowledge, only \cite{barigozziluciani2022} have proved that the approximation error is asymptotically negligible but under a different EM algorithm and a different identification scheme from ours (\citet[p.19]{barigozziluciani2022}). %However, in a breakthrough study, \cite{barigozziluciani2022} showed that the estimate obtained by an EM algorithm, under the assumptions of their paper, converges to the MLE fast enough so that they are asymptotically equivalent.
%Theorem \ref{thm qmleres1} in Section \ref{sec theory qmle all res} proves that for the QML-all-res estimator, we need not worry about the issue of a finite number of iterations. 

%We propose several estimators of our model (\ref{random29}): the MLE-one day, the QMLE-res, the QMLE, the QMLE-md, and the Bayesian. The MLE-one day estimator is the usual MLE estimator of our model. The QMLE-res estimator is the MLE estimator of the two-day representation of our model while maintaining the working independence hypothesis (see Section \ref{sec two-day representation}); in this article we shall refer a likelihood-based estimator obtained under the working independence hypothesis as the QMLE rather than the MLE. The QMLE estimator differs from the QMLE-res in the sense that only a specific subset of restrictions implied by our model is imposed. Since the QMLE is inefficient, we propose an improved estimator, the QMLE-md, which uses the QMLE in the first step and incorporates an additional finite number of restrictions implied by our model via the minimum distance method in the second step. Last, the Bayesian estimator uses the Gibbs sampling to estimate the model (\ref{random29}). However, the Gibbs sampling is computationally intensive and feasible only for a not-so-large number of entities.

The large sample results of the aforementioned studies are not applicable to
our model because the proofs of these results are identification-scheme
dependent. In particular, \cite{baili2012}, \cite{baili2016} established their
results under five popular identification schemes, none of which is consistent
with our model. 

%We contribute to methodology by providing a new modelling framework for daily global stock market returns. Our framework could easily handle a large number of stocks and at the same time take into account the time-zone differences. Under a mild fixed-signs assumption, our model is identified and has a structural interpretation. We also contribute to theory by deriving the asymptotic results of the QMLE and the QMLE-md. The machinery is based on the theoretical results of \cite{baili2012}, but we demonstrate how one could obtain their results for almost \textit{any} identified dynamic factor model. This is an important contribution as many dynamic factor models, like ours, are motivated by different economic theories and might not be compatible with the five identification schemes of \cite{baili2012}. We last contribute to the applied literature by proposing several practically usable estimators and validate their performances in the Monte Carlo simulations. When applying our model to two real data sets, we draw some new insights about linkages among different stock markets.

\subsection{Roadmap and Notation}

The rest of the paper is structured as follows. Section \ref{sec Model} introduces the model. We explain our estimators and provide their large-sample theories in Sections \ref{sec Estimation} and \ref{sec Theories}, respectively. Section \ref{sec Monte Carlo} conducts the Monte Carlo simulations to assess the performance of QML-all-res, and Section \ref{sec application} presents empirical applications. %Section \ref{sec conclusion} concludes. 
Additional materials and proofs are in the
Online Supplement (OS in what follows).

%\subsection{Notation}
%\label{sec notation}

%Let $\mathbb{R}^{n}$ denote the $n$-dimensional Euclidean space. For $x\in\mathbb{R}^{n}$, let $\Vert x\Vert:=\sqrt{\sum_{i=1}^{n}x_{i}^{2}}$ and $\|x\|_{\infty}:=\max_{1\leq i\leq n}|x_{i}|$ denote the Euclidean ($\ell_{2}$) and element-wise maximum ($\ell_{\infty}$) norms, respectively. 
Given a vector $x$, let $\Vert x\Vert$ and $\diag(x)$ denote the Euclidean  norm and  the diagonal matrix whose diagonal elements are $x$, respectively. %We use $P(\cdot|\cdot)$ to denote a conditional probability density function.
Notation $\ve A$ means the vector obtained by stacking columns of matrix $A$. %The \textit{commutation matrix} $K_{m,n}$ is an $mn\times mn$ \textit{orthogonal} matrix which translates $\ve A$ to $\ve(A^{\intercal})$, i.e., $\ve(A^{\intercal})=K_{m,n}\ve(A)$. 
If $A$ is an $n\times n$ symmetric  matrix, its supradiagonal elements are redundant in the sense that they can be deduced from symmetry. If we eliminate these redundant elements from $\ve A$, we obtain an  $n(n+1)/2\times1$ vector, $\vech A$. They are related by the full-column-rank \textit{duplication matrix} $D_{n}$: $\ve A=D_{n}\vech A$. Conversely, $\vech A=D_{n}^{+}\ve A$, where $D_{n}^{+}$ is  the Moore-Penrose generalized inverse of $D_{n}$. Let $\lambda_{\max}(\cdot)$ and $\lambda_{\min}(\cdot) $ denote the maximum and minimum eigenvalues of some real symmetric matrix, respectively. For any real matrix $A$, let $\Vert A\Vert$ denote the Frobenius  norm. 
Given a square matrix $A$, $\diag(A)$ sets the off-diagonal elements of $A$ to zero. 
Landau notation  should be interpreted in the sense that $N,T\rightarrow\infty$ jointly. Let $C, C^{\prime}$ to denote absolute positive constants; %(i.e., constants independent of anything which is a function of $N$ and/or $T$)
identities of such $C$s might change from one place to another. Let $a\wedge b$ denote $\min (a,b)$.

\section{Model}
\label{sec Model}

\subsection{Setup}
\label{sec model setup}

\begin{figure}[ptb]
\centering
\begin{tikzpicture}[domain=-7:7]
\draw[black,->](-7,1.1)--(7,1.1);
\draw  [densely dotted,line width=0.5pt](-6,-2)--(-6,2);
\draw  [densely dotted,line width=0.5pt](-4,-2)--(-4,2);
\draw  [densely dotted,line width=0.5pt](-2,-2)--(-2,2);
\draw  [densely dotted,line width=0.5pt](0,-2)--(0,2);
\draw  [densely dotted,line width=0.5pt](2,-2)--(2,2);
\draw  [densely dotted,line width=0.5pt](4,-2)--(4,2);
\draw  [densely dotted,line width=0.5pt](6,-2)--(6,2);
\draw  [<->,line width=0.5pt](-6,0.2)--(-4.03,0.2);
\draw  [<->,line width=0.5pt](-4,0.2)--(-2.03,0.2);
\draw  [<->,line width=0.5pt](-2,0.2)--(-0.03,0.2);
\draw  [<->,line width=0.5pt](0,0.2)--(1.97,0.2);
\draw  [<->,line width=0.5pt](2,0.2)--(3.97,0.2);
\draw  [<->,line width=0.5pt](4,0.2)--(5.97,0.2);
\draw  [<->,line width=0.5pt](-2,-1)--(4,-1);
\node[align=left,font=\fontsize{8}{10}\selectfont] at (-6,2.5){Amer. closes\\Day $s-2$};
\node[align=left,font=\fontsize{8}{10}\selectfont] at (-4,2.5){Asia closes\\Day $s-1$};
\node[align=left,font=\fontsize{8}{10}\selectfont] at (-2,2.5){Euro. closes\\Day $s-1$};
\node[align=left,font=\fontsize{8}{10}\selectfont] at (0,2.5){Amer. closes\\Day $s-1$};
\node[align=left,font=\fontsize{8}{10}\selectfont] at (2,2.5){Asia closes\\Day $s$};
\node[align=left,font=\fontsize{8}{10}\selectfont] at (4,2.5){Euro. closes\\Day $s$};
\node[align=left,font=\fontsize{8}{10}\selectfont] at (6,2.5){Amer. closes\\Day $s$};
\node[align=left,font=\fontsize{8}{10}\selectfont] at (-5,-0.3){$f^{G}_{s-1,\text{as\_time}}$};
\node[align=left,font=\fontsize{8}{10}\selectfont] at (-3,-0.3){$f^{G}_{s-1,\text{eu\_time}}$};
\node[align=left,font=\fontsize{8}{10}\selectfont] at (-1,-0.3){$f^{G}_{s-1,\text{am\_time}}$};
\node[align=left,font=\fontsize{8}{10}\selectfont] at (1,-0.3){$f^{G}_{s,\text{as\_time}}$};
\node[align=left,font=\fontsize{8}{10}\selectfont] at (3,-0.3){$f^{G}_{s,\text{eu\_time}}$};
\node[align=left,font=\fontsize{8}{10}\selectfont] at (5,-0.3){$f^{G}_{s,\text{am\_time}}$};
\node[align=left,font=\fontsize{8}{10}\selectfont] at (1,-1.3){$\boldsymbol{R}_{\text{eu},s}$};
\node[align=left,font=\fontsize{8}{10}\selectfont] at (1,-1.6){${f}_{\text{eu},s}$};
\node[align=left,font=\fontsize{8}{10}\selectfont] at (1,-1.9){$\boldsymbol{\epsilon}_{\text{eu},s}$};
\end{tikzpicture}
\caption{Returns and factors}%
\label{fig:return}%
\end{figure}
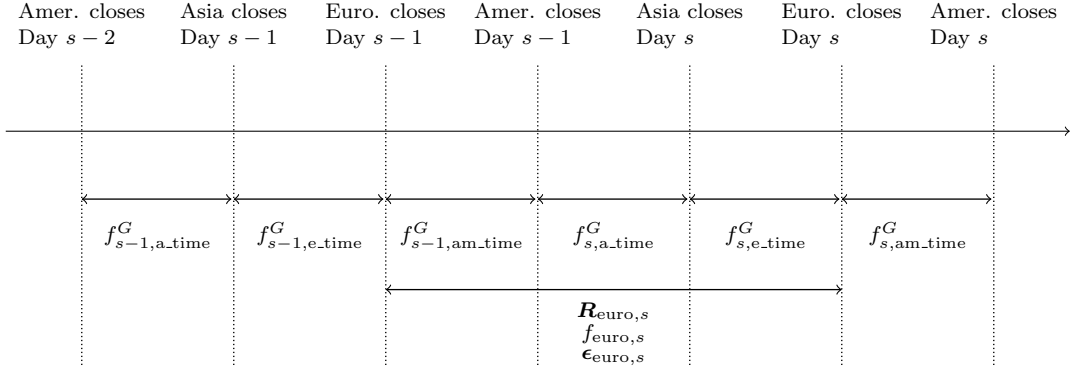

Our model is based on the closing prices of various stocks in three continents, namely Asia, Europe and America. All the closing prices of the stocks of a continent are assumed to be observed at their continental closing time. On a calendar day, the Asian markets close first, followed by the European markets, and finally the American markets close. Let $\boldsymbol{R}_{\text{as},s},\boldsymbol{R}_{\text{eu},s}$, and $\boldsymbol{R}_{\text{am},s}$ represent the logarithmic 24-hour close-to-close returns on day $s$ for stocks traded in Asia, Europe, and America, respectively. Without loss of generality, we assume that each continent has the same number of stocks: $\boldsymbol{R}_{\text{as},s},\boldsymbol{R}_{\text{eu},s}$, and $\boldsymbol{R}_{\text{am},s}$ are all $N\times1$ dimensional vectors. We formulate the following model:
\begin{align}
& \boldsymbol{R}_{\text{as},s}=\boldsymbol{\beta}_{\text{as,as\_time}}^{G}f^{G}{}_{s,\text{as\_time}}+\boldsymbol{\beta}_{\text{as,am\_time}}^{G}f^{G}{}_{s-1,\text{am\_time}}+\boldsymbol{\beta}_{\text{as,eu\_time}}^{G}f^{G}{}_{s-1,\text{eu\_time}}+\boldsymbol{\beta}_{\text{as}}f{}_{\text{as},s}+\boldsymbol{\epsilon}_{\text{as},s},\nonumber\\
& \boldsymbol{R}_{\text{eu},s}=\boldsymbol{\beta}_{\text{eu,eu\_time}}^{G}f^{G}{}_{s,\text{eu\_time}}+\boldsymbol{\beta}_{\text{eu,as\_time}}^{G}f^{G}{}_{s,\text{as\_time}}+\boldsymbol{\beta}_{\text{eu,am\_time}}^{G}f^{G}{}_{s-1,\text{am\_time}}+\boldsymbol{\beta}_{\text{eu}}f{}_{\text{eu},s}+\boldsymbol{\epsilon}_{\text{eu},s},\nonumber\\
& \boldsymbol{R}_{\text{am},s}=\boldsymbol{\beta}_{\text{am,am\_time}}^{G}f^{G}_{s,\text{am\_time}}+\boldsymbol{\beta}_{\text{am,eu\_time}}^{G}f^{G}_{s,\text{eu\_time}}+\boldsymbol{\beta}_{\text{am,as\_time}}^{G}f^{G}_{s,\text{as\_time}}+\boldsymbol{\beta}_{\text{am}}f{}_{\text{am},s}+\boldsymbol{\epsilon}_{\text{am},s},\label{main_r6}
\end{align}
where $f^{G}_{s,\text{as\_time}},f^{G}_{s,\text{eu\_time}}$, and $f^{G}_{s,\text{am\_time}}$ are the values of the scalar unobserved global factor in the three sub-periods of day $s$ and $\boldsymbol{\beta}^{G}_{\cdot,\cdot}$s are the corresponding factor loadings. The three sub-periods are: the Asian sub-period (from the American close on day $s-1$ to the Asian close on day $s$), the European sub-period (from the Asian close on day $s$ to the European close on day $s$), and the American sub-period (from the European close on day $s$ to the American close on day $s$). Here, $f{}_{\text{as},s},$ $f{}_{\text{eu},s}$, and $f{}_{\text{am},s}$ are the scalar unobserved 24-hour Asian, European, and American factors on day $s$, respectively, and $\boldsymbol{\beta}_{\text{as}},\boldsymbol{\beta}_{\text{eu}}$, and $\boldsymbol{\beta}_{\text{am}}$ are the corresponding factor loadings. $\boldsymbol{\epsilon}_{\text{as},s},\boldsymbol{\epsilon}_{\text{eu},s}$, and $\boldsymbol{\epsilon}_{\text{am},s}$ are the 24-hour Asian, European, and American idiosyncratic components, respectively. The global factor, along with $\boldsymbol{R}_{\text{eu},s},f{}_{\text{eu},s}$ and $\boldsymbol{\epsilon}_{\text{eu},s}$, is illustrated in Figure \ref{fig:return}. We suppose for simplicity that the global factor follows an AR(1) process:
\begin{align}
f^{G}_{s,\text{am\_time}}  &  =\phi f^{G}_{s,\text{eu\_time}}+\eta_{s,\text{am\_time}},\nonumber\\
f^{G}_{s,\text{eu\_time}}  &  =\phi f^{G}_{s,\text{as\_time}}+\eta_{s,\text{eu\_time}},\label{main_r5}\\
f^{G}_{s,\text{as\_time}}  &  =\phi f^{G}{}_{s-1,\text{am\_time}}+\eta_{s,\text{as\_time}},\nonumber
\end{align}
%where $f^{G}{}_{s,\text{b\_time}}=0$ for $s\leq0$, $b=\text{as},\text{eu},\text{am}$.
%
%Here we specify that the global factor follows an AR(1) process. 
where the parameters $\phi$ are the same in each equation.\footnote{Display (\ref{main_r5}) implies
\begin{align*}
f^{G}_{s,\text{am\_time}}=\phi^3 f^{G}_{s-1,\text{am\_time}}+\phi^2\eta_{s,\text{as\_time}}+\phi\eta_{s,\text{eu\_time}}+\eta_{s,\text{am\_time}}.
\end{align*}
We do not directly use this representation though. Instead, we shall cast our model in a two-day form, in which case a better way to describe the dynamics of the global factor is (\ref{res_r91}).}
Alternatively, we could model the global factor as a VAR(1) process with respect to the global factor vector $(f^{G}_{s,\text{am\_time}}, f^{G}_{s,\text{eu\_time}}, f^{G}_{s,\text{as\_time}})^{\intercal}$ to allow for greater flexibility. Since we focus on drawing financial implications from the empirical applications rather than contructing a model to better fit data, an AR(1) specification may yield more robust estimates of parameters. %However, our preliminary investigations suggest that with a VAR(1) specification, using slightly different time periods of the data could lead to notable changes in the parameter estimates, particularly those of the factor loading matrix. 

Our model postulates that stock returns are driven by both a global factor and a continental factor. In particular, with the logarithmic 24-hour close-to-close returns of the three continents, we can identify the global factor within a sub-period (roughly a 1/3 day), since we observe the closing prices of some continent at the end of that sub-period. However, we are unable to identify the Asian continental factor, say, within the European sub-period since the closing prices of Europe do not provide any information about the Asian continental factor (see Assumption \ref{assu eta}(iv)),
%. For instance, the closing prices of the European or American markets do not contain any information about the Asian continental factor,
but they do contain information about the global factor. The correlation of stock returns in different continents is solely due to the international news represented by the global factor.\footnote{We explore the issue of including multiple global factors in our model in Section \ref{sec more global factors} of the OS.}

There is a fairly long history of nearly 24-hour trading of S\&P 500 index futures. Our model can be modified to include this additional information. We can treat the 24-hour S\&P 500 index future as a special asset, and utilize its prices at the Asian, European, and American closing times. Then we can calculate its logarithmic 1/3-day close-to-close returns of each sub-period of a day. Such a return will have a non-zero loading on the global factor for the corresponding sub-period. For example, the return for the Asian sub-period of a S\&P 500 index future has a non-zero loading on $f^G_{s,\text{as\_time}}$ while zero loadings on $f^G_{s,\text{eu\_time}}$ and $f^G_{s,\text{am\_time}}$. In addition, all these logarithmic 1/3-day close-to-close returns are affected by the American continental factor and not by the other two continental factors because the index futures concern S\&P 500 index.\footnote{If we have data on 24-hour index futures of the Asian and European continents, we could identify all the continental factors within a sub-period of a day.} The dynamic equation (\ref{main_r5}) of the global factor  can remain unchanged. Estimation of this modified model should be straightforward.  %However, in the current version, we only consider a very basic model that relies on the closing price. Due to the potential complexity it may introduce to our theoretical framework, we have decided to defer this extension to future work.

\subsection{The Two-Day Representation}
\label{sec two day representation}

We shall cast our model in a two-day form where every two calendar days is treated as a single time unit. Calendar days $1,2$ are treated as time unit 1, and calendar days $3,4$ are treated as time unit $2$, and so forth. Estimation and large sample theories are based on this two-day representation. Specifically, we have
\begin{align}
\underbrace{\boldsymbol{y}_{t}}_{6N\times1}=\underbrace{\Lambda}_{6N\times 14}\underbrace{\boldsymbol{f}_{t}}_{14\times1}+\underbrace{\boldsymbol{e}_{t}}_{6N\times1} ,\label{r64}%
\end{align}
for $t=1,2,\ldots,T,$ where
%\[\boldsymbol{y}_{t}:=\left[\begin{array}[c]{c}%
%\boldsymbol{R}_{\text{as},s}\\
%\boldsymbol{R}_{\text{eu},s}\\
%\boldsymbol{R}_{\text{am},s}\\
%\boldsymbol{R}_{\text{as},s+1}\\
%\boldsymbol{R}_{\text{eu},s+1}\\
%\boldsymbol{R}_{\text{am},s+1}%
%\end{array}
%\right]  \qquad\boldsymbol{e}_{t}:=\left[
%\begin{array}[c]{c}%
%\boldsymbol{\epsilon}_{\text{as},s}\\
%\boldsymbol{\epsilon}_{\text{eu},s}\\
%\boldsymbol{\epsilon}_{\text{am},s}\\
%\boldsymbol{\epsilon}_{\text{as},s+1}\\
%\boldsymbol{\epsilon}_{\text{eu},s+1}\\
%\boldsymbol{\epsilon}_{\text{am},s+1}\\
%\end{array}
%\right]  \qquad\boldsymbol{f}_{t}:=\left[
%\begin{array}[c]{c}%
%f^{G}{}_{s+1,\text{am\_time}}\\
%f^{G}{}_{s+1,\text{eu\_time}}\\
%f^{G}{}_{s+1,\text{as\_time}}\\
%f^{G}{}_{s,\text{am\_time}}\\
%f^{G}{}_{s,\text{eu\_time}}\\
%f^{G}{}_{s,\text{as\_time}}\\
%f^{G}{}_{s-1,\text{am\_time}}\\
%f^{G}{}_{s-1,\text{eu\_time}}\\
%f_{\text{am},s+1}\\
%f_{\text{eu},s+1}\\
%f_{\text{as},s+1}\\
%f_{\text{am},s}\\
%f_{\text{eu},s}\\
%f_{\text{as},s}%
%\end{array}\right],\]
%The loading matrix $\Lambda$ is%

\[\boldsymbol{y}_{t}:=\left[
\begin{array}[c]{c}%
\boldsymbol{R}_{\text{as},2(t-1)+1}\\
\boldsymbol{R}_{\text{eu},2(t-1)+1}\\
\boldsymbol{R}_{\text{am},2(t-1)+1}\\
\boldsymbol{R}_{\text{as},2t}\\
\boldsymbol{R}_{\text{eu},2t}\\
\boldsymbol{R}_{\text{am},2t}%
\end{array}
\right]  \qquad\boldsymbol{e}_{t}:=\left[
\begin{array}[c]{c}%
\boldsymbol{\epsilon}_{\text{as},2(t-1)+1}\\
\boldsymbol{\epsilon}_{\text{eu},2(t-1)+1}\\
\boldsymbol{\epsilon}_{\text{am},2(t-1)+1}\\
\boldsymbol{\epsilon}_{\text{as},2t}\\
\boldsymbol{\epsilon}_{\text{eu},2t}\\
\boldsymbol{\epsilon}_{\text{am},2t}\\
\end{array}
\right]  \qquad\boldsymbol{f}_{t}:=\left[
\begin{array}[c]{c}%
f^{G}{}_{2t,\text{am\_time}}\\
f^{G}{}_{2t,\text{eu\_time}}\\
f^{G}{}_{2t,\text{as\_time}}\\
f^{G}{}_{2(t-1)+1,\text{am\_time}}\\
f^{G}{}_{2(t-1)+1,\text{eu\_time}}\\
f^{G}{}_{2(t-1)+1,\text{as\_time}}\\
f^{G}{}_{2(t-1),\text{am\_time}}\\
f^{G}{}_{2(t-1),\text{eu\_time}}\\
f_{\text{am},2t}\\
f_{\text{eu},2t}\\
f_{\text{as},2t}\\
f_{\text{am},2(t-1)+1}\\
f_{\text{eu},2(t-1)+1}\\
f_{\text{as},2(t-1)+1}%
\end{array}\right]  ,\]

\begin{align}
& \Lambda := \nonumber \\
& \resizebox{\linewidth}{!}{$\left[\begin{array}{cccccccccccccc}
\mathbf{0} & \mathbf{0} & \mathbf{0} & \mathbf{0} & \mathbf{0} &
\boldsymbol{\beta}^{G}_{\text{as,as}} & \boldsymbol{\beta}^{G}_{\text{as,am}} &
\boldsymbol{\beta}^{G}_{\text{as,eu}} & \mathbf{0} & \mathbf{0} & \mathbf{0} &
\mathbf{0} & \mathbf{0} & \boldsymbol{\beta}_{\text{as}} \\
\mathbf{0} & \mathbf{0} & \mathbf{0} & \mathbf{0} & \boldsymbol{\beta}^{G}_{\text{eu,eu}} & \boldsymbol{\beta}^{G}_{\text{eu,as}} & \boldsymbol{\beta}^{G}_{\text{eu,am}} & \mathbf{0} & \mathbf{0} & \mathbf{0} & \mathbf{0} & \mathbf{0} & \boldsymbol{\beta}_{\text{eu}} & \mathbf{0} \\
\mathbf{0} & \mathbf{0} & \mathbf{0} & \boldsymbol{\beta}^{G}_{\text{am,am}} & \boldsymbol{\beta}^{G}_{\text{am,eu}} & \boldsymbol{\beta}^{G}_{\text{am,as}} & \mathbf{0} & \mathbf{0} & \mathbf{0} & \mathbf{0} & \mathbf{0} & \boldsymbol{\beta}_{\text{am}} & \mathbf{0} & \mathbf{0} \\
\mathbf{0} & \mathbf{0} & \boldsymbol{\beta}^{G}_{\text{as,as}} & \boldsymbol{\beta}^{G}_{\text{as,am}} & \boldsymbol{\beta}^{G}_{\text{as,eu}} & \mathbf{0} & \mathbf{0} & \mathbf{0} & \mathbf{0} & \mathbf{0} & \boldsymbol{\beta}_{\text{as}} & \mathbf{0} & \mathbf{0} & \mathbf{0} \\
\mathbf{0} & \boldsymbol{\beta}^{G}_{\text{eu,eu}} & \boldsymbol{\beta}^{G}_{\text{eu,as}} & \boldsymbol{\beta}^{G}_{\text{eu,am}} & \mathbf{0} & \mathbf{0} & \mathbf{0} & \mathbf{0} & \mathbf{0} & \boldsymbol{\beta}_{\text{eu}} & \mathbf{0} & \mathbf{0} & \mathbf{0} & \mathbf{0} \\
\boldsymbol{\beta}^{G}_{\text{am,am}} & \boldsymbol{\beta}^{G}_{\text{am,eu}} & \boldsymbol{\beta}^{G}_{\text{am,as}} & \mathbf{0} & \mathbf{0} & \mathbf{0} & \mathbf{0} & \mathbf{0} & \boldsymbol{\beta}_{\text{am}} & \mathbf{0} & \mathbf{0} & \mathbf{0} & \mathbf{0} & \mathbf{0}%
\end{array}
\right]
$}\label{r5}
\end{align}
where $\boldsymbol{\beta}^{G}_{\text{am,as}}$ is an abbreviation for $\boldsymbol{\beta}_{\text{am,as\_time}}^{G}$, and so forth. We show in Section \ref{sec one day representation} of the OS that a one-day representation is not identifiable under the case $\phi=0$, when we ignore auto-correlation among the stacked factors when setting up the likelihood.

The loading matrix $\Lambda$ has six row blocks of dimension $N\times14$. Let $\boldsymbol{\lambda}_{k,j}^{\intercal}$ denote the $j$th row of the $k$th row block of $\Lambda$, for $k=1,\ldots, 6$ and $j=1,\ldots,N$. In other words, $\boldsymbol{\lambda}_{1,j}^{\intercal}$ is the factor loading for the $j$th Asian stock in ``day one", while $\boldsymbol{\lambda}_{5,j}^{\intercal}$ is the factor loading for the $j$th European stock in ``day two". In the theoretical sections, we will use $\{\boldsymbol{\lambda}_{k,j}\}$ to denote factor loadings without referring to notation like $\boldsymbol{\beta}_{\text{eu,as}}^{G}$ or $\boldsymbol{\beta}_{\text{as}}$. Each $\boldsymbol{\lambda}_{k,j}$ contains ten zero elements and four non-zero elements. 

With the two-day representation laid out, we first state our assumptions on the idiosyncratic components $\{\boldsymbol{e}_t\}$.\footnote{In this paper we shall state our assumptions based on the two-day representation. It is possible to state similar assumptions based on (\ref{main_r6}), but the notation will become more cluttered.}

%\bigskip

\begin{assumption}
\label{assu model}
\begin{enumerate}[(i)]

\item Suppose that $\{\boldsymbol{e}_{t}\}_{t=1}^{\infty}$ are strictly stationary, with $\mathbb{E}[\boldsymbol{e}_{t}]=\boldsymbol{0}$ and $\Xi:=\mathbb{E}[\boldsymbol{e}_{t}\boldsymbol{e}_{t}^{\intercal}]$ being positive definite. Let $\Sigma_{ee}$ be the diagonal matrix containing the diagonal elements of $\Xi$:
\begin{align}
&  \Sigma_{ee}:=\diag(\mathbb{E}[\boldsymbol{e}_{t}\boldsymbol{e}_{t}^{\intercal}])=\label{eq:sigma}\\
%&  =\diag(\sigma_{as,1}^{2},\ldots,\sigma_{as,N}^{2},\sigma_{eu,1}^{2},\ldots,\sigma_{eu,N}^{2},\sigma_{am,1}^{2},\ldots,\sigma_{am,N}^{2},\sigma_{as,1}^{2},\ldots,\sigma_{as,N}^{2},\sigma_{eu,1}^{2},\ldots,\sigma_{eu,N}^{2},\sigma_{am,1}^{2},\ldots,\sigma_{am,N}^{2})\nonumber\\
&  =:\diag(\sigma_{1,1}^{2},\ldots,\sigma_{1,N}^{2},\sigma_{2,1}^{2},\ldots,\sigma_{2,N}^{2},\sigma_{3,1}^{2},\ldots,\sigma_{3,N}^{2},\sigma_{4,1}^{2},\ldots,\sigma_{4,N}^{2},\sigma_{5,1}^{2},\ldots,\sigma_{5,N}^{2},\sigma_{6,1}^{2},\ldots,\sigma_{6,N}^{2}),\nonumber
\end{align}
where $\sigma_{1,j}^{2}=\sigma_{4,j}^2$, $\sigma_{2,j}^{2}=\sigma_{5,j}^2$ and $\sigma_{3,j}^{2}=\sigma_{6,j}^2$ for $j=1,\ldots, N$. Assume $C^{-1}\leq \sigma_{k,j}^{2}\leq C$ for $k=1,2,3$ and $j=1,\ldots, N$.

\item Let $\Xi_{i,j}$ denote the $(i,j)$th term of $\Xi$ for $i,j=1,\ldots, 6N$. Assume $|\Xi_{i,j}|\leq \xi_{i,j}$ for some $\xi_{i,j}>0$. In addition, $\sum_{i=1}^{N}\xi_{i,j}\leq C$ for any $j$.

\item Define $\Xi_{(t,t')}:=\mathbb{E}[\boldsymbol{e}_{t}\boldsymbol{e}_{t'}^{\intercal}]$. Let $\Xi_{(t,t'),i,i}$ denote the $(i,i)$th term of $\Xi_{(t,t')}$ for $i=1,\ldots, 6N$. Suppose that $|\Xi_{(t,t'),i,i}|\leq \rho_{t,t'}$ for some $\rho_{t,t'}>0$ and for all $i=1,\ldots, 6N$. In addition, $T^{-1}\sum_{t=1}^{T}\sum_{t'=1}^{T}\rho_{t,t'}\leq C$.
\end{enumerate} 
\end{assumption}

%\bigskip

Assumption \ref{assu model} is analogous to Assumption C of \cite{baili2016}. Assumption \ref{assu model}(i) assumes that the strictly stationary idiosyncratic components $\{\boldsymbol{e}_t\}$ have mean zero and positive definite covariance matrix. Moreover, we restrict the idiosyncratic variances to be in a compact set. Assumption \ref{assu model}(ii) controls the cross-sectional correlation of the idiosyncratic components across the three continents and the two calendar days contained in a time unit. Assumption \ref{assu model}(iii) controls the serial correlation of $\{e_{t,k,j}\}$, where $e_{t,k,j}$ is the $[(k-1)N+j]$th element of $\boldsymbol{e}_t$. 

%\bigskip

\begin{assumption}
\label{assu strong mixing and exponential tail}
\begin{enumerate}[(i)]
%For all $s=1,\ldots, 2T$, $c=asia, euro, amer$, $i=1,\ldots, N$, there exist absolute constants $K_{1}>1, K_{2}>0, 2\geq r_{1}>0$ such that%\footnote{"Absolute constants" mean constants that are independent of both $n$ and $T$.}
%\[ \mathbb{E}\left[ \exp\left( K_2 |\epsilon_{c,s,i}|^{r_1}\right) \right] \leq K_{1}.\]

\item For all $t=1,2,\ldots$, $k=1,\ldots, 6$, $j=1,\ldots, N$, there exist absolute constants $K_{1}>1, K_{2}>0,  r_{1}>0$ such that
\begin{align}
\mathbb{E}\left[ \exp\left( K_2 |e_{t,k,j}|^{r_1}\right) \right] \leq K_{1}.\label{main_r21}
\end{align}

%\item For all $c=asia, euro, amer$, $i=1,\ldots, N$, There exist absolute constants $K_{3}>0$ and $r_{2}>0$ such that for all $h\in\mathbb{N}$
%\[\alpha_{c,i}(h)\leq\exp(-K_3h^{r_2}),\]
%where $\alpha_{c,i}(h)$ is the strong mixing coefficients of $\{\epsilon_{c,s,i}\}_{s=1}^{2T}$.

\item For all $k=1,\ldots,6$, $j=1,\ldots, N$, there exist absolute constants $K_{3}>0$ and $r_{2}>0$ such that for all $h\in\mathbb{N}$
\begin{align}
\label{main_r22}
\alpha_{k,j}(h)\leq\exp(-K_3h^{r_2}),
\end{align}
where $\alpha_{k,j}(h)$ is the strong mixing coefficients of $\{e_{t,k,j}\}_{t=1}^{\infty}$.
\end{enumerate}
\end{assumption}

Assumption \ref{assu strong mixing and exponential tail} is quite standard in the literature of the large factor models (\cite{fanliaomincheva2011},  \cite{barigozziluciani2022} etc.). In essence, Assumption \ref{assu strong mixing and exponential tail}(i) assumes that the idiosyncratic components have exponential-type tail probabilities which, together with Assumption \ref{assu strong mixing and exponential tail}(ii), allows us to invoke a version of the Bernstein's inequality for strong mixing time series. Although an exponential-type tail probability implies the existence of all moments, these moments can be quite large (see Remark 2 of \cite{barigozziluciani2022}). The parameter $r_{1}$ restricts size of tails of the idiosyncratic components - the smaller $r_{1}$, the heavier the tails. %When $r_{1}=2$, $e_{t,k,j}$ is said to be \textit{subgaussian}, when $r_{1}=1$, $e_{t,k,j}$ is said to be \textit{subexponential}, and when $0<r_{1}<1$, $e_{t,k,j}$ is said to be \textit{semiexponential}. 

Assumption \ref{assu strong mixing and exponential tail}(ii) assumes that the idiosyncratic components are strong mixing because $\alpha_{k,j}(h)\to 0$ as $h\to\infty$ for all $k,j$. In fact, we require them to decrease at an exponential rate. The bigger $r_{2}$ gets, the faster the decay rate and the less dependence the idiosyncratic components exhibit. This assumption covers a wide range of time series. It is well known that both classical ARMA and GARCH processes are strong mixing with mixing coefficients which decrease to zero at an exponential rate (see Section 2.6.1 of \cite{fanyao2003} and the references therein).

In Section \ref{sec Sufficient conditions of mainr8} of the OS, we prove that Assumption \ref{assu strong mixing and exponential tail} implies that
\begin{align}
\left| \frac{1}{\sqrt{T}}\sum_{t=1}^{T}\left[e_{t,k,j}e_{t,k',j'}-\mathbb{E}[e_{t,k,j}e_{t,k',j'}] \right] \right|^4=O_p(1), \label{main_r8}
\end{align}
pointwise for $k,k'=1,\ldots,6$, $j,j'=1,\ldots, N$. Note that (\ref{main_r8}) is not the same as Assumption C.5 of \cite{baili2016}, which requires the left hand side of (\ref{main_r8}) being bounded in expectation, but (\ref{main_r8}) can achieve the same purposes in the proofs.

We next state the assumptions on the innovations of the global factor and on the continental factors. To facilitate statement of the assumptions, we re-define the components of $\boldsymbol{f}_t$:
\begin{align}
%\boldsymbol{f}_{t}:=\left[
%\begin{array}[c]{c}%
%f^{G}{}_{2t,\text{am\_time}}\\
%f^{G}{}_{2t,\text{eu\_time}}\\
%f^{G}{}_{2t,\text{as\_time}}\\
%f^{G}{}_{2(t-1)+1,\text{am\_time}}\\
%f^{G}{}_{2(t-1)+1,\text{eu\_time}}\\
%f^{G}{}_{2(t-1)+1,\text{as\_time}}\\
%f^{G}{}_{2(t-1),\text{am\_time}}\\
%f^{G}{}_{2(t-1),\text{eu\_time}}\\
%f_{\text{am},2t}\\
%f_{\text{eu},2t}\\
%f_{\text{as},2t}\\
%f_{\text{am},2(t-1)+1}\\
%f_{\text{eu},2(t-1)+1}\\
%f_{\text{as},2(t-1)+1}%
%\end{array}
%\right]=:
\boldsymbol{f}_t&=\left[
\begin{array}{cccccccccccccc}
f_{g,\ell+5}, f_{g,\ell+4}, f_{g,\ell+3}, f_{g,\ell+2}, f_{g,\ell+1}, f_{g,\ell}, f_{g,\ell-1}, f_{g,\ell-2}, f_{c,\ell+5}, f_{c,\ell+4}, f_{c,\ell+3}, f_{c,\ell+2}, f_{c,\ell+1}, f_{c,\ell}%
\end{array}
\right]^{\intercal},\notag\\
 f_{g,\ell} &  =\phi f_{g,\ell-1}+\eta_{g,\ell-1},\label{res_r91}
\end{align} 
for $\ell:=6(t-1)+1$, $t=1,\ldots, T$.

%\bigskip

\begin{assumption}
\label{assu eta}
\begin{enumerate}[(i)]
\item $|\phi|<1$.

\item Suppose that $\{\eta_{g,\ell}\}_{\ell=0}^{\infty}$ are strictly stationary with $\mathbb{E}[\eta_{g,\ell}]=0$, $\mathbb{E}[\eta_{g,\ell}^2]=1$, and satisfy the assumption of exponential-type tail probabilities as in (\ref{main_r21}). %$\eta_{g,\ell}\stackrel{i.i.d.}{\sim} (0,1)$, and
%\begin{align}
%\label{main_r13}
%\int_{\mathbb{R}} |P_{\eta}(x+z)-P_{\eta}(x)|dx\leq C|z|,
%\end{align}
%for any $z\in \mathbb{R}$, where $P_{\eta}(\cdot)$ is the probability density function of $\{\eta_{g,\ell}\}$. In addition, 
%there exist absolute constants $K_{1}>1, K_{2}>0, 2\geq r_{1}>0$ such that 
%\begin{align*}
%\mathbb{E}\left[ \exp\left( K_2 |\eta_{g, \ell}|^{r_1}\right) \right] \leq K_{1},
%\end{align*}
%where $K_1, K_2, r_1$ are the same as those in Assumption \ref{assu strong mixing and exponential tail} without loss of generality. 
Assume $\eta_{g,\ell}=0$ for $\ell <0$.
%where $\{\eta_{g,\ell}\}_{\ell=1}^{6T-5}$ are the innovations of the global factors (see (\ref{res_r91})). 

\item 
%For all $h\in\mathbb{N}$,
%\[\alpha_{f_g}(h)\leq\exp(-K_3h^{r_2}),\]
%where $\alpha_{f_g}(h)$ is the strong mixing coefficient of $\{f_{g,\ell}\}_{\ell=1}^{6T}$, and $K_3, r_2$ are the same as those in Assumption \ref{assu strong mixing and exponential tail} without loss of generality. 
$\{f_{g,\ell}\}_{\ell=1}^{\infty}$ are a strong mixing process with mixing coefficients satisfying (\ref{main_r22}). Assume $f_{g,\ell}=0$ for $\ell\leq 0$.

\item For all $\ell=1,2,\ldots$, $f_{c,\ell}\stackrel{i.i.d.}{\sim} (0,1)$, and satisfy the assumption of exponential-type tail probabilities as in (\ref{main_r21}). %there exist absolute constants $K_{1}>1, K_{2}>0, 2\geq r_{1}>0$ such that 
%\[ \mathbb{E}\left[ \exp\left( K_2 |f_{c, \ell}|^{r_1}\right) \right] \leq K_{1},\]
%where $K_1, K_2, r_1$ are the same as those in Assumption \ref{assu strong mixing and exponential tail} without loss of generality. 
Assume $f_{c,\ell}=0$ for $\ell\leq 0$. In addition, the continental factors are independent of the innovations of the global factor. The innovations of the global factor and the continental factors are independent of the idiosyncratic components $\{\boldsymbol{e}_t\}$.
\end{enumerate}
\end{assumption}

Assumption \ref{assu eta}(i) says that $|\phi|<1$. If the weak-form efficient market hypothesis (along with a time invariant risk premium) holds for stock markets across different time-zones, then one would expect $\phi=0$. On the other hand, if certain international news takes more than a 1/3 day to fully unfold, it is possible that $\phi\neq0$. Assumption \ref{assu eta}(ii) assumes that the innovations of the global factor $\{\eta_{g,\ell}\}_{\ell=0}^{\infty}$ are strictly stationary, and have mean zero, variance one, and some exponential-type tail probability. %Assumption \ref{assu eta}(ii) assumes that the innovations of the global factor $\{\eta_{g,\ell}\}$ are i.i.d. across $\ell$, with mean zero and variance one for normalisation, and have some exponential-type tail probability. The independence assumption of $\{\eta_{g,\ell}\}$ and the integral Lipschitz condition (\ref{main_r13}) ensure that $\{f_{g,\ell}\}$ are a strong mixing process with an exponentially decaying mixing coefficient (\cite{phamtran1985}). 
Assumption \ref{assu eta}(iii) says that the global factors $\{f_{g,\ell}\}$ are a strong mixing process with an exponentially decaying mixing coefficient. Assumption \ref{assu eta}(iii) could be implied by an independence assumption of $\{\eta_{g,\ell}\}$ and an integral Lipschitz condition (see \cite{phamtran1985}). Alternatively, Assumption \ref{assu eta}(iii) could be implied by an assumption of geometric ergodicity of $\{\eta_{g,\ell}\}$, which would allow for GARCH effects (\cite{francq2006mixing}, \cite{barigozziluciani2022}, \cite{barigozzi2022PC}).

Assumption \ref{assu eta}(iv) says that the continental factors $\{f_{c,\ell}\}$ have mean zero, variance one for normalisation, and some exponential-type tail probability. It next assumes that the global and continental factors are independent so that a continental factor contains the information local to the continent itself. In this way, the factors have clear economic interpretations -- so the global and continental factors are truly the ``global" and ``continental" factors, respectively. Indeed, empiricists often orthogonalize global and regional factors (\cite{bekaert2014global}). Assumption \ref{assu eta}(iv) also maintains that the global and continental factors are independent of the idiosyncratic components $\{\boldsymbol{e}_{t}\}$. This is a standard assumption in the literature of the dynamic factor models.

We have two reasons for modelling the continental factors being independent across continents and over time as in Assumption \ref{assu eta}(iv). First, since a continental factor measures the 24-hour continental information, orthogonal to the international news, we believe that investors can efficiently process such information at daily frequency, resulting in very weak auto-correlation of a continental factor. In contrast, the international news represented by the global factor, which measured at 1/3-day frequency, is more plausible to exhibit auto-correlation. Second, the focus of this paper is to study interdependence of stock markets across different time zones, so we abstain from putting unnecessary structures on the continental factors to complicate estimation or the large-sample analysis. Specifying an autoregressive structure for the continental factors may lead to less robust parameter estimates because of the possible convergence to a mediocre local maximum when estimating the model.

%The covariance matrix of $\boldsymbol{f}_{t}$ is denoted by 
Define $M:=\mathbb{E}[\boldsymbol{f}_{t}\boldsymbol{f}_{t}^{\intercal}]$, which is parametrized by $\phi$:
\begin{equation}
\underbrace{M}_{14\times14}=\left[
\begin{array}
[c]{cc}%
\Phi & 0\\
0 & I_{6}%
\end{array}
\right]  ,\qquad\Phi:=\frac{1}{1-\phi^{2}}\left[
\begin{array}
[c]{cccccccc}%
1 & \phi & \phi^{2} & \phi^{3} & \phi^{4} & \phi^{5} & \phi^{6} & \phi^{7}\\
\phi & 1 & \phi & \phi^{2} & \phi^{3} & \phi^{4} & \phi^{5} & \phi^{6}\\
\phi^{2} & \phi & 1 & \phi & \phi^{2} & \phi^{3} & \phi^{4} & \phi^{5}\\
\phi^{3} & \phi^{2} & \phi & 1 & \phi & \phi^{2} & \phi^{3} & \phi^{4}\\
\phi^{4} & \phi^{3} & \phi^{2} & \phi & 1 & \phi & \phi^{2} & \phi^{3}\\
\phi^{5} & \phi^{4} & \phi^{3} & \phi^{2} & \phi & 1 & \phi & \phi^{2}\\
\phi^{6} & \phi^{5} & \phi^{4} & \phi^{3} & \phi^{2} & \phi & 1 & \phi\\
\phi^{7} & \phi^{6} & \phi^{5} & \phi^{4} & \phi^{3} & \phi^{2} & \phi & 1
\end{array}
\right]. \label{eq:m}%
\end{equation}
The specific structures present in the loading matrix $\Lambda$ and the covariance matrix of factors $M$ distinguish our model from those considered by \cite{baili2016}.

\section{Estimation}
\label{sec Estimation}

%We propose two estimators of the parameters of the model: QML-just-identified and QML-all-res. In Section \ref{sec loss function}, we introduce the objective function for estimation. Sections \ref{sec esti QMLE just identified} and \ref{sec esti QMLE all res} provide detailed expositions of these two estimators.

\subsection{The QML Objective Function}
\label{sec loss function}

Consider the two-day representation (\ref{r64}). Define $S_{yy}:=\frac{1}{T}\sum_{t=1}^{T}\boldsymbol{y}_{t}\boldsymbol{y}_{t}^{\intercal}$ and $\Sigma_{yy}:=\Lambda M\Lambda^{\intercal}+\Sigma_{ee}$. Treating $\{\boldsymbol{f}_{t}, \boldsymbol{e}_t\}_{t=1}^{T}$ as i.i.d over $t$, we can write down the quasi-log-likelihood of $\{ \boldsymbol{y}_{t}\}_{t=1}^{T}$ (scaled by $1/(NT)$):
\begin{equation}
\frac{1}{NT}\ell\left(  \{ {\boldsymbol{y}}_{t}\}_{t=1}^{T};\boldsymbol{\theta}\right) =-3\log(2\pi)-\frac{1}{2N}\log|\Sigma_{yy}|-\frac{1}{2N}\tr(S_{yy}\Sigma_{yy}^{-1}), \label{qmle_loss}
\end{equation}
where $\boldsymbol{\theta}:=\{\Lambda, M, \Sigma_{ee}\}$. We call this the \textit{quasi}-log-likelihood, and the estimator that maximizes (\ref{qmle_loss}) the QML estimator, rather than the maximum likelihood estimator (MLE), because of three reasons. First, $\{\boldsymbol{y}_t\}$ are not Gaussian. Second, we use $\Sigma_{ee}$ instead of $\Xi$; that is, we ignore the cross-sectional correlation of $\{\boldsymbol{e}_t\}$. When $N$ is large, treating an approximate factor model as an exact factor model when setting up the likelihood, albeit incorrectly, will not destroy consistency or the asymptotic normality of the QML estimator (\cite{dozgiannonereichlin2012}, \cite{baili2016}, \cite{barigozzi2023}). Third, $\{\boldsymbol{f}_{t},\boldsymbol{e}_t\}_{t=1}^{T}$ are treated as i.i.d. over $t$ when setting up the likelihood. That is, we ignore the auto-correlation between $\boldsymbol{f}_{t}$ and $\boldsymbol{f}_{t-1}$, as well as that between $\boldsymbol{e}_{t}$ and $\boldsymbol{e}_{t-1}$.
%\footnote{\cite{baili2012} also ignored the the auto-correlation between $\boldsymbol{f}_{t}$ and $\boldsymbol{f}_{t-1}$ by treating factors as fixed constants.}
This is the idea of working independence (\cite{panconnett2002}), and is essential for setting up a simple likelihood function of the two-day representation and obtaining an analytic formula for the estimated factors. We stress that $\phi$ enters into the quasi-log-likelihood only through $M$ because $\{\boldsymbol{f}_{t}\}_{t=1}^{T}$ are treated as i.i.d. over $t$ when setting up the likelihood. Correctly setting up the likelihood would result in the MLE. For instance, the MLE could be computed using the \textit{prediction error decomposition} form of the likelihood when $N$ is small (see Chapter 3.4 of  \cite{harvey1990forecasting}). However, when $N$ is large, estimation and inference of the MLE are complicated issues.

\subsection{QML-Just-Identified}

\label{sec esti QMLE just identified}

Utilizing the information that $M$ is symmetric, positive definite and that $\Sigma_{ee}$ is diagonal, we show in Section \ref{sec FOC} of the OS that the resulting first-order conditions (FOCs) of (\ref{qmle_loss}) could not separately identify $\Lambda$ or $M$. We hence need to impose $14^{2}$ identification restrictions on the estimates of $\Lambda$ and $M$ to rule out the rotational indeterminacy. After imposing $14^{2}$ identification restrictions and a set of column-sign restrictions (Assumption \ref{assu ident1}), we obtain the \textit{QML-just-identified} estimators $\hat{\Lambda},\hat{M}$, and $\hat{\Sigma}_{ee}$
which are the solution of the following FOCs:\footnote{See Section
\ref{sec FOC} of the OS for details.}
\begin{align}
\label{r2}\hat{\Lambda}^{\intercal}\hat{\Sigma}_{yy}^{-1}(S_{yy}-\hat{\Sigma
}_{yy})  & =0,\nonumber\\
\diag(\hat{\Sigma}_{yy}^{-1})  &  =\diag(\hat{\Sigma}_{yy}^{-1}S_{yy}%
\hat{\Sigma}_{yy}^{-1}),
\end{align}
where $\hat{\Sigma}_{yy}:=\hat{\Lambda}\hat{M}\hat{\Lambda}^{\intercal}%
+\hat{\Sigma}_{ee}$. These FOCs are identical to (2.7) and (2.8) in \cite{baili2012}, which serve as the cornerstone for establishing the large sample theories. Although $\Lambda$, $M$, $\Sigma_{ee}$ imply more than $14^{2}$
restrictions on $\hat{\Lambda},\hat{M}, \hat{\Sigma}_{ee}$, imposing more than
$14^{2}$ restrictions might result in no solution for (\ref{r2}).

How to select these $14^{2}$ restrictions from those implied by the model is crucial because we cannot impose a restriction that is not instrumental for identification of the two-day representation or establishment of the large-sample theories later. Let $\hat{\boldsymbol{\lambda}}_{k,j}^{\intercal}$ denote the loading of the $j$th asset in the $k$th row block of $\hat{\Lambda}$. Display (\ref{main_r1}) lists the rows of $\hat{\Lambda}$ that have restrictions. Positions marked with ``*'' have no restrictions; positions marked with ``0'' are restricted to be 0; positions marked with ``='' are restricted to be equal to their repetitions in the two-day representation (see (\ref{r5})) while positions marked with ``--'' are these repetitions.
%Note that $\hat{\Lambda}$ has six row blocks, each with $N$ rows and 14 columns, and $\hat{\boldsymbol{\lambda}}_{k,i}^{\intercal}$ represents the loading of ith asset in the $k$th row block. 
Rows of $\hat{\Lambda}$ without any restriction are not listed. In total, we impose 177 restrictions on $\hat{\Lambda}$ (i.e., 156 zero and 21 equality restrictions).
\begin{align}
\left[\begin{array}[c]{c}
\hat{\boldsymbol{\lambda}}_{1,1}^{\intercal}\\
\hat{\boldsymbol{\lambda}}_{1,2}^{\intercal}\\
\hat{\boldsymbol{\lambda}}_{1,3}^{\intercal}\\
\hat{\boldsymbol{\lambda}}_{1,4}^{\intercal}\\
\hat{\boldsymbol{\lambda}}_{2,1}^{\intercal}\\
\hat{\boldsymbol{\lambda}}_{2,2}^{\intercal}\\
\hat{\boldsymbol{\lambda}}_{3,1}^{\intercal}\\
\hat{\boldsymbol{\lambda}}_{3,2}^{\intercal}\\
\hat{\boldsymbol{\lambda}}_{3,3}^{\intercal}\\
\hat{\boldsymbol{\lambda}}_{3,4}^{\intercal}\\
\hat{\boldsymbol{\lambda}}_{4,1}^{\intercal}\\
\hat{\boldsymbol{\lambda}}_{4,2}^{\intercal}\\
\hat{\boldsymbol{\lambda}}_{4,3}^{\intercal}\\
\hat{\boldsymbol{\lambda}}_{4,4}^{\intercal}\\
\hat{\boldsymbol{\lambda}}_{5,1}^{\intercal}\\
\hat{\boldsymbol{\lambda}}_{5,2}^{\intercal}\\
\hat{\boldsymbol{\lambda}}_{6,1}^{\intercal}\\
\hat{\boldsymbol{\lambda}}_{6,2}^{\intercal}\\
\hat{\boldsymbol{\lambda}}_{6,3}^{\intercal}\\
\hat{\boldsymbol{\lambda}}_{6,4}^{\intercal}
\end{array}
\right]  =\left[
\begin{array}[c]{cccccccccccccc}%
0 & 0 & 0 & 0 & 0 & = & = & = & 0 & 0 & 0 & 0 & 0 & -\\
0 & 0 & 0 & 0 & 0 & = & = & = & 0 & 0 & 0 & 0 & 0 & -\\
0 & 0 & 0 & 0 & 0 & = & = & * & 0 & 0 & 0 & 0 & 0 & *\\
0 & 0 & 0 & 0 & 0 & * & * & * & 0 & 0 & 0 & 0 & 0 & *\\
0 & 0 & 0 & 0 & * & * & * & 0 & 0 & 0 & 0 & 0 & * & 0\\
* & * & * & 0 & * & * & * & 0 & * & * & 0 & * & * & 0\\
0 & 0 & 0 & - & - & - & 0 & 0 & 0 & 0 & 0 & - & 0 & 0\\
* & * & 0 & - & - & - & 0 & 0 & * & * & * & - & * & 0\\
* & * & 0 & * & - & - & * & * & * & * & * & * & * & *\\
* & * & * & * & * & - & * & * & * & * & * & * & * & *\\
0 & 0 & - & - & - & 0 & 0 & 0 & 0 & 0 & = & 0 & 0 & 0\\
0 & 0 & - & - & - & 0 & 0 & 0 & 0 & 0 & = & 0 & 0 & 0\\
0 & 0 & - & - & * & 0 & 0 & 0 & 0 & 0 & * & 0 & 0 & 0\\
0 & 0 & * & * & * & * & * & * & 0 & * & * & * & * & *\\
0 & * & * & * & 0 & 0 & 0 & 0 & 0 & * & 0 & 0 & 0 & 0\\
0 & * & * & * & 0 & * & * & * & 0 & * & 0 & * & * & *\\
= & = & = & 0 & 0 & 0 & 0 & 0 & = & 0 & 0 & 0 & 0 & 0\\
= & = & = & 0 & 0 & 0 & 0 & 0 & = & 0 & 0 & 0 & 0 & 0\\
* & = & = & 0 & 0 & 0 & 0 & 0 & * & 0 & 0 & 0 & 0 & 0\\
* & * & = & 0 & 0 & 0 & 0 & 0 & * & 0 & 0 & 0 & 0 & 0
\end{array}
\right] . \label{main_r1}
\end{align}
The remaining 19 (=196-177) restrictions imposed on $\hat{M}$ are:
\begin{align}
\hat{M}_{2,10},\hat{M}_{3,10},\hat{M}_{4,10},\hat{M}%
_{3,11},\hat{M}_{4,11},\hat{M}_{5,11}  &  =0,\notag\\
\hat{M}_{4,12},\hat{M}_{5,12},\hat{M}_{6,12},\hat{M}%
_{5,13},\hat{M}_{6,13},\hat{M}_{7,13}  &  =0,\notag\\
\hat{M}_{10,10},\hat{M}_{11,11},\hat{M}_{12,12},\hat
{M}_{13,13}  &  =1, \label{main_r2}\\
\hat{M}_{4,4}=\hat{M}_{5,5}&=\hat{M}_{6,6},\notag\\ \hat{M}_{4,4}-\hat{M}_{6,4}&=1,\notag
\end{align}
where $\hat{M}_{2,10}$ is the $(2,10)$th element of $\hat{M}$, and so forth. We now state the set of column-sign restrictions.

%\begin{align*}
%\hat{M}_{2,10},\hat{M}_{3,10},\hat{M}_{4,10},\hat{M}_{3,11},\hat{M}_{4,11},\hat{M}_{5,11}  &  =0\\
%\hat{M}_{4,12},\hat{M}_{5,12},\hat{M}_{6,12},\hat{M}_{5,13},\hat{M}_{6,13},\hat{M}_{7,13}  &  =0\\
%\hat{M}_{10,10},\hat{M}_{11,11},\hat{M}_{12,12},\hat{M}_{13,13}  &  =1
%\end{align*}
%and%
%\[\hat{M}_{4,4}=\hat{M}_{5,5}=\hat{M}_{6,6},\hat{\text{ }M}_{4,4}-\hat{M}_{6,4}=1.\]

%\bigskip

\begin{assumption}
\label{assu ident1} 
Suppose that the column signs of $\boldsymbol{\beta}^{G}_{\text{as,as}}, \boldsymbol{\beta}^{G}_{\text{eu,eu}}, \boldsymbol{\beta}^{G}_{\text{am,am}}, \boldsymbol{\beta}_{\text{as}}, \boldsymbol{\beta}_{\text{eu}}, \boldsymbol{\beta}_{\text{am}}$ are known.
\end{assumption}

%\bigskip

Assumption \ref{assu ident1} implies that all the column signs of $\Lambda$ are known. \cite{baili2012} have made similar assumption as an implicit part of their identification schemes (IC2, IC3 and IC5) (see \citet[][p.445, p.463]{baili2012}). Also see Assumption 6(c) of \cite{barigozziluciani2022}. Assumption \ref{assu ident1} is also reasonable in finance. Take $\boldsymbol{\beta}_{\text{as,as}}^G$ as an example. For normalisation, we shall assume that most elements of $\boldsymbol{\beta}_{\text{as,as}}^G$ are positive. That is, most Asian stocks react positively to the global factor during the Asian sub-period, $f_{s,\text{as\_time}}^G$. In this sense, the column sign of $\boldsymbol{\beta}_{\text{as,as}}^G$ is known. However, we do not normalise the column sign of $\boldsymbol{\beta}_{\text{eu,as}}^G$ or of $\boldsymbol{\beta}_{\text{am,as}}^G$, allowing for the possibility of most European and American stocks reacting negatively to $f_{s,\text{as\_time}}^G$. Assumption \ref{assu ident1} ensures identification of QML-just-identified for both $\phi\neq 0$ and $\phi=0$.  %It ensures that the model can obtain the necessary normalization conditions for identification in both $\phi$ equals 0 and not equal to 0, similar to the normalization conditions in traditional factor models.  %Likewise, we allow the possibility of negative reactions from Asian and American stock markets during the European time period, as well as negative reactions from Asian and European stock markets during the American time period.
Section \ref{sec compute just} of the OS explains how to compute QML-just-identified in practice.

\subsection{QML-All-Res}
\label{sec esti QMLE all res}

We now improve upon the QML-just-identified estimator and propose a \textit{QML-all-res} estimator that takes all the restrictions implied by $\Lambda$, $M$, and $\Sigma_{ee}$ into account. The idea is that all the restrictions could be incorporated in an EM algorithm. Define $\boldsymbol{Y}:=(\boldsymbol{y}_1,\ldots, \boldsymbol{y}_T)$ and $\boldsymbol{F}:=(\boldsymbol{f}_1,\ldots, \boldsymbol{f}_T)$. Let $\hat{\boldsymbol{\theta}}^{(p)}:=(\hat{\Lambda}^{(p)}, \hat{\Sigma}_{ee}^{(p)}, \hat{M}^{(p)})$ denote the estimate of $\boldsymbol{\theta}$ in the $p$th iteration of the EM algorithm for $p=0,1,2,\ldots$. The quasi-log-likelihood of $\boldsymbol{Y}$ could be decomposed into 
\begin{align}
\ell(  \boldsymbol{Y};\boldsymbol{\theta})=\ell(  \boldsymbol{Y}, \boldsymbol{F};\boldsymbol{\theta})-\ell(\boldsymbol{F}|\boldsymbol{Y};\boldsymbol{\theta}),\label{main_r27}
\end{align}
where $\ell(  \boldsymbol{Y}, \boldsymbol{F};\boldsymbol{\theta})$ is the complete quasi-log-likelihood of the two-day representation. Taking expectation of (\ref{main_r27}) with respect to the conditional density $P( \boldsymbol{F}|\boldsymbol{Y}) $ at $\hat{\boldsymbol{\theta}}^{(p)}$, we have
\begin{align*}
\ell(  \boldsymbol{Y};\boldsymbol{\theta})\equiv\mathbb{E}_{\boldsymbol{F}|\boldsymbol{Y},\hat{\boldsymbol{\theta}}^{(p)}}\left[ \ell(  \boldsymbol{Y};\boldsymbol{\theta})\right] =\mathbb{E}_{\boldsymbol{F}|\boldsymbol{Y},\hat{\boldsymbol{\theta}}^{(p)}}\left[\ell(  \boldsymbol{Y}, \boldsymbol{F};\boldsymbol{\theta})\right] -\mathbb{E}_{\boldsymbol{F}|\boldsymbol{Y},\hat{\boldsymbol{\theta}}^{(p)}}\left[\ell(\boldsymbol{F}|\boldsymbol{Y};\boldsymbol{\theta})\right]. 
\end{align*}
In the $(p+1)$th iteration, the EM algorithm maximises $\mathbb{E}_{\boldsymbol{F}|\boldsymbol{Y},\hat{\boldsymbol{\theta}}^{(p)}}\left[\ell(  \boldsymbol{Y}, \boldsymbol{F};\boldsymbol{\theta})\right]$ subject to all the restrictions implied by $\Lambda$, $M$, and $\Sigma_{ee}$; let $\hat{\boldsymbol{\theta}}^{(p+1)}$ denote this constrained maximum. Because $\mathbb{E}_{\boldsymbol{F}|\boldsymbol{Y},\hat{\boldsymbol{\theta}}^{(p)}}[\ell(\boldsymbol{F}|\boldsymbol{Y};\boldsymbol{\theta})]\leq \mathbb{E}_{\boldsymbol{F}|\boldsymbol{Y},\hat{\boldsymbol{\theta}}^{(p)}}[\ell(\boldsymbol{F}|\boldsymbol{Y};\hat{\boldsymbol{\theta}}^{(p)})]$ for any $\boldsymbol{\theta}$ (see (4) of \cite{wu1983convergence}), we have $\ell(  \boldsymbol{Y};\hat{\boldsymbol{\theta}}^{(p+1)})\geq \ell(  \boldsymbol{Y};\hat{\boldsymbol{\theta}}^{(p)})$ subject to all the restrictions implied by $\Lambda$, $M$, and $\Sigma_{ee}$. In other words, as the number of iterations increases, $\ell(  \boldsymbol{Y};\boldsymbol{\theta})$ monotonically increases subject to the restrictions.
%The complete quasi-log-likelihood of the two-day representation is
%\begin{align}
%&\ell( \{\boldsymbol{y}_t\}_{t=1}^{T}, \{\boldsymbol{f}_{t}\}_{t=1}^{T};\boldsymbol{\theta}) 
%=\ell( \{\boldsymbol{y}_t\}_{t=1}^T|\{\boldsymbol{f}_{t}\}_{t=1}^T;\boldsymbol{\theta})+\ell( \{\boldsymbol {f}_{t}\}_{t=1}^T;\boldsymbol{\theta}) =\sum_{t=1}^{T}\ell( \boldsymbol{y}_t|\boldsymbol{f}_{t};\boldsymbol{\theta})+\sum_{t=1}^{T}\ell( \boldsymbol {f}_{t};\boldsymbol{\theta}),\notag\\
%&\sum_{t=1}^{T}\ell( \boldsymbol{y}_t|\boldsymbol{f}_{t};\boldsymbol{\theta})=-\frac{6NT}{2}\log(2\pi)-\frac{T}{2}\log|\Sigma_{ee}|-\frac{1}{2}\sum_{t=1}^{T}(\boldsymbol{y}_{t}-\Lambda\boldsymbol{f}_{t})^{\intercal}\Sigma_{ee}^{-1}(\boldsymbol{y}_{t}-\Lambda\boldsymbol{f}_{t}),\notag\\
%&\sum_{t=1}^{T}\ell( \boldsymbol {f}_{t};\boldsymbol{\theta}) =-\frac{14T}{2}\log(2\pi)-\frac{T}{2}\log|M|-\frac{1}{2}\sum_{t=1}^{T}\boldsymbol{f}_{t}^{\intercal}M^{-1}\boldsymbol{f}_{t}.\notag
%\end{align}
%
Because we treat $\{\boldsymbol{f}_t,\boldsymbol{e}_t\}_{t=1}^T$ as i.i.d. over $t$ when setting up the likelihood, we have the following simplification:
\begin{align*}
\mathbb{E}_{\boldsymbol{F}|\boldsymbol{Y},\hat{\boldsymbol{\theta}}^{(p)}}\left[\ell(  \boldsymbol{Y}, \boldsymbol{F};\boldsymbol{\theta})\right]=\sum_{t=1}^{T}\mathbb{E}_{\boldsymbol{F}|\boldsymbol{Y},\hat{\boldsymbol{\theta}}^{(p)}}\left[\ell(  \boldsymbol{y}_t, \boldsymbol{f}_t;\boldsymbol{\theta})\right]=\sum_{t=1}^{T}\mathbb{E}_{\boldsymbol{f}_t|\boldsymbol{y}_t,\hat{\boldsymbol{\theta}}^{(p)}}\left[\ell(  \boldsymbol{y}_t, \boldsymbol{f}_t;\boldsymbol{\theta})\right],
\end{align*}
where $\mathbb{E}_{\boldsymbol{f}_t|\boldsymbol{y}_t,\hat{\boldsymbol{\theta}}^{(p)}}[\cdot]$ denotes the expectation with respect to the conditional density
$P\left( \boldsymbol{f}_{t}|\boldsymbol{y}_t\right) $ at $\hat{\boldsymbol{\theta}}^{(p)}$. However, we do not rely on $\{\boldsymbol{f}_t,\boldsymbol{e}_t\}_{t=1}^T$ being i.i.d. over $t$ in the derivation of the large sample theories. The procedure of the EM algorithm is detailed as follows:
\begin{enumerate}[(i)]
\item {Initialization}. Let $\hat{\boldsymbol{\theta}}^{(0)}:=\{\hat{\hat{\Lambda}}, \hat{\Sigma}_{ee}, \hat{M}\}$, where $\hat{\Sigma}_{ee}$ and $\hat{M}$ are the QML-just-identified estimators for $\Sigma_{ee}$ and $M$, respectively, and $\hat{\hat{\Lambda}}$ is an estimator built upon $\hat{\Lambda}$, the QML-just-identified estimator for $\Lambda$. We call $\hat{\hat{\Lambda}}$ the \textit{QML-minus} estimator and give its definition in Section \ref{sec main paper qml minus}. %(see Section \ref{sec def of minus} of the OS for the definition of $\hat{\hat{\Lambda}}$).
  
\item {E-step.}  %Let $\hat{\boldsymbol{\theta}}^{(p)}:=(\hat{\Lambda}^{(p)}, \hat{\Sigma}_{ee}^{(p)}, \hat{M}^{(p)})$ denote the estimate of $\boldsymbol{\theta}$ in the $p$th iteration. 
The E-step in the $(p+1)$th iteration is to compute $\sum_{t=1}^T \mathbb{E}_{\boldsymbol{f}_t|\boldsymbol{y}_t,\hat{\boldsymbol{\theta}}^{(p)}}\left[ \ell( \boldsymbol{y}_t, \boldsymbol{f}_{t};\boldsymbol{\theta})\right]  $.

\item {M-step.} Let $\{\boldsymbol{\lambda}_{k,j,R}\}$ denote the four non-zero elements in $\{\boldsymbol{\lambda}_{k,j}\}$. %, which are the parameters of interest of loadings in the case of QML-all-res. 
Find values of $\{\boldsymbol{\lambda}_{k,j,R}\},\{\sigma_{k,j}^2\},\phi$ to maximize  $\sum_{t=1}^T \mathbb{E}_{\boldsymbol{f}_t|\boldsymbol{y}_t,\hat{\boldsymbol{\theta}}^{(p)}}\left[ \ell
( \boldsymbol{y}_t, \boldsymbol{f}_{t};\boldsymbol{\theta})\right] $. All the restrictions implied by $\Lambda$, $M$, and $\Sigma_{ee}$ will be taken into account. We denote them as $\{\hat{\boldsymbol{\lambda}}_{k,j,R}^{(p+1)}\}$, $\{\hat{\sigma}^{2,(p+1)}_{k,j}\}$, $\hat{\phi}^{(p+1)}$. Form $\hat{\Lambda}^{(p+1)}, \hat{\Sigma}_{ee}^{(p+1)}, \hat{M}^{(p+1)}$ using $\{\hat{\boldsymbol{\lambda}}_{k,j,R}^{(p+1)}\}$, $\{\hat{\sigma}^{2,(p+1)}_{k,j}\}$, $\hat{\phi}^{(p+1)}$ via (\ref{r5}), (\ref{eq:sigma}), and (\ref{eq:m}), respectively. This is the estimate of $\boldsymbol{\theta}$ in the $(p+1)$th iteration, denoted $\hat{\boldsymbol{\theta}}^{(p+1)}$.

%\item \textbf{E-step}.  In the E-step of the $\left(  j+1\right)  $th iteration, we calculate,%
%\begin{align*}
%E_{f|Y,\hat{\theta}^{(j)}}\boldsymbol{f}_{t}  &  =\left[  \left(
%\hat{M}^{(j)}\right)  ^{-1}+\left(  \hat{\Lambda}^{(j)}\right)^{\intercal}\left(  \hat{\Sigma}_{ee}^{(j)}\right)  ^{-1}\left(
%\hat{\Lambda}^{(j)}\right)  \right]  ^{-1}\left(  \hat{\Lambda}^{(j)}\right)  ^{\intercal}\left(  \hat{\Sigma}_{ee}^{(j)}\right)^{-1}\boldsymbol{y}_{t}\\
%E_{f|Y,\hat{\theta}^{(j)}}\boldsymbol{f}_{t}\boldsymbol{f}_{t}%
%{}^{\intercal}  &  =\left[  \left(  \hat{M}^{(j)}\right)  ^{-1}+\left(
%\hat{\Lambda}^{(j)}\right)  ^{\intercal}\left(  \hat{\Sigma}%
%_{ee}^{(j)}\right)  ^{-1}\left(  \hat{\Lambda}^{(j)}\right)  \right]
%^{-1}+\left(  E_{f|Y,\hat{\theta}^{(j)}}f_{t}\right)  \left(
%E_{f|Y,\hat{\theta}^{(j)}}f_{t}\right)  ^{\intercal}%
%\end{align*}

%\item M-step: Obtain the $\hat{\theta}^{(j+1)}$ by maximizing the expected QMLE loss function with respect to $\theta$ while imposing all the restrictions on ${\Lambda}$ and ${M}$. That gives $\hat{\Lambda}^{(j+1)}$, $\hat{M}^{(j+1)}$and $\hat{\Sigma}_{ee}^{(j+1)}.$

\item Iterate steps (ii) and (iii) until convergence. We denote the final iteration estimates as $\{\hat{\boldsymbol{\lambda}}_{k,j,R}^{*}\}, \{\hat{\sigma}^{2,*}_{k,j}\}, \hat{\phi}^{*}$. Form $\hat{\boldsymbol{\theta}}^*:=(\hat{\Lambda}^{*}, \hat{\Sigma}_{ee}^{*}, \hat{M}^{*})$ using $\{\hat{\boldsymbol{\lambda}}_{k,j,R}^{*}\}, \{\hat{\sigma}^{2,*}_{k,j}\}, \hat{\phi}^{*}$ via (\ref{r5}), (\ref{eq:sigma}), and (\ref{eq:m}), respectively, and name them the \textit{QML-all-res} estimators.
\end{enumerate}

%\bigskip

We have some remarks about our EM algorithm. First, in the literature, there are two types of EM algorithms for estimation of factor models. The first type was used by \cite{dozgiannonereichlin2012}, and originally introduced in the 1980s-1990s, say, by \cite{watsonengle1983} and \cite{quahsargent1992} (also see \citet[][p.23]{barigozzi2023} for a brief discussion on this). The second type was used by \cite{baili2012} and originally proposed by \cite{rubinthayer1982}. Both types of EM algorithms use, in terms of our notation, $\mathbb{E}_{\boldsymbol{F}|\boldsymbol{Y}, \hat{\boldsymbol{\theta}}^{(p)}}[\ell( \boldsymbol{Y}, \boldsymbol{F};\boldsymbol{\theta})]$ in the E-step. The difference is that in the second type of EM algorithm because $\{\boldsymbol{f}_t,\boldsymbol{y}_t\}_{t=1}^T$ are independent over $t$, or treated as independent when setting up the likelihood,  $\mathbb{E}_{\boldsymbol{F}|\boldsymbol{Y}, \hat{\boldsymbol{\theta}}^{(p)}}[\ell( \boldsymbol{Y}, \boldsymbol{F};\boldsymbol{\theta})]$ is reduced to $\sum_{t=1}^{T}\mathbb{E}_{\boldsymbol{f}_t|\boldsymbol{y}_t, \hat{\boldsymbol{\theta}}^{(p)}}[\ell( \boldsymbol{y}_t, \boldsymbol{f}_{t};\boldsymbol{\theta})]$ in the E-step, whereas in the first type of EM algorithm, $\{\boldsymbol{y}_t\}_{t=1}^T$ are correlated over $t$. Our EM algorithm is of the second type.

Second, in practice, the number of iterations is finite, so there are two gaps between the output of an EM algorithm and the QML estimator. The first gap is the one between the output of an EM algorithm and the local maximum of the quasi-log-likelihood, to which an EM algorithm usually converges when iterations are infinite under some conditions (see \cite{wu1983convergence} and Section E of the supplement of \cite{baili2012}). The second gap is the one between the local maximum and the QML estimator (i.e., the global maximum). To the best of our knowledge, only \cite{barigozziluciani2022} have proved that these two gaps are asymptotically negligible but under a different EM algorithm and a different identification scheme (\citet[p.19]{barigozziluciani2022}). In this paper, we do not directly consider convergence of QML-all-res to the global maximum of $\ell( \boldsymbol{Y}; \boldsymbol{\theta})$ subject to all the restrictions implied by $\Lambda$, $M$ and $\Sigma_{ee}$. Nevertheless, a sensible starting value allows the local maximum, to which an EM algorithm converges when iterations are infinite, to be close to the global maximum.  Using $\hat{\boldsymbol{\theta}}^{(0)}$ as the starting value for the QML-all-res estimator allows us to establish the large-sample theories of QML-all-res because the rates of convergence of the starting value serve as a stepping stone in the proof. This has some flavor of the one-step estimator (\cite{barigozziluciani2022}).

Third, we provide more details about step (iii) (i.e., the M-step). We derived in Section \ref{sec res intro2} of the OS that
\begin{align}
\hat{\boldsymbol{\lambda}}_{k,j,R}^{(p+1)}\equiv\hat{\boldsymbol{\lambda}}_{k+3,j,R}^{(p+1)}&=\left(  L_k\frac{1}{T}\sum_{t=1}^{T}\mathbb{E}_{\boldsymbol{f}_t|\boldsymbol{y}_t,\hat{\boldsymbol{\theta}}^{(p)}}[\boldsymbol{f}_{t}\boldsymbol{f}_{t}^{\intercal}]L_k^{\intercal}+L_{k+3}\frac{1}{T}\sum_{t=1}^{T}\mathbb{E}_{\boldsymbol{f}_t|\boldsymbol{y}_t,\hat{\boldsymbol{\theta}}^{(p)}}[\boldsymbol{f}_{t}\boldsymbol{f}_{t}^{\intercal}]L_{k+3}^{\intercal}\right)  ^{-1}\notag\\
&\qquad \times \left(  L_k\frac{1}{T}\sum_{t=1}^{T}\mathbb{E}_{\boldsymbol{f}_t|\boldsymbol{y}_t,\hat{\boldsymbol{\theta}}^{(p)}}[\boldsymbol{f}_{t}]y_{t,k,j}+L_{k+3}\frac{1}{T}\sum_{t=1}^{T}\mathbb{E}_{\boldsymbol{f}_t|\boldsymbol{y}_t,\hat{\boldsymbol{\theta}}^{(p)}}[\boldsymbol{f}_{t}]y_{t,k+3,j}\right),\label{main_r12}
\end{align}
for $j=1,\ldots, N$ and $k=1,2,3$, where $y_{t,k,j}$ is the $[(k-1)N+j]$th element of $\boldsymbol{y}_t$, and $L_{k}$ is a $4 \times 14$ selection matrix to pick out the non-zero elements in $\boldsymbol{\lambda}_{k,j}$ (i.e., $\boldsymbol{\lambda}_{k,j,R}=L_k\boldsymbol{\lambda}_{k,j}$, $\boldsymbol{\lambda}_{k,j}=L_k^{\intercal}\boldsymbol{\lambda}_{k,j,R}$, and $L_kL_k^{\intercal}=I_4$). Moreover, we have 
%\begin{align*}
%\hat{\sigma}_{k,j}^{2,(p+1)}  &  =\frac{1}{2T}\sum_{t=1}^{T}\left[y_{t,k,j}^{2}-2\mathbb{E}_{\boldsymbol{f}_t|\boldsymbol{y}_t,\hat{\boldsymbol{\theta}}^{(p)}}[\boldsymbol{f}_{t}^{\intercal}]\hat{\boldsymbol{\lambda}}_{k,j}^{(p+1)}y_{t,k,j}+\hat{\boldsymbol{\lambda}}_{k,j}^{(p+1),\intercal}\mathbb{E}_{\boldsymbol{f}_t|\boldsymbol{y}_t,\hat{\boldsymbol{\theta}}^{(p)}}[\boldsymbol{f}_{t}\boldsymbol{f}_{t}^{\intercal}]\hat{\boldsymbol{\lambda}}_{k,j}^{(p+1)}\right]\\
%& \qquad +\frac{1}{2T}\sum_{t=1}^{T}\left[  y_{t,k+3,j}^{2}-2\mathbb{E}_{\boldsymbol{f}_t|\boldsymbol{y}_t,\hat{\boldsymbol{\theta}}^{(p)}}[\boldsymbol{f}_{t}^{\intercal}]\hat{\boldsymbol{\lambda}}_{k+3,j}^{(p+1)}y_{t,k+3,j}+\hat{\boldsymbol{\lambda}}_{k+3,j}^{(p+1), \intercal}\mathbb{E}_{\boldsymbol{f}_t|\boldsymbol{y}_t,\hat{\boldsymbol{\theta}}^{(p)}}[\boldsymbol{f}_{t}\boldsymbol{f}_{t}^{\intercal}]\hat{\boldsymbol{\lambda}}_{k+3,j}^{(p+1)}\right],
%\end{align*}
\begin{align}
\hat{\sigma}_{k,j}^{2,(p+1)}=\frac{1}{2T}\sum_{t=1}^{T}\mathbb{E}_{\boldsymbol{f}_t|\boldsymbol{y}_t,\hat{\boldsymbol{\theta}}^{(p)}}\left[  ({y}_{t,k,j}-\hat
{\boldsymbol{\lambda}}_{k,j}^{(p+1),\intercal}\boldsymbol{f}_{t})^{2}+({y}_{t,k+3,j}-\hat{\boldsymbol{\lambda}}_{k+3,j}^{(p+1),\intercal}\boldsymbol{f}_{t})^{2}\right]\label{main_r11}
\end{align}
for $j=1,\ldots, N$ and $k=1,2,3$. As for $\hat{\phi}^{(p+1)}$:
\begin{align}
\hat{\phi}^{(p+1)}&=\arg\max_{|\phi|<1} \sum_{t=1}^T \mathbb{E}_{\boldsymbol{f}_t|\boldsymbol{y}_t,\hat{\boldsymbol{\theta}}^{(p)}}\left[ \ell( \boldsymbol{y}_t, \boldsymbol{f}_{t};\boldsymbol{\theta})\right]=\arg\max_{|\phi|<1} \sum_{t=1}^T \mathbb{E}_{\boldsymbol{f}_t|\boldsymbol{y}_t,\hat{\boldsymbol{\theta}}^{(p)}}\left[ \ell
( \boldsymbol{f}_{t};\boldsymbol{\theta})\right] \notag \\
&  =\arg \min_{|\phi|<1}\left[\log|M\left(  \phi \right)
|+\frac{1}{T}\sum_{t=1}^{T}\tr\left(  \mathbb{E}_{\boldsymbol{f}_t|\boldsymbol{y}_t,\hat{\boldsymbol{\theta}}^{(p)}}[\boldsymbol{f}_{t}\boldsymbol{f}_{t}^{\intercal}]M(  \phi)^{-1}\right)\right].\label{main_r10}
\end{align}
The FOC of (\ref{main_r10}) is given by (\ref{mle4_10}) in the OS. %, from which one could so  LVe for $\hat{\phi}^{(p+1)}$. 
Alternatively, one could obtain $\hat{\phi}^{(p+1)}$ with numerical optimization, whose computational burden is almost negligible. We remark that alternative objective functions other than (\ref{main_r10}) could be used to compute $\phi$ in the $(p+1)$th iteration. %In fact, the QML-all-res estimator for $\phi$ could be set to any consistent estimator of $\phi$ in the step of initialization, and in the following iterations we do not update it. 
When evaluating (\ref{main_r12}), (\ref{main_r11}), and (\ref{main_r10}), one needs to compute  $\mathbb{E}_{\boldsymbol{f}_t|\boldsymbol{y}_t,\hat{\boldsymbol{\theta}}^{(p)}}[\boldsymbol{f}_t]$ and $\mathbb{E}_{\boldsymbol{f}_t|\boldsymbol{y}_t,\hat{\boldsymbol{\theta}}^{(p)}}[\boldsymbol{f}_t\boldsymbol{f}_t^{\intercal}]$ of the E-step. We derive in Section \ref{sec res intro} of the OS that under Gaussianity of $\{\boldsymbol{f}_t, \boldsymbol{y}_t\}$:
\begin{align}
\mathbb{E}_{\boldsymbol{f}_t|\boldsymbol{y}_t,\hat{\boldsymbol{\theta}}^{(p)}}[\boldsymbol{f}_t]&=[\hat{M}^{(p),-1}+\hat{\Lambda}^{(p),\intercal}\hat{\Sigma}_{ee}^{(p),-1}\hat{\Lambda}^{(p)}]^{-1}\hat{\Lambda}^{(p),\intercal}\hat{\Sigma}_{ee}^{(p),-1}\boldsymbol{y}_t,\label{main_r3}\\
\mathbb{E}_{\boldsymbol{f}_t|\boldsymbol{y}_t,\hat{\boldsymbol{\theta}}^{(p)}}[\boldsymbol{f}_t\boldsymbol{f}_t^{\intercal}]
&=(M^{(p),-1}+\Lambda^{(p),\intercal}\Sigma_{ee}^{(p),-1}\Lambda^{(p)})^{-1}+\mathbb{E}_{\boldsymbol{f}_t|\boldsymbol{y}_t,\hat{\boldsymbol{\theta}}^{(p)}}[\boldsymbol{f}_t]\mathbb{E}_{\boldsymbol{f}_t|\boldsymbol{y}_t,\hat{\boldsymbol{\theta}}^{(p)}}[\boldsymbol{f}_t^{\intercal}].\label{main_r9}
\end{align}
%
%The maximization of $\hat{\Lambda}^{(j+1)}$ and $\hat{\Sigma}_{ee}^{(j+1)}$ in the M-step (iii) can be obtained analytically, as follows:
%\begin{align*}
%\hat{\boldsymbol{\lambda}}_{k,i}^{(j+1)}  &  =L_{k}^{\intercal}\left(L_{k}\frac{1}{T}\sum_{t=1}^{T}E_{f|Y,\hat{\theta}^{(j)}}\boldsymbol{f}_{t}\boldsymbol{f}_{t}{}^{\intercal}L_{k}^{\intercal}+L_{k+3}\frac{1}{T}\sum_{t=1}^{T}E_{f|Y,\hat{\theta}^{(j)}}\boldsymbol{f}_{t}\boldsymbol{f}_{t}{}^{\intercal}L_{k+3}^{\intercal}\right)  ^{-1}\\
%&  \qquad \times \left(  L_{k}\frac{1}{T}\sum_{t=1}^{T}\left(  E_{f|Y,\hat{\theta}^{(p)}}\boldsymbol{f}_{t}\right)  y_{t,k,i}+L_{k+3}\frac{1}{T}\sum_{t=1}^{T}\left(  E_{f|Y,\hat{\theta}^{(p)}}\boldsymbol{f}_{t}\right)y_{t,k+3,i}\right)  ,\\
%\end{align*}
%

Without Gaussianity of $\{\boldsymbol{f}_t, \boldsymbol{y}_t\}$, the equalities in (\ref{main_r3}) and (\ref{main_r9}) do not necessarily hold. For instance, the right hand side of (\ref{main_r3}) equals the Kalman smoother of $\boldsymbol{f}_t$ given $\boldsymbol{Y}$ evaluated at the parameter $\hat{\boldsymbol{\theta}}^{(p)}$ or equivalently the linear projection of $\boldsymbol{f}_t$ on $\boldsymbol{Y}$ and a constant. Without Gaussianity, this Kalman smoother/linear projection does not necessarily equal the conditional moment $\mathbb{E}_{\boldsymbol{f}_t|\boldsymbol{y}_t,\hat{\boldsymbol{\theta}}^{(p)}}[\boldsymbol{f}_t]$. To ensure the equalities in (\ref{main_r3}) and (\ref{main_r9}) under non-Gaussianity, we make the following assumption:

\begin{assumption}
\label{assu conditional expection linear}
$\mathbb{E}_{\boldsymbol{f}_t|\boldsymbol{y}_t,\boldsymbol{\theta}}[\boldsymbol{f}_t]$ is linear in $\boldsymbol{y}_t$ for all $\boldsymbol{\theta}$ and $t$.
\end{assumption}

Assumption \ref{assu conditional expection linear} is the same as Assumption 4 of \cite{barigozziluciani2022}, and similar to that made by \citet[p.292]{quahsargent1992}. Gaussianity implies Assumption \ref{assu conditional expection linear}.\footnote{Under non-Gaussianity, although the estimators defined in (\ref{main_r12})-(\ref{main_r10}) with $\mathbb{E}_{\boldsymbol{f}_t|\boldsymbol{y}_t,\hat{\boldsymbol{\theta}}^{(p)}}[\boldsymbol{f}_t]$ and $\mathbb{E}_{\boldsymbol{f}_t|\boldsymbol{y}_t,\hat{\boldsymbol{\theta}}^{(p)}}[\boldsymbol{f}_t\boldsymbol{f}_t^{\intercal}]$ replaced by their respective estimators, that is, the right hand sides of (\ref{main_r3}) and (\ref{main_r9}), remain consistent and asymptotically normal, the resulting algorithm is not an EM algorithm anymore and hence the relationship between its output and the QML estimator remains unclear.}

\subsection{QML-Minus}
\label{sec main paper qml minus}

In this section, we give a definition of the QML-minus estimator $\hat{\hat{\Lambda}}$, which is an auxiliary estimator used in $\hat{\boldsymbol{\theta}}^{(0)}$ in the initialiation of QML-all-res. Moreover, as we shall mention in Section \ref{sec asy rep of qmle-minus} of the OS, the uniform asymptotic representation of QML-minus is a crucial step for establishing the large sample theories for QML-all-res.  %The definition of QML-minus is quite cumbersome, so for readers wishing to get the gist of this paper quickly could skip this section in the first pass of reading.

For loadings, we could build upon the QML-just-identified estimator $\hat{\boldsymbol{\lambda}}_{k,j}$ to define a QML-minus estimator $\hat{\hat{\boldsymbol{\lambda}}}_{k,j}$. In (\ref{r61}) of the OS, we show that the asymptotic representation of the QML-just-identified estimator $\hat{\boldsymbol{\lambda}}_{k,j}$ under $\sqrt{T}/N\to 0$ is:
\begin{align*}
\sqrt{T}(\hat{\boldsymbol{\lambda}}_{k,j}-\boldsymbol{\lambda}_{k,j})=\sqrt{T}(\boldsymbol{\lambda}_{k,j}^{\intercal}\otimes I_{14})\ve(A)+M^{-1}\frac{1}{\sqrt{T}}\sum_{t=1}^{T}\boldsymbol{f}_{t}e_{t,k,j}+o_{p}(1),
\end{align*}
pointwise for $k=1,\ldots, 6$ and $j=1,\ldots, N$, where $A:=(\hat{\Lambda}-\Lambda)^{\intercal}\hat{\Sigma}_{ee}^{-1}\hat{\Lambda}(\hat{\Lambda}^{\intercal}\hat{\Sigma}_{ee}^{-1}\hat{\Lambda})^{-1}$ and $e_{t,k,j}$ is the $[(k-1)N+j]$th element of $\boldsymbol{e}_t$. QML-just-identified is inefficient because the
term $\sqrt{T}(\boldsymbol{\lambda}_{k,j}^{\intercal}\otimes I_{14})\ve(A)$ introduces additional estimation noise compared to the infeasible OLS estimator of the loading which treats $\{\boldsymbol{f}_{t}\}$ as known. If we could eliminate this term, we could obtain an estimator that is as efficient as the infeasible OLS estimator. This is the intuition behind QML-minus.
For $k=1,\ldots, 6$ and $j=1,\ldots, N$, %we define
\begin{align*}
\hat{\hat{\boldsymbol{\lambda}}}_{k,j}:=\hat{\boldsymbol{\lambda}}_{k,j}-( \hat{\boldsymbol{\lambda}}_{k,j}^{\intercal}\otimes I_{14}) (\hat{B}_{3}\hat{\hat{\boldsymbol{a}}}),
\end{align*}
where $\hat{B}_{3}$ and $\hat{\hat{\boldsymbol{a}}}$ are defined in Section \ref{sec def of minus} of the OS. Form $\hat{\hat{\Lambda}}$ using $\{\hat{\hat{\boldsymbol{\lambda}}}_{k,j}\}$ via (\ref{r5}). While the QML-just-identified estimator $\hat{\Lambda}$ only utilizes $14^2$ restrictions implied by $\Lambda$ and $M$, the QML-minus estimator $\hat{\hat{\Lambda}}$ not only utilizes these $14^2$ restrictions, but also utilizes all other restrictions in $\Lambda$ when computing $\hat{\hat{\boldsymbol{a}}}$. In some sense, $\hat{\hat{\Lambda}}$ is the projection of $\hat{\Lambda}$ onto the subspace that embodies all the restrictions in $\Lambda$ and is therefore more efficient.  %, and $\hat{\hat{\Lambda}}$ is used in $\hat{\boldsymbol{\theta}}^{(0)}$ in the initialization of QML-all-res. 

\section{Large Sample Theories}
\label{sec Theories}

We make the following assumption. %Let $C$ be a large constant.

%\bigskip

\begin{assumption}
\label{assu estimated within compact}

\begin{enumerate}[(i)]
\item The factor loadings $\{ \boldsymbol{\lambda}_{k,j}\}$ satisfy $\|\boldsymbol{\lambda}_{k,j}\|\leq C$ for all $k$ and $j$. All elements of $\hat{\boldsymbol{\lambda}}_{k,j}$ are restricted to a compact set $[-C, C]$ for all $k$ and $j$.

\item %Assume $C^{-1}\leq \sigma_{k,j}^{2}\leq C$ for all $k$ and $j$. 
Also $\hat{\sigma}_{k,j}^{2}$ is restricted to a compact set $[C^{-1}, C]$ for all
$k$ and $j$.

\item $\hat{M}$ is restricted to be in a set consisting of all positive definite matrices with all the elements bounded in the interval $[-C, C]$.

\item Suppose that $Q:=\lim_{N\to \infty}\frac{1}{N}\Lambda^{\intercal}%
\Sigma_{ee}^{-1}\Lambda$ is a positive definite matrix.
\end{enumerate}
\end{assumption}

%\bigskip

The first half of Assumption \ref{assu estimated within compact}(i) and Assumption \ref{assu estimated within compact}(iv) are the same as Assumption B of \cite{baili2016}, whereas the second half of Assumption \ref{assu estimated within compact}(i) requires that all the elements of the QML-just-identified estimator for loadings, $\{\hat{\boldsymbol{\lambda}}_{k,j}\}$, are estimated within a compact set. Assumption \ref{assu estimated within compact}(ii) is the same as the first sentence of Assumption D of \cite{baili2016}. %Assumption \ref{assu estimated within compact}(ii) implies that the moment generating functions of $\{e_{t,k,j}\}$ are uniformly (in $t,k,j$) bounded.  
Assumption \ref{assu estimated within compact}(iii) is part of Assumption D of \cite{baili2016} and implies that $\hat{M}=O_p(1)$ and $\hat{M}^{-1}=O_p(1)$. Assumption \ref{assu estimated within compact}(iv) implies that the loadings in display (\ref{main_r6}) are non-sparse.\footnote{In Section \ref{sec pervasive and weak factors} of the OS, we discuss the links between our factors and strong/weak factors.}

%is standard in the literature of factor models and has been taken from the assumptions of \cite{baili2012}.

%Assumption \ref{assu estimated within compact}(ii) and (iv) together imply that the separate identification of $\Lambda \boldsymbol{f}_t$ and $e_t$ when $N\to \infty$ and it is one of the basis for large factor models (\cite{baing2002}, \cite{dozgiannonereichlin2012}, \cite{barigozziluciani2020} etc).

%\bigskip

\begin{assumption}
\label{assu wt}
Define $\boldsymbol{w}_{t}:=\frac{1}{\sqrt{6N}}\sum_{p=1}^{6}\sum_{q=1}^{N}\frac{\boldsymbol{\lambda}_{p,q}}{\sigma_{p,q}^{2}}{e}_{t,p,q}$. Let $w_{t,i}$ denote the $i$th component of $\boldsymbol{w}_{t}$ for $i=1,\ldots, 14$. Suppose that for all $t=1,2,\ldots$, $i=1,\ldots, 14$, $w_{t,i}$ satisfies the assumption of exponential-type tail probabilities as in (\ref{main_r21}).
%\[ \mathbb{E}\left[ \exp\left( K_2 |w_{t,i}|^{r_1}\right) \right] \leq K_{1},\]
%where $K_1, K_2, r_1$ are the same as those in Assumption \ref{assu strong mixing and exponential tail} without loss of generality.
\end{assumption}

Assumption \ref{assu wt} implies that $\boldsymbol{w}_t$ has an exponential-type tail probability. Together with $\{\boldsymbol{w}_t\}$ being a strong mixing process with an exponentially decaying coefficient, it allows us to invoke a version of the Bernstein’s inequality for strong mixing time series. Such an assumption is quite common in the literature (e.g., see Assumption 5(b) of \cite{barigozziluciani2022}).

\subsection{QML-Just-Identified}
\label{sec theory qmle just identified}

We now present the large sample theories of the QML-just-identified estimator. The idea of the proof is based on that  of \cite{baili2016}, but is considerably more involved because our identification scheme is non-standard.

\begin{proposition}
\label{prop5.2} 
Suppose that Assumptions \ref{assu model}%\ref{assu strong mixing and exponential tail}, \ref{assu eta}, \ref{assu ident1} and 
--\ref{assu estimated within compact} hold. When $N,T\to \infty$, with the $196$ particular restrictions imposed on $\hat{\Lambda}$ and $\hat{M}$ as in (\ref{main_r1}) and (\ref{main_r2}), we have
\begin{enumerate}[(i)]
\item \[\Vert \hat{\boldsymbol{\lambda}}_{k,j}-\boldsymbol{\lambda}_{k,j}\Vert^{2}=O_{p}(T^{-1})+O_p(N^{-2}),\]
pointwise for $k=1,\ldots,6,j=1,\ldots,N$.

\item 
\[\frac{1}{6N}\sum_{k=1}^{6}\sum_{j=1}^{N}\frac{1}{\hat{\sigma}_{k,j}^{2}}\Vert \hat{\boldsymbol{\lambda}}_{k,j}-\boldsymbol{\lambda}_{k,j}\Vert^{2}  =O_{p}(T^{-1})+O_p(N^{-2}).\]

\item \begin{equation}
\frac{1}{6N}\sum_{k=1}^{6}\sum_{j=1}^{N}(\hat{\sigma}_{k,j}^{2}-\sigma_{k,j}^{2})^{2}    =O_{p}(T^{-1})+O_p(N^{-2}).\label{5.2b}
\end{equation}

\item \[\hat{M}-M    =O_{p}(T^{-1/2})+O_p(N^{-1}).\]
\end{enumerate}
\end{proposition}

%\bigskip

Proposition \ref{prop5.2}(ii)-(iv) resemble Proposition 1 of \cite{baili2016} and establish some average rates of convergence for $\{\hat{\boldsymbol{\lambda}}_{k,j}\}$ and $\{\hat{\sigma}_{k,j}^2\}$, and the rate of convergence for $\hat{M}$. Note that dimensions of $\Lambda, M, \Sigma_{ee}$ are $6N\times 14$, $14\times 14$, and $6N\times 6N$, respectively. As $N$ increases, the number of parameters to be estimated by QML-just-identified also increases. That is when the notion of average rates of convergence %$\{\hat{\boldsymbol{\lambda}}_{k,j}\}$ and $\{\hat{\sigma}_{k,j}^2\}$ 
comes in handy.  Proposition \ref{prop5.2}(i) establishes the rate of convergence for the individual loading estimator $\hat{\boldsymbol{\lambda}}_{k,j}$. For estimation of each $\boldsymbol{\lambda}_{k,j}$, $k=1,\ldots, 6$, $j=1,\ldots, N$, we do not incur the curse of dimensionality of $N$ because  we can rely on $T\to \infty$.

We now establish the asymptotic distributions of the QML-just-identified estimators. Define
{\footnotesize
\begin{align}
\Upsilon_{k,j}&:=\lim_{T\to \infty}\frac{1}{T}\sum_{t=1}^T\sum_{t'=1}^T\left( \begin{array}{cc}
\mathbb{E}[\boldsymbol{f}_{t}\boldsymbol{f}_{t'}^{\intercal}]\mathbb{E}[e_{t,k,j}e_{t',k,j}]  &  \mathbb{E}[\boldsymbol{f}_{t}\boldsymbol{f}_{t'}^{\intercal}]\mathbb{E}[e_{t,k,j}e_{t',k+3,j}]\\
\mathbb{E}[\boldsymbol{f}_{t}\boldsymbol{f}_{t'}^{\intercal}]\mathbb{E}[e_{t,k+3,j}e_{t',k,j}]     & \mathbb{E}[\boldsymbol{f}_{t}\boldsymbol{f}_{t'}^{\intercal}]\mathbb{E}[e_{t,k+3,j}e_{t',k+3,j}]
\end{array} \right) \label{just_r33}\\
\Upsilon_{e,k,j}&:=\lim_{T\to \infty}\frac{1}{T}\sum_{t=1}^{T}\sum_{t'=1}^{T}\left( \begin{array}{cc}
\mathbb{E}\left[ (e_{t,k,j}^2-\sigma^2_{k,j})(e_{t',k,j}^2-\sigma^2_{k,j})\right] & \mathbb{E}\left[ (e_{t,k,j}^2-\sigma^2_{k,j})(e_{t',k+3,j}^2-\sigma^2_{k+3,j})\right]\\
\mathbb{E}\left[ (e_{t,k+3,j}^2-\sigma^2_{k+3,j})(e_{t',k,j}^2-\sigma^2_{k,j})\right] & \mathbb{E}\left[ (e_{t,k+3,j}^2-\sigma^2_{k+3,j})(e_{t',k+3,j}^2-\sigma^2_{k+3,j})\right]
\end{array}\right), \label{just_r34}
\end{align}}
for $k=1,2,3$ and $j=1,\ldots, N$. Moreover, assume $\Upsilon_{k+3,j}=\Upsilon_{k,j}$, and $\Upsilon_{e,k+3,j}=\Upsilon_{e,k,j}$ for $k=1,2,3$ and $j=1,\ldots, N$.  The limits in (\ref{just_r33}) and (\ref{just_r34}) exist under our assumptions of exponentially decaying strong mixing coefficients and exponential-type tail probabilities (cf. Theorem 1.7 of \cite{ibragimov1962some}).  %We remark that  $\Upsilon$ and $\Upsilon_{e}$ are $k,j$ specific, but we suppress such dependence. 
We state an additional assumption. 

%\bigskip

\begin{assumption}
\label{assu CLT}
\begin{enumerate}[(i)]
\item For each $k=1,2,3$, $j=1,\ldots, N$, as $T\to \infty$,  
\begin{align*}
&\frac{1}{\sqrt{T}}\sum_{t=1}^{T}\left(
\begin{array}[c]{c}%
\boldsymbol{f}_{t}e_{t,k,j}\\
\boldsymbol{f}_{t}e_{t,k+3,j}%
\end{array}
\right)  \xrightarrow{d} N(\boldsymbol{0}_{28}, \Upsilon_{k,j}).
\end{align*}

\item For each $k=1,2,3$, $j=1,\ldots, N$, as $T\to \infty$,
\begin{align*}
&\frac{1}{\sqrt{T}}\sum_{t=1}^{T}\left(\begin{array}{c}
e_{t,k,j}^2-\sigma_{k,j}^2\\
e_{t,k+3,j}^2-\sigma_{k+3,j}^2
\end{array} \right) \xrightarrow{d} N \left(\boldsymbol{0}_2, \Upsilon_{e,k,j} \right).  
\end{align*}
\end{enumerate}
\end{assumption}

%\bigskip

%Assumption \ref{assu CLT} is almost the same as Assumption F of \cite{baili2016}. 
Assumption \ref{assu CLT}(i) is quite common in the literature (e.g., Assumption F4 of \cite{bai2003}, Assumption F.1 of \cite{baili2016}, and \citet[p.12]{barigozziluciani2022}). Assumption \ref{assu CLT}(ii) is also almost the same as Assumption F.2 of \cite{baili2016}.  Further define
\begin{align*}
\mathcal{M}&:=4D_{14}^{+}(I_{14}\otimes M)\Gamma\left[ \lim_{T\to \infty}\frac{1}{T}\sum_{t=1}^T\sum_{t'=1}^T\left(\mathbb{E}[\boldsymbol{e}_t^{\dagger}\boldsymbol{e}_{t'}^{\dagger \intercal}]\otimes \mathbb{E}[\boldsymbol{f}_t\boldsymbol{f}_{t'}^{\intercal}] \right)\right]\Gamma^{\intercal}(I_{14}\otimes M)D_{14}^{+\intercal}\\
\Upsilon_{C,k,j,11}&:=(\boldsymbol{\lambda}_{k,j}^{\intercal}\otimes I_{14})\Gamma\left[ \lim_{T\to \infty}\frac{1}{T}\sum_{t=1}^T\sum_{t'=1}^T\left(\mathbb{E}[\boldsymbol{e}_t^{\dagger}\boldsymbol{e}_{t'}^{\dagger \intercal}]\otimes \mathbb{E}[\boldsymbol{f}_t\boldsymbol{f}_{t'}^{\intercal}] \right)\right]\Gamma^{\intercal} (\boldsymbol{\lambda}_{k,j}\otimes I_{14})\\
\Upsilon_{C,k,j,22}&:=M^{-1}\left[ \lim_{T\to \infty}\frac{1}{T}\sum_{t=1}^T\sum_{t'=1}^T\mathbb{E}[\boldsymbol{f}_{t}\boldsymbol{f}_{t'}^{\intercal}]\mathbb{E}[e_{t,k,j}e_{t',k,j}]\right] M^{-1}\\
\Upsilon_{C,k,j,12}&:=(\boldsymbol{\lambda}_{k,j}^{\intercal}\otimes I_{14})\Gamma\left[\lim_{T\to \infty}\frac{1}{T}\sum_{t=1}^{T}\sum_{t'=1}^{T}\left( \mathbb{E}[\boldsymbol{e}_t^{\dagger}e_{t',k,j}]\otimes \mathbb{E}[\boldsymbol{f}_t\boldsymbol{f}_{t'}^{\intercal}]\right)  \right]M^{-1},
\end{align*}
and $\Upsilon_{C,k,j,21}:=\Upsilon_{C,k,j,12}^{\intercal}$, for for $k=1,\ldots, 6$ and $j=1,\ldots, N$, where $\boldsymbol{e}_{t}^{\dagger}$ is a $24\times1$ vector consisting of $e_{t,p,q}$ for $p=1,\ldots, 6$ and $q=1,\ldots,4$, and $\Gamma$ is a $196\times336$ matrix, whose elements are known (but complicated) linear functions of elements of (inverted) sub-matrices of $\Lambda$ and $M$ (see (\ref{just_r32}) of the OS). %We remark that  $\Upsilon_{J,11}$, $\Upsilon_{J,22}$ and $\Upsilon_{J,12}$ are $k,j$ specific, but we suppress such dependence.

%\bigskip

\begin{thm}
\label{thm5.2} 
Suppose that Assumptions \ref{assu model}%\ref{assu strong mixing and exponential tail}, \ref{assu eta}, \ref{assu ident1} and 
--\ref{assu estimated within compact}, and \ref{assu CLT} hold.
\begin{enumerate}[(i)]
\item If $\sqrt{T}/N\to 0$ as $N,T\rightarrow \infty$, we have
\[\sqrt{T}(\hat{\sigma}_{k,j}^{2}-\sigma_{k,j}^{2})\xrightarrow{d}N(0,\Upsilon_{e,k,j,11}),\]
for $k=1,\ldots, 6$ and $j=1,\ldots, N$, where $\Upsilon_{e,k,j,11}$ is the $(1,1)$th element of $\Upsilon_{e,k,j}$, which is defined in (\ref{just_r34}). Under the additional assumption of serial independence of $\{e_{t,k,j}\}$, we have the following simplifications $\Upsilon_{e,k,j,11}=\mathbb{E}[e_{t,k,j}^4]-\sigma_{k,j}^4$.

\item If $\sqrt{T}/N\to 0$ as $N,T\rightarrow \infty$, %for $k=1,\ldots,6$, $j=5,\ldots,N$, we have
%\[\sqrt{T}(\hat{\boldsymbol{\lambda}}_{k,j}-\boldsymbol{\lambda}_{k,j})\xrightarrow{d}N\del[2]{\boldsymbol{0}, (\boldsymbol{\lambda}_{k,j}^{\intercal}\otimes I_{14})\Gamma (\Sigma_{ee}^{\dagger}\otimes M)\Gamma^{\intercal}(\boldsymbol{\lambda}_{k,j}\otimes I_{14})+M^{-1}\sigma_{k,j}^2},\]
%and for $k=1,\ldots,6$, $j=1,\ldots,4$, we have
%\[\sqrt{T}(\hat{\boldsymbol{\lambda}}_{k,j}-\boldsymbol{\lambda}_{k,j})\xrightarrow{d}N\del[2]{\boldsymbol{0}, (\boldsymbol{\lambda}_{k,j}^{\intercal}\otimes I_{14})\Gamma (\Sigma_{ee}^{\dagger}\otimes M)\Gamma^{\intercal}(\boldsymbol{\lambda}_{k,j}\otimes I_{14})+M^{-1}\sigma_{k,j}^2+\text{cov}_{k,j}+\text{cov}_{k,j}^{\intercal}},\]

%where $\Gamma$ is a $196\times 336$ matrix, whose elements are known (but complicated) linear functions of elements of (inverted) submatrices of $\Lambda$ and $M$, satisfying
%\[\ve A = \Gamma \times \frac{1}{T}\sum_{t=1}^{T}\left(\boldsymbol{e}_t^{\dagger}\otimes \boldsymbol{f}_t\right)+o_{p}(T^{-1/2}),\]
%where $\boldsymbol{e}_{t}^{\dagger}$ is a $24\times1$ vector consisting of $e_{t,k,j}$ for $k=1,\ldots, 6$ and $j=1,\ldots,4$, $\Sigma_{ee}^{\dagger}$ is the covariance matrix of $\boldsymbol{e}_{t}^{\dagger}$, and $\text{cov}_{k,j}:=(\boldsymbol{\lambda}_{k,j}^{\intercal}\otimes I_{14})\Gamma \sbr[1]{\boldsymbol{\iota}_{k,j}\otimes I_{14}}\sigma_{k,j}^{2}$, with $\boldsymbol{\iota}_{k,j}$ being a $24\times1$ zero vector with its $[4(k-1)+j]$th element replaced by one.

\[\sqrt{T}(\hat{\boldsymbol{\lambda}}_{k,j}-\boldsymbol{\lambda}_{k,j})\xrightarrow{d} N\left( \boldsymbol{0}_{14}, \Upsilon_{C,k,j,11}+\Upsilon_{C,k,j,12}+\Upsilon_{C,k,j,21}+\Upsilon_{C,k,j,22}\right) ,\]
for $k=1,\ldots, 6$ and $j=1,\ldots, N$. Under the additional assumption of serial uncorrelatedness of $\{e_{t,k,j}\}$, we have the following simplifications:
\begin{align*}
\Upsilon_{C,k,j,11}&=(\boldsymbol{\lambda}_{k,j}^{\intercal}\otimes I_{14})\Gamma (\mathbb{E}[\boldsymbol{e}_t^{\dagger}\boldsymbol{e}_t^{\dagger\intercal}]\otimes M)\Gamma^{\intercal}(\boldsymbol{\lambda}_{k,j}\otimes I_{14})\\
\Upsilon_{C,k,j,22}&=M^{-1}\sigma_{k,j}^2\\
\Upsilon_{C,k,j,12}&=(\boldsymbol{\lambda}_{k,j}^{\intercal}\otimes I_{14})\Gamma\left[ \mathbb{E}[\boldsymbol{e}_t^{\dagger}e_{t,k,j}]\otimes M\right]M^{-1}.
\end{align*}

\item If $\sqrt{T}/N\to 0$ as $N,T\rightarrow \infty$, we have
\[\sqrt{T}\vech(\hat{M}-M)\xrightarrow{d}N\left(0, \mathcal{M}\right).\]
Under the additional assumption of serial uncorrelatedness of $\{e_{t,k,j}\}$, we have the following simplification: \begin{align*}
\mathcal{M}=4D_{14}^{+}(I_{14}\otimes M)\Gamma \left(\mathbb{E}[\boldsymbol{e}_t^{\dagger}\boldsymbol{e}_{t}^{\dagger \intercal}]\otimes \mathbb{E}[\boldsymbol{f}_t\boldsymbol{f}_{t}^{\intercal}] \right)\Gamma^{\intercal}(I_{14}\otimes M)D_{14}^{+\intercal}.%\label{just_r38}
\end{align*}

\end{enumerate}
\end{thm}%

We see that the QML-just-identified estimators are asymptotically normal, albeit inefficient. Proposition \ref{prop5.2} is more important than  Theorem \ref{thm5.2} because the rates of convergence are the stepping stone for obtaining the large-sample theories for the QML-all-res estimator, which is the estimator we recommend in practice.

\subsection{QML-All-Res}
\label{sec theory qmle all res}

We will now present the large sample theories of the QML-all-res estimator. We first define the infeasible OLS estimators of $\{\boldsymbol{\lambda}_{k,j,R}\}, \{\sigma^2_{k,j}\}, \phi$, which not only treat $\{\boldsymbol{f}_t\}$ as known, but also incorporate all the restrictions implied by $\Lambda, M$, and $\Sigma_{ee}$.
\begin{align*}
&\underbrace{\hat{\boldsymbol{\lambda}}_{k,j,R}^{ols}}_{4\times 1}:=\left(  \frac{L_k}{T}\sum_{t=1}^{T}\boldsymbol{f}_t\boldsymbol{f}_{t}^{\intercal}L_k^{\intercal}+\frac{L_{k+3}}{T}\sum_{t=1}^{T}\boldsymbol{f}_t\boldsymbol{f}_{t}^{\intercal}L_{k+3}^{\intercal}\right)  ^{-1}\left(  \frac{L_k}{T}\sum_{t=1}^{T}\boldsymbol{f}_t y_{t,k,j}+\frac{L_{k+3}}{T}\sum_{t=1}^{T}\boldsymbol{f}_t y_{t,k+3,j}\right),
\end{align*}
\begin{align*}
\hat{\sigma}_{k,j}^{2, ols}  &:=\frac{1}{2T}\sum_{t=1}^{T}\left[  (y_{t,k,j}-\hat{\boldsymbol{\lambda}}_{k,j}^{ols, \intercal}\boldsymbol{f}_t)^2+(y_{t,k+3,j}-\hat{\boldsymbol{\lambda}}_{k+3,j}^{ols,\intercal}\boldsymbol{f}_t)^2\right] ,
\end{align*}
\begin{align*}
\hat{\phi}^{ols}  &  :=\arg \min_{|\phi|<1}\left[\log|M\left(  \phi \right)
|+\frac{1}{T}\sum_{t=1}^{T}\tr\left(  
\boldsymbol{f}_{t}\boldsymbol{f}_{t}^{\intercal}M(  \phi)^{-1}
\right)\right],
\end{align*}
where $\{L_k\}$ are defined in the sentences following (\ref{main_r12}), and $\hat{\boldsymbol{\lambda}}_{k,j}^{ols}:=L_k^{\intercal}\hat{\boldsymbol{\lambda}}_{k,j,R}^{ols}$. That is, these OLS estimators are the maximizers of the complete quasi-log-likelihood of the two-day representation $\sum_{t=1}^{T}\ell( \boldsymbol{y}_t, \boldsymbol{f}_{t};\boldsymbol{\theta}) $, treating $\{\boldsymbol{f}_t\}_{t=1}^T$ as known and taking into account of all the restrictions implied by $\Lambda$, $M$, and $\Sigma_{ee}$. 
%
%Theorem \ref{thm qmleres1} in Section \ref{sec proof of res rate} of the OS shows that the QML-all-res estimators are asymptotically equivalent to these OLS estimators. In particular, the QML-all-res estimators converge to the OLS estimators sufficiently fast so that they share the same asymptotic distributions if $\sqrt{T}/N\to 0$ as $N,T\to \infty$. This is consistent with \cite{barigozzi2023}. Moreover, 

%\bigskip

\begin{thm}
\label{thm qmleres1} 
Suppose that Assumptions \ref{assu model} -- \ref{assu wt} hold. %In addition, assume $\{\eta_{g,\ell}\}\stackrel{i.i.d.}{\sim}(0,1)$. 
If $\log N/T\to 0$ as $N,T\to \infty$, then we have
\begin{enumerate}[(i)]
\item \[\max_{\substack{1\leq k\leq 3\\1\leq j\leq N}}\Vert \hat{\boldsymbol{\lambda}}_{k,j,R}^{\ast}-\hat{\boldsymbol{\lambda}}_{k,j,R}^{ols}\Vert  =O_{p}(N^{-1})+O_{p}(\Delta),\]

\item \[\max_{\substack{1\leq k\leq 6\\1\leq j\leq N}}|\hat{\sigma}_{k,j}^{2,\ast}-\hat{\sigma}_{k,j}^{2,ols}| =O_{p}(N^{-1})+O_{p}(\Delta),\]

\item \[|\hat{\phi}^{\ast}-\hat{\phi}^{ols}| =O_{p}(N^{-1})+O_{p}(\Delta),\]

\item \[\left[  \frac{1}{6N}\sum_{k=1}^{k=6}\sum_{j=1}^{N}(\hat{\sigma}_{k,j}^{2,*}-\sigma_{k,j}^{2})^{2}\right]  ^{1/2}  =O_{p}(N^{-1})+O_{p}(T^{-1/2}),\]
\end{enumerate}
where $O_{p}(\Delta)=o_p(T^{-1/2})$ in general and $O_{p}(\Delta)=o_p(T^{-1/2}\wedge N^{-1/2})$ under the additional assumption $N\log^3 N/T^2\to 0$. These results also hold for $\{\hat{\boldsymbol{\lambda}}_{k,j,R}^{(p)}\}, \{\hat{\sigma}_{k,j}^{2,(p)}\},\hat{\phi}^{(p)}$ for $p=1,2,3,\ldots$.
\end{thm}

%\bigskip

Theorem \ref{thm qmleres1}(i) uniformly bounds the distance between $\{\hat{\boldsymbol{\lambda}}_{k,j,R}^{\ast}\}$ and their corresponding OLS counterparts; one can interpret Theorem \ref{thm qmleres1}(ii) similarly for $\{\hat{\sigma}_{k,j}^{2,\ast}\}$. Establishing the uniform bound is technically more involved than establishing the pointwise one. Theorem \ref{thm qmleres1}(iii) bounds the distance between $\hat{\phi}^*$ and its OLS corresponding counterpart. We see that the QML-all-res and OLS estimators are asymptotically equivalent so that they share the same asymptotic distributions if $\sqrt{T}/N\to 0$ as $N,T\to \infty$. This is consistent with \cite{barigozzi2023}. Theorem \ref{thm qmleres1}(iv) establishes an average rate of convergence of $\{\hat{\sigma}_{k,j}^{2,\ast}\}$. Corollary \ref{coro rate} in Section \ref{sec coro} of the OS establishes the uniform rates of convergence for the QML-all-res estimators. Such uniform rates of convergence are new to the literature.

%\bigskip

\begin{thm}
\label{thm qmleres} Suppose that Assumptions \ref{assu model} -- \ref{assu CLT} hold. If $\log N/T\to 0$ and $\sqrt{T}/N\to 0$ as $N,T\to \infty$, we have
\begin{enumerate}[(i)]
\item 
\begin{align}
\sqrt{T}(\hat{\boldsymbol{\lambda}}_{k,j,R}^{\ast}-\boldsymbol{\lambda}_{k,j,R})&\xrightarrow{d}N\left(  \boldsymbol{0},\left(  L_{k}ML_{k}^{\intercal}+L_{k+3}ML_{k+3}^{\intercal}\right)  ^{-1}\Upsilon_{L,k,j}\left(  L_{k}ML_{k}^{\intercal}+L_{k+3}ML_{k+3}^{\intercal}\right) ^{-1}  \right) \label{main_r19}
\end{align}
where $\Upsilon_{L,k,j}:=(L_{k}\quad  L_{k+3})\Upsilon_{k,j} (L_{k} \quad L_{k+3})^{\intercal}$, $L_{k}$ is the $4 \times 14$ selection matrix to pick out the non-zero elements in $\boldsymbol{\lambda}_{k,j}$ for $k=1, 2, 3$ and $j=1,\ldots,N$, and $\Upsilon_{k,j}$ is defined in (\ref{just_r33}). Under the additional assumptions of 
\begin{align}
\mathbb{E}[e_{t,k,j}e_{t',k,j}]&=0, \qquad t\neq t', \qquad k=1,\ldots,6,\qquad j=1,\ldots, N, \label{main_r17}\\
\mathbb{E}[e_{t,k,j}e_{t',k+3,j}]&=0, \qquad t,t'=1,\ldots,T, \quad k=1,2,3, \quad j=1,\ldots, N, \label{main_r18}
\end{align}
$\Upsilon_{L,k,j}$ is reduced to $(  L_{k}ML_{k}^{\intercal}+L_{k+3}ML_{k+3}^{\intercal})  \sigma_{k,j}^{2}$. Thus,
\begin{align}
\label{main_r15}
\sqrt{T}(\hat{\boldsymbol{\lambda}}_{k,j,R}^{*}-\boldsymbol{\lambda}_{k,j,R})\xrightarrow{d}N\left(  \boldsymbol{0},(  L_{k}ML_{k}^{\intercal}+L_{k+3}ML_{k+3}^{\intercal})  ^{-1}\sigma_{k,j}^{2} \right) .
\end{align}

\item 
\begin{align}
\sqrt{T}(\hat{\sigma}_{k,j}^{2,*}-\sigma_{k,j}^{2})\xrightarrow{d}N\left(0, \frac{1}{4}\iota^{\intercal}\Upsilon_{e,k,j}\iota \right),\label{main_r7}
\end{align}
for $k=1,2,3$ and $j=1,\ldots,N$, where $\iota=(1,1)^{\intercal}$ and $\Upsilon_{e,k,j}$ is defined in (\ref{just_r34}). Under the additional assumptions of $\{e_{t,k,j}\}$ being independent across $t$ for $k=1,\ldots, 6$ and $j=1,\ldots, N$, and independence between $e_{t,k,j}$ and $e_{t',k+3,j}$ for $t,t'=1,\ldots, T$, $k=1,2,3$, and $j=1,\ldots, N$, $\iota^{\intercal}\Upsilon_{e,k,j}\iota/4$ is reduced to
\begin{align}
\label{main_r16}
%\frac{1}{4}\left(\mathbb{E}[(e_{t,k,j}^2-\sigma_{k,j}^2)^2 ]+ \mathbb{E}[(e_{t,k+3,j}^2-\sigma_{k,j}^2)^2 ] \right). 
v_{e,k,j}:=\frac{1}{4}\left(\mathbb{E}[e_{t,k,j}^4]+\mathbb{E}[e_{t,k+3,j}^4] \right) -\frac{1}{2}\sigma_{k,j}^4.
\end{align}
Under the additional assumption of Gaussianity of $e_{t,k,j}$ and $e_{t,k+3,j}$, (\ref{main_r16}) is reduced to $\sigma_{k,j}^4$.

\item In addition, assume $\{\eta_{g,\ell}\}_{\ell=1}^{\infty}\stackrel{i.i.d.}{\sim}(0,1)$. We have
\[\sqrt{T}(\hat{\phi}^{\ast}-\phi)\xrightarrow{d}N(0,v),\]
\begin{align*}
v=\frac{(1-\phi^{2})^{2}}{(7-5\phi^{2})^{2}}\left[  9-7\phi^{2}+\frac{4(\phi^{12}-\phi^{14}+\phi^{2})}{1-\phi^{12}}+\frac{ \phi^2\gamma_g}{1+\phi^2}+\frac{2\phi^{14}\gamma_g}{1+\phi^2}\cdot \frac{1}{1-\phi^{12}}\right],
\end{align*}
and $\gamma_g$ is the excess kurtosis of $\eta_{g,\ell}$: $\gamma_g:=\mathbb{E}[\eta_{g,\ell}^4]/[\var(\eta_{g,\ell})]^2-3=\mathbb{E}[\eta_{g,\ell}^4]-3$. Under Gaussianity of $\{\eta_{g,\ell}\}_{\ell=1}^{\infty}$, $v$ is reduced to
\[\frac{(1-\phi^{2})^{2}}{(7-5\phi^{2})^{2}}\left[  9-7\phi^{2}+\frac{4(\phi^{12}-\phi^{14}+\phi^{2})}{1-\phi^{12}}\right].\]
\end{enumerate}
\end{thm}

%\bigskip

Theorem \ref{thm qmleres} presents the asymptotic distributions of the QML-all-res estimators. This is a consequence of the OLS estimators having these asymptotic distributions and the asymptotic equivalence of the QML-all-res and OLS estimators (Theorem \ref{thm qmleres1}). $\{\hat{\boldsymbol{\lambda}}_{k,j,R}^*\}$ and $\{\hat{\sigma}_{k,j}^*\}$ are more efficient than the QML-just-identified counterparts because QML-all-res incorporates more information from the model than QML-just-identified does. This is perhaps easier to see under the special cases of (\ref{main_r15}) and (\ref{main_r16}): $\{\hat{\boldsymbol{\lambda}}_{k,j,R}^*\}$ and $\{\hat{\sigma}_{k,j}^*\}$ have smaller asymptotic variances.   %Comparing Theorems \ref{thm qmleres} and \ref{thm5.2} in the OS,  because the asymptotic variances are smaller.  
Moreover, instead of having an estimate for matrix $M$ as QML-just-identified does, we now have an estimate for the only unknown parameter in $M$: $\phi$. For simplicity, to work out consistency and the asymptotic distribution of $\hat{\phi}^{ols}$, we add the assumption of $\{\eta_{g,\ell}\}\stackrel{i.i.d.}{\sim}(0,1)$. This is because in the AR(1) structure of the global factor (\ref{res_r91}), correlated innovations of the global factor make consistent estimation of $\phi$ less straightforward. For example, one needs to specify a parametric structure for the innovations in order to implement the generalised least squares (GLS) estimation. In practice, we recommend to use the QML-all-res estimators.

In the special case (\ref{main_r15}), consistent estimate of the asymptotic variance of $\hat{\boldsymbol{\lambda}}_{k,j,R}^{\ast}$ is 
\begin{align}
\label{main_r25}
(  L_{k}M(\hat{\phi}^*)L_{k}^{\intercal}+L_{k+3}M(\hat{\phi}^*)L_{k+3}^{\intercal})  ^{-1}\hat{\sigma}_{k,j}^{2,*}.
\end{align}
In the special case (\ref{main_r16}), consistent estimate of the asymptotic variance of $\hat{\sigma}_{k,j}^{2,*}$ is
\begin{align}
\label{main_r26}
\hat{v}_{e,k,j}:=\frac{1}{4}\left(  \frac{1}{T}\sum_{t=1}^{T}(y_{t,k,j}-\hat{\boldsymbol{\lambda}}_{k,j,R}^{\ast, \intercal}L_k\hat{\boldsymbol{f}}_t)^4+\frac{1}{T}\sum_{t=1}^{T}(y_{t,k+3,j}-\hat{\boldsymbol{\lambda}}_{k,j,R}^{\ast, \intercal}L_{k+3}\hat{\boldsymbol{f}}_t)^4 \right) -\frac{1}{2}\hat{\sigma}_{k,j}^{4,*}
\end{align}
where $\hat{\boldsymbol{f}}_t$ is defined in (\ref{main_r23}) (see Corollary \ref{coro consistency of standard errors}).  The additional assumptions of (\ref{main_r17}) and (\ref{main_r18}) (or the additional independence assumptions leading to (\ref{main_r16})) say that the idiosyncratic component of every stock return is uncorrelated (or independent) across the calendar days, something quite reasonable in finance. Thus, in the Monte Carlo simulations and empirical studies of this paper, the standard errors of $\{\hat{\boldsymbol{\lambda}}_{k,j,R}^{\ast}\}$ and $\{\hat{\sigma}_{k,j}^{2,*}\}$ are computed based on (\ref{main_r25}) and (\ref{main_r26}), respectively. For the more general cases (\ref{main_r19}) and (\ref{main_r7}), one needs to use a HAC-type estimator (\cite{newey1987simple}, \cite{andrews1991heteroskedasticity}, Section 7 of \cite{barigozziluciani2022}). %The same reasoning applies to $\{\hat{\sigma}_{k,j}^{2,*}\}$. 

\begin{figure}
\centering
\includegraphics[scale=0.65]{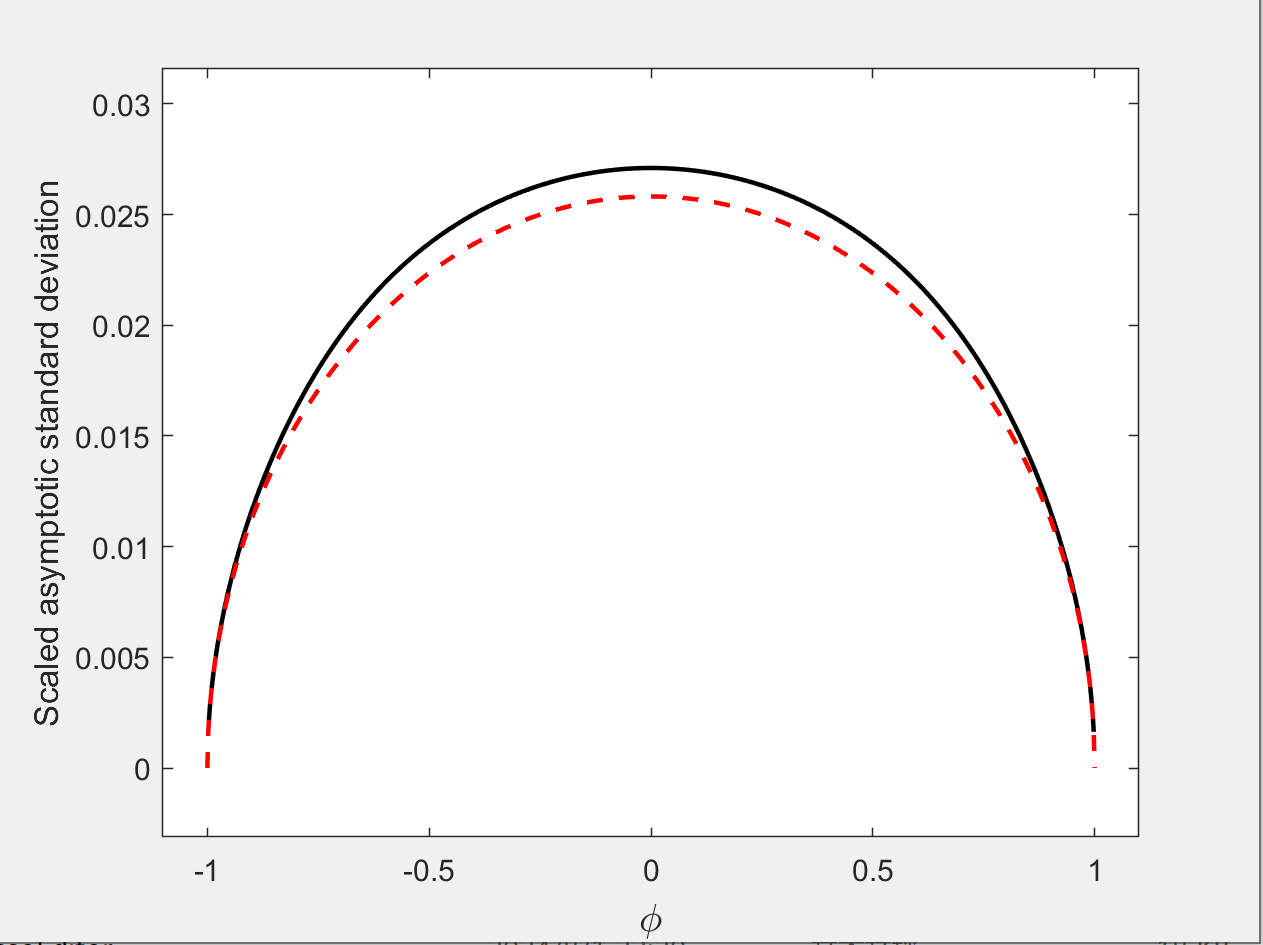}
\caption{{\protect \small Set $T=250$. The black solid line denotes the scaled asymptotic standard deviation of $\hat{\phi}^*$ (or $\hat{\phi}^{ols}$). The red dashed line denotes the scaled asymptotic standard deviation of $\hat{\phi}^{\diamond}$.}}
\label{figure standard error}
\end{figure}

We would like to discuss a bit about the complicated formula of $v$. The scaled asymptotic standard deviation for $\hat{\phi}^*$ (or $\hat{\phi}^{ols}$) is $\sqrt{v/T}$. Recall that both $\hat{\phi}^*$ and $\hat{\phi}^{ols}$ ignore the auto-correlation between $\boldsymbol{f}_t$ and $\boldsymbol{f}_{t-1}$ when setting up the objective functions. If one could observe $\{\boldsymbol{f}_t\}$ and takes the auto-correlation between $\boldsymbol{f}_t$ and $\boldsymbol{f}_{t-1}$ into account when setting up the objective function, the most efficient estimator of $\phi$ would be $\hat{\phi}^{\diamond}:=(\sum_{\ell=1}^{6T}f_{g,\ell-1}^2)^{-1}\sum_{\ell=1}^{6T}f_{g,\ell-1}f_{g,\ell}$. %  $\hat{\phi}^{\diamond}:=(\sum_{\ell=1}^{6T-1}f_{g,\ell}^2)^{-1}\sum_{\ell=1}^{6T-1}f_{g,\ell}f_{g,\ell+1}$. 
It is easy to show that the scaled asymptotic standard deviation for $\hat{\phi}^{\diamond}$ is $\sqrt{(1-\phi^2)/(6T)}$. Figure \ref{figure standard error} plots the scaled asymptotic standard deviations of $\hat{\phi}^*$ (or $\hat{\phi}^{ols}$) and $\hat{\phi}^{\diamond}$ against all the possible values of $\phi$ at $T=250$ under Gaussianity. The black solid line denotes the scaled asymptotic standard deviation of $\hat{\phi}^*$ (or $\hat{\phi}^{ols}$). The red dashed line denotes the scaled asymptotic standard deviation of $\hat{\phi}^{\diamond}$. We see that $\hat{\phi}^{\diamond}$ is marginally more efficient than $\hat{\phi}^*$ (or $\hat{\phi}^{ols}$). The gap between the black solid and red dashed lines is the cost of ignoring the auto-correlation between vectors $\boldsymbol{f}_t$ and $\boldsymbol{f}_{t-1}$. The gap is largest when $|\phi|$ is smallest. This is intuitive because when $|\phi|$ is large, the covariance matrix of $\boldsymbol{f}_t$, $M$, already contains enough information about $\phi$ so ignoring the cross-$t$-period correlation of $\{\boldsymbol{f}_t\}_{t=1}^T$ incurs minimal efficiency loss. In the Monte Carlo simulations and empirical studies of this paper, we estimate $v$ by 
\begin{align*}
\hat{v}&=\frac{(1-\hat{\phi}^{*,2})^{2}}{(7-5\hat{\phi}^{*,2})^{2}}\left[  9-7\hat{\phi}^{*,2}+\frac{4(\hat{\phi}^{*,12}-\hat{\phi}^{*,14}+\hat{\phi}^{*,2})}{1-\hat{\phi}^{*,12}}+\frac{ \hat{\phi}^{*,2}\gamma_g}{1+\hat{\phi}^{*,2}}+\frac{2\hat{\phi}^{*,14}\gamma_g}{1+\hat{\phi}^{*,2}}\cdot \frac{1}{1-\hat{\phi}^{*,12}}\right],\\
\hat{\gamma}_g&:=\frac{1}{6T}\sum_{\ell=1}^{6T}(\hat{f}_{g,\ell}-\hat{\phi}^*\hat{f}_{g,\ell-1})^4-3,
\end{align*}
where $\{\hat{f}_{g,\ell}\}$ are elements of $\{\hat{\boldsymbol{f}}_t\}$ (see (\ref{main_r23}) and Corollary \ref{coro consistency of standard errors}). 

\subsection{The Estimated Factors}

We now provide a large-sample result for the estimated factors: 
\begin{equation}
\hat{\boldsymbol{f}}_t:=[\hat{M}^{*,-1}+\hat{\Lambda}^{*,\intercal}\hat{\Sigma}_{ee}^{*,-1}\hat{\Lambda}^{*}]^{-1}\hat{\Lambda}^{*,\intercal}\hat{\Sigma}_{ee}^{*,-1}\boldsymbol{y}_t,\label{main_r23}
\end{equation}
for $t=1,\ldots, T$. We note that $\hat{\boldsymbol{f}}_t$ is some sort of the cross-sectional average. We make an additional assumption.

%\bigskip

\begin{assumption}
\label{assu factor}
\begin{enumerate}[(i)]
\item For all $t$, $t=1,2,\ldots,T$, and $k=1,2,3$,
\begin{align*}
\mathbb{E}\left\| \frac{1}{\sqrt{6NT}}\sum_{s=1}^{T}\sum_{j=1}^{N}\frac{1}{\sigma^{2}_{k,j}}  \boldsymbol{f}_{s}\left[ e_{s,k,j} e_{t,k,j}-\mathbb{E}[e_{s,k,j} e_{t,k,j}]\right]\right\|&\leq C\\
\mathbb{E}\left\| \frac{1}{\sqrt{6NT}}\sum_{s=1}^{T}\sum_{j=1}^{N}\frac{1}{\sigma^{2}_{k,j}}  \boldsymbol{f}_{s}\left[ e_{s,k+3,j} e_{t,k,j}-\mathbb{E}[e_{s,k+3,j} e_{t,k,j}]\right]\right\|&\leq C.
\end{align*}
For all $t$, $t=1,2,\ldots,T$, and $k=4,5,6$,
\begin{align*}
\mathbb{E}\left\| \frac{1}{\sqrt{6NT}}\sum_{s=1}^{T}\sum_{j=1}^{N}\frac{1}{\sigma^{2}_{k,j}}  \boldsymbol{f}_{s}\left[ e_{s,k,j} e_{t,k,j}-\mathbb{E}[e_{s,k,j} e_{t,k,j}]\right]\right\|&\leq C\\
\mathbb{E}\left\| \frac{1}{\sqrt{6NT}}\sum_{s=1}^{T}\sum_{j=1}^{N}\frac{1}{\sigma^{2}_{k,j}}  \boldsymbol{f}_{s}\left[ e_{s,k-3,j} e_{t,k,j}-\mathbb{E}[e_{s,k-3,j} e_{t,k,j}]\right]\right\|&\leq C.
\end{align*}

\item For all $t$, $t=1,2,\ldots,T$, 
\begin{align*}
\mathbb{E}\left\|\frac{1}{\sqrt{6NT}}\sum_{k=1}^{3}\sum_{j=1}^{N}\sum_{s=1}^{T}\frac{\boldsymbol{\lambda}_{k,j}}{\sigma^4_{k,j}} (e_{s,k,j}^{2}-\sigma_{k,j}^{2}) e_{t,k,j}\right\|&\leq C\\
\mathbb{E}\left\|\frac{1}{\sqrt{6NT}}\sum_{k=1}^{3}\sum_{j=1}^{N}\sum_{s=1}^{T}\frac{\boldsymbol{\lambda}_{k,j}}{\sigma^4_{k,j}} (e_{s,k+3,j}^{2}-\sigma_{k+3,j}^{2}) e_{t,k,j}\right\|&\leq C\\
\mathbb{E}\left\|\frac{1}{\sqrt{6NT}}\sum_{k=4}^{6}\sum_{j=1}^{N}\sum_{s=1}^{T}\frac{\boldsymbol{\lambda}_{k,j}}{\sigma^4_{k,j}} (e_{s,k,j}^{2}-\sigma_{k,j}^{2}) e_{t,k,j}\right\|&\leq C\\
\mathbb{E}\left\|\frac{1}{\sqrt{6NT}}\sum_{k=4}^{6}\sum_{j=1}^{N}\sum_{s=1}^{T}\frac{\boldsymbol{\lambda}_{k,j}}{\sigma^4_{k,j}} (e_{s,k-3,j}^{2}-\sigma_{k-3,j}^{2}) e_{t,k,j}\right\|&\leq C.
\end{align*}

\item Define
\begin{align*}
\Upsilon_F:=\lim_{N\to \infty} \frac{1}{6N}\sum_{k=1}^{6}\sum_{j=1}^{N}\sum_{k'=1}^{6}\sum_{j'=1}^{N}\frac{1}{\sigma_{k,j}^2\sigma_{k',j'}^2}\boldsymbol{\lambda}_{k,j}\boldsymbol{\lambda}_{k',j'}^{\intercal}\Xi_{(k-1)N+j, (k'-1)N+j'}.
\end{align*}
For each $t$, as $N\to \infty$,
\begin{align*}
&\frac{1}{\sqrt{6N}}\Lambda^{\intercal}\Sigma_{ee}^{-1}\boldsymbol{e}_t=\frac{1}{\sqrt{6N}}\sum_{k=1}^{6}\sum_{j=1}^{N}\frac{\boldsymbol{\lambda}_{k,j}e_{t,k,j}}{\sigma_{k,j}^2} \xrightarrow{d} N\left( \boldsymbol{0}, \Upsilon_F\right) . 
\end{align*}
\end{enumerate}
\end{assumption}

Assumption \ref{assu factor}(i) is similar to Assumption E.4 of \cite{baili2016}, while Assumption \ref{assu factor}(ii) is similar to Assumption E.6 of \cite{baili2016}. These two assumptions ensure that $\hat{\Lambda}^{*,\intercal}\hat{\Sigma}_{ee}^{*,-1}\boldsymbol{e}_t/\sqrt{6N}$ is a consistent estimator for $\Lambda^{\intercal}\Sigma_{ee}^{-1}\boldsymbol{e}_t/\sqrt{6N}$ as $N,T\to \infty$. Assumption \ref{assu factor}(iii) is the same as Assumption F.3 of \cite{baili2016}. 

%\bigskip

\begin{thm}
\label{thm clt factors}
Suppose that Assumptions \ref{assu model} -- \ref{assu wt}, \ref{assu factor} hold. If 
\begin{align}
\label{main_r24}
\frac{N\cdot\left[ \log^3 N \vee (\log N)^{1+2/r_1}\vee (\log T)^{2/r_1}\right] }{T^2}\to 0,
\end{align}
%$N\cdot[ \log^3 N\vee \log^2 T] /T^2\to 0$ 
as $N,T\to \infty$, then we have
\begin{equation}
\sqrt{N}\left(\hat{\boldsymbol{f}}_t-\boldsymbol{f}_t \right)\xrightarrow{d}N(0, Q^{-1}\Upsilon_FQ^{-1}), \label{main_r20}
\end{equation}
for $t=1,\ldots, T$.
\end{thm}

%\bigskip

Theorem \ref{thm clt factors} establishes the asymptotic distributions for the estimated factors. These asymptotic distributions are the same as Theorem 2(IC3) of \cite{baili2016} and Proposition 5 of \cite{barigozzi2023}. For estimation of each $\boldsymbol{f}_{t}$, $t=1,\ldots, T$, we do not incur the curse of dimensionality of $T$ because we can rely on $N\to \infty$.
%
%As \cite{barigozzi2023} has pointed out, we need $N\to \infty$ to consistently estimate $\boldsymbol{f}_t$, as $\hat{\boldsymbol{f}}_t$ is some sort of the cross-sectional average, and 
At the same time we need $T\to \infty$ to consistently estimate $\boldsymbol{f}_t$, as we use the QML-all-res estimators $\hat{\boldsymbol{\theta}}^*$, which are consistent only if $T\to \infty$. For both Theorems \ref{thm qmleres} and \ref{thm clt factors} holding, we need $T^{1/2+\delta}\leq N\leq T^{2-\delta}$ for some $0<\delta\leq 3/4$ as $N,T\to \infty$, which is in line with Remark 9 of \cite{barigozziluciani2022}. Under the additional assumption of cross-sectional uncorrelatedness of $\boldsymbol{e}_t$, the asymptotic covariance matrix is reduced to $Q^{-1}$, which could be consistently estimated by $(\frac{1}{6N}\hat{\Lambda}^{*,\intercal}\hat{\Sigma}_{ee}^{*,-1}\hat{\Lambda}^{*})^{-1}$ (see (\ref{f_r15}) of the OS). For the more general case (\ref{main_r20}), one needs to use a HAC-type estimator (\cite{newey1987simple}, \cite{andrews1991heteroskedasticity}, Section 7 of \cite{barigozziluciani2022}).

\begin{assumption}
\label{assu cross sectional Bernstein corollary}
\begin{align*}
&\max_{1\leq t\leq T}\left\| \frac{1}{6N}\Lambda^{\intercal}\Sigma_{ee}^{-1}\boldsymbol{e}_t\right\| =\max_{1\leq t\leq T}\left\| \frac{1}{6N}\sum_{k=1}^{6}\sum_{j=1}^{N}\frac{\boldsymbol{\lambda}_{k,j}e_{t,k,j}}{\sigma_{k,j}^2}\right\|  =o_p(1) ,
\end{align*}
as $N,T\to \infty$.
\end{assumption} 

Assumption \ref{assu cross sectional Bernstein corollary} is a high-level assumption and kind of a uniform law of large numbers along the cross-sectional dimension. To state the primitive conditions for Assumption \ref{assu cross sectional Bernstein corollary}, one needs to introduce some cross-sectional ordering. Under the special case of cross-sectional independence of the idiosyncratic components, one could derive Assumption \ref{assu cross sectional Bernstein corollary} using the Bernstein's inequality.

\begin{lemma}
\label{lemma uniform convergence for f}
Suppose that Assumptions \ref{assu model} -- \ref{assu wt}, \ref{assu cross sectional Bernstein corollary} hold. If (\ref{main_r24}) and 
\begin{align}
\label{f_r23}
\frac{(\log N)^{1/2}(\log (NT))^{1/r_1}}{\sqrt{N\wedge T}}\to 0
\end{align}
hold, as $N,T\to \infty$, then we have $\max_{1\leq t\leq T}\|\hat{\boldsymbol{f}}_t-\boldsymbol{f}_t\|=o_p(1)$.
\end{lemma}

Lemma \ref{lemma uniform convergence for f} establishes uniform consistency for the estimated factors; such uniform consistency is new to the literature. Lemma \ref{lemma uniform convergence for f} is instrumental in establishing consistency of the estimated asymptotic variances of $\{\hat{\sigma}_{k,j}^{2,*}\}$ and $\hat{\phi}^*$: $\hat{v}_{e,k,j}$ and $\hat{v}$, which we state the results as a corollary.

\begin{corollary}
\label{coro consistency of standard errors}
Suppose that Assumptions \ref{assu model} -- \ref{assu CLT}, \ref{assu cross sectional Bernstein corollary} hold. If (\ref{main_r24}), (\ref{f_r23}), and $\sqrt{T}/N\to 0$ hold as $N,T\to \infty$, 
\begin{enumerate}[(i)]
\item we have $\hat{v}_{e,k,j}\xrightarrow{p}v_{e,k,j}$
for $k=1,2,3$, and $j=1,\ldots, N$.

\item In addition, assume $\{\eta_{g,\ell}\}_{\ell=1}^{\infty}\stackrel{i.i.d.}{\sim}(0,1)$. We have $\hat{v}\xrightarrow{p}v$.
\end{enumerate}
\end{corollary}

\section{Monte Carlo Simulations}
\label{sec Monte Carlo}

%\subsection{Gaussian Data}
%\label{sec baseline simu}

In this section, we conduct Monte Carlo simulations to evaluate performance of QML-all-res under an exact factor model setting.\footnote{In Section of \ref{sec simu Non-gaussian} of the OS, we conduct Monte Carlo simulations of the non-Gaussian data.}  We specify $N=100,200$ and $T=250,750$, and a burn-in sample of 250 periods.  Note that $T=250$ and $T=750$ correspond to two and six years' daily data, respectively. The off-diagonal elements of $\Xi$ are set to zero. Assume that $\{\boldsymbol{e}_t\}_{t=1}^T$ are i.i.d. Gaussian across $t$, $\{\eta_{g,\ell}\}_{\ell=0}^{6T}\stackrel{i.i.d.}{\sim}N(0,1)$, and $\{f_{c,\ell}\}\stackrel{i.i.d.}{\sim}N(0,1)$, where $\eta_{g,\ell}$ and $f_{c,\ell}$ are defined in (\ref{res_r91}). Set $\phi=0.2$. The variances of the idiosyncratic components $\{\sigma^2_{k,j}\}$ are all drawn from uniform$[1,1.5]$. For $b=\text{as},\text{eu},\text{am}$, 
\begin{align*}
& \beta_{b,j}=0.6a_{b,j}+0.4d_{b}-0.1, \\ 
& \beta_{b,\text{as\_time},j}^G=0.6a_{b,\text{as\_time},j}+0.4d_{b,\text{as\_time}}-0.1, \\ 
& \beta_{b,\text{eu\_time},j}^G=0.6a_{b,\text{eu\_time},j}+0.4d_{b,\text{eu\_time}}-0.1, \\ 
&\beta_{b,\text{am\_time},j}^G=0.6a_{b,\text{am\_time},j}+0.4d_{b,\text{am\_time}}-0.1 ,
\end{align*}
where $\beta_{b,j}, \beta_{b,\text{as\_time},j}^G, \beta_{b,\text{eu\_time},j}^G, \beta_{b,\text{am\_time},j}^G$ are the $j$th component of $\boldsymbol{\beta}_{b}$, $\boldsymbol{\beta}_{b,\text{as\_time}}^G$, $\boldsymbol{\beta}_{b,\text{eu\_time}}^G$, $\boldsymbol{\beta}_{b,\text{am\_time}}^G$, respectively, for $j=1,\ldots, N$, and $a_{b,j}$, $a_{b,\text{as\_time},j}$, $a_{b,\text{eu\_time},j}$, $a_{b,\text{am\_time},j}$, $d_{b}$, $d_{b,\text{as\_time}}$, $d_{b,\text{eu\_time}}$ and $d_{b,\text{am\_time}}$ are all drawn from uniform$[0,1]$. The econometrician observes those logarithmic 24-hour close-to-close returns. He is aware of the structure of the true model, but does not know the values of those
parameters. 

The number of the Monte Carlo samples is chosen to be 1000. From these 1000
Monte Carlo samples, we calculate the following three quantities for evaluation:

\begin{enumerate}
[(i)]

\item the root mean square error (RMSE) of the quantity of interest,

\item the average of the standard errors (Ave.se) of the estimated elements of the quantity of interest.
%The standard error of the EM estimator is calculated according to Section ...

\item the average of the coverage probabilities (Cove) of the confidence intervals, formed by the point estimates $\pm$ 1.96$\times$the standard errors, of the elements of the quantity of interest. 
\end{enumerate}

Table \ref{table Para_1} presents the RMSE, Ave.se, and Cove of the QML-all-res estimator. % in our simulation. Due to limitations in presenting individual parameter results for each asset, we report the average values across all assets (or within the same market). Specifically, we first compute the RMSE, Ave.se, and Cove for each parameter, then report their average value across all assets (within a particular market). For the loadings, the reported values represent the average across all assets within a particular market. 
We see that the RMSE and Ave.se of the factor loadings are in general quite close to each other. However, in the case of $N=100$, the Ave.se is smaller than the RMSE, indicating that the standard errors are slightly underestimated. When $N$ increases
to 200, the Ave.se gets closer to the RMSE. When $N$ increases from 100 to 200, the Cove improves considerably. %When $N$ increases to 200, the RMSE gets closer to the Ave.se, indicating a smaller finite-sample bias of QML-all-res, when $N$ increases from 100 to 200. The Cove also improves considerably when $N$ increases from 100 to 200. 
We consider size of QML-all-res pretty good as it is well known in the literature that size of the large-dimensional factor models is difficult to control.

%Regarding the idiosyncratic variance $\sigma_{c,i}$, the reported value represents the average across all assets in all markets.
The performance of QML-all-res for $\Sigma_{ee}$ is very good, with the Cove being very close to 0.95. Obtaining a good estimate of $\Sigma_{ee}$ is relatively straightforward. In unreported Monte Carlo simulations, even the QML-just-identified estimator produces a reasonably accurate estimate of $\Sigma_{ee}$. %This is perhaps supported by that $\hat{\sigma}^{2,*}_{k,j}$ (QML-all-res) and $\hat{\sigma}^2_{k,j}$ (QML-just-identified) have the same asymptotic distribution (see Theorems \ref{thm5.2}(i) and \ref{thm qmleres}(ii)). 
When $N=100$, $T=250$ or 750, the standard error of $\hat{\phi}^*$ is under-estimated to some extent, resulting in a slightly lower Cove. However, the Cove improves as $N$ increases to 200.

 \begin{table}[ptb]
 	\centering
 	\begin{tabular}{lcccccccc}
 \toprule[0.3mm] 
 \toprule[0.3mm] 
 &	\multicolumn{3}{c}{$N=100,  T=250$} &       & \multicolumn{3}{c}{$N=100, T=750$} \\
 &\multicolumn{1}{c}{RMSE} & \multicolumn{1}{c}{Ave.se} & \multicolumn{1}{c}{Cove} &       & \multicolumn{1}{c}{RMSE} & \multicolumn{1}{c}{Ave.se} & \multicolumn{1}{c}{Cove} \\
 \midrule[0.3mm]
           $\boldsymbol{\beta}_{\text{as,as\_time}}^G$    & 0.0634  & 0.0528  & 0.8946  &       & 0.0363  & 0.0306  & 0.8997  \\
    $\boldsymbol{\beta}_{\text{as,am\_time}}^G$   & 0.0662  & 0.0555  & 0.8919  &       & 0.0376  & 0.0323  & 0.9002  \\
    $\boldsymbol{\beta}_{\text{as,eu\_time}}^G$   & 0.0683  & 0.0557  & 0.8842  &       & 0.0384  & 0.0323  & 0.8953  \\
    $\boldsymbol{\beta}_{\text{as}}$       & 0.0575  & 0.0520  & 0.9224  &       & 0.0332  & 0.0303  & 0.9242  \\
    $\boldsymbol{\beta}_{\text{eu,eu\_time}}^G$   & 0.0668  & 0.0520  & 0.8757  &       & 0.0381  & 0.0308  & 0.8829  \\
    $\boldsymbol{\beta}_{\text{eu,as\_time}}^G$   & 0.0714  & 0.0549  & 0.8595  &       & 0.0397  & 0.0314  & 0.8766  \\
    $\boldsymbol{\beta}_{\text{eu,am\_time}}^G$   & 0.0649  & 0.0530  & 0.8891  &       & 0.0372  & 0.0311  & 0.8949  \\
    $\boldsymbol{\beta}_{\text{eu}}$   & 0.0613  & 0.0538  & 0.9149  &       & 0.0349  & 0.0313  & 0.9230  \\
    $\boldsymbol{\beta}_{\text{am,am\_time}}^G$   & 0.0631  & 0.0536  & 0.9003  &       & 0.0364  & 0.0306  & 0.9043  \\
    $\boldsymbol{\beta}_{\text{am,eu\_time}}^G$   & 0.0643  & 0.0544  & 0.8982  &       & 0.0370  & 0.0314  & 0.9016  \\
    $\boldsymbol{\beta}_{\text{am,as\_time}}^G$   & 0.0632  & 0.0525  & 0.8943  &       & 0.0362  & 0.0304  & 0.8991  \\
    $\boldsymbol{\beta}_{\text{am}}$   & 0.0614  & 0.0559  & 0.9187  &       & 0.0352  & 0.0319  & 0.9236  \\
     $\Sigma_{ee}$  & 0.0820  & 0.0766  & 0.9254  &       & 0.0471  & 0.0445  & 0.9327  \\
    $\phi$    & 0.0422  & 0.0272  & 0.8520  &       & 0.0210  & 0.0157  & 0.8760  \\
 \midrule[0.3mm]
 \midrule[0.3mm]
 	&	\multicolumn{3}{c}{$N=200, T=250$} &       & \multicolumn{3}{c}{$N=200, T=750$} \\
 	&	\multicolumn{1}{c}{RMSE} & \multicolumn{1}{c}{Ave.se} & \multicolumn{1}{c}{Cove} &       & \multicolumn{1}{c}{RMSE} & \multicolumn{1}{c}{Ave.se} & \multicolumn{1}{c}{Cove} \\
 \midrule[0.3mm]
$\boldsymbol{\beta}_{\text{as,as\_time}}^G$    & 0.0552  & 0.0515  & 0.9301  &       & 0.0319  & 0.0298  & 0.9309  \\
    $\boldsymbol{\beta}_{\text{as,am\_time}}^G$   & 0.0564  & 0.0537  & 0.9328  &       & 0.0323  & 0.0307  & 0.9352  \\
    $\boldsymbol{\beta}_{\text{as,eu\_time}}^G$   & 0.0555  & 0.0522  & 0.9339  &       & 0.0320  & 0.0303  & 0.9359  \\
    $\boldsymbol{\beta}_{\text{as}}$       & 0.0558  & 0.0538  & 0.9402  &       & 0.0322  & 0.0313  & 0.9424  \\
    $\boldsymbol{\beta}_{\text{eu,eu\_time}}^G$   & 0.0570  & 0.0505  & 0.9186  &       & 0.0327  & 0.0296  & 0.9222  \\
    $\boldsymbol{\beta}_{\text{eu,as\_time}}^G$   & 0.0580  & 0.0526  & 0.9207  &       & 0.0334  & 0.0304  & 0.9225  \\
    $\boldsymbol{\beta}_{\text{eu,am\_time}}^G$   & 0.0557  & 0.0510  & 0.9273  &       & 0.0321  & 0.0297  & 0.9303  \\
    $\boldsymbol{\beta}_{\text{eu}}$   & 0.0558  & 0.0531  & 0.9369  &       & 0.0321  & 0.0307  & 0.9396  \\
    $\boldsymbol{\beta}_{\text{am,am\_time}}^G$   & 0.0587  & 0.0520  & 0.9137  &       & 0.0337  & 0.0303  & 0.9167  \\
    $\boldsymbol{\beta}_{\text{am,eu\_time}}^G$   & 0.0587  & 0.0529  & 0.9172  &       & 0.0339  & 0.0305  & 0.9176  \\
    $\boldsymbol{\beta}_{\text{am,as\_time}}^G$   & 0.0585  & 0.0510  & 0.9109  &       & 0.0336  & 0.0294  & 0.9137  \\
    $\boldsymbol{\beta}_{\text{am}}$   & 0.0564  & 0.0530  & 0.9319  &       & 0.0325  & 0.0305  & 0.9344  \\
     $\Sigma_{ee}$  & 0.0810  & 0.0775  & 0.9333  &       & 0.0466  & 0.0452  & 0.9395  \\
    $\phi$    & 0.0309  & 0.0272  & 0.9220  &       & 0.0180  & 0.0157  & 0.9180  \\
 \bottomrule[0.3mm]
 \bottomrule[0.3mm]
 \end{tabular}%
 \caption{{\protect \small RMSE, Ave.se and Cove stand for the root mean square error of, average of the standard errors of the estimated elements of, and average of the coverage probabilities of the confidence intervals, formed by the point estimates $\pm$ 1.96$\times$the standard errors, of the elements of the quantity of interest, respectively.}}%
 	\label{table Para_1}%
 \end{table}%

%  In Figures \ref{figure qq loading} and \ref{figure qq phi}, we present the QQ plots of the 10th element of the standardized QML-all-res of $\boldsymbol{\beta}^G_{\text{as,as\_time}}$ and standardized $\hat{\phi}^*$, respectively, for illustrative purposes. %In the QQ plots, the $x$-axis represents the quantiles of the standard normal, while the $y$-axis represents the quantiles of the standardized estimator. 
% The standardized estimator is calculated as the difference between the estimator and true value, divided by the standard error. These two plots demonstrate that our standardized estimators closely resemble the standard normal.

% \begin{figure}[ptb]
% \centering
% \includegraphics[scale=0.8]{z_plot.pdf}
% \caption{{\protect \small QQ plots: the 10th element of the standardized QML-all-res estimator of $\boldsymbol{\beta}^G_{\text{as,as\_time}}$.}}
% \label{figure qq loading}
% \end{figure}

% \begin{figure}[ptb]
% \centering
% \includegraphics[scale=0.8]{phi_plot.pdf}
% \caption{{\protect \small  QQ plots: the standardized $\hat{\phi}^*$.}}
% \label{figure qq phi}
% \end{figure}

\section{Empirical Application}
\label{sec application}

We apply our model to daily MSCI equity indices of 42 developed and emerging markets.\footnote{The markets include Australia, Austria, Belgium, Brazil, Canada, Chile, China's mainland, Colombia, Czech Republic, Denmark, Finland, France, Germany, Greece, Hong Kong of China, Hungary, India, Indonesia, Ireland, Italy, Japan, Malaysia, Mexico, Netherlands, New Zealand, Norway, Peru, Philippines, Poland, Portugal, Singapore, South Africa, South Korea, Spain, Sweden, Switzerland, Taiwan of China, Thailand, Turkey, United Kingdom, United States.} The data are from Eikon and from May 5th, 2014 (the date when most indices are available) to June 23rd, 2023 ($T=1179$). For most markets, there are 6 indices: Large-Growth, Mid-Growth, Small-Growth, Large-Value, Mid-Value, and Small-Value. We have 212 indices in total. %(the cross-sectional dimension of the two-day representation is 424). 
Dates with more than 50 missing indices are excluded. We categorize these markets into three continents based on their closing times: Asia-Pacific, Europe, and Americas. We exclude the Israeli market from our sample because its closing time is far away from those of Asia-Pacific or Europe. All the returns have been demeaned. To better compare magnitudes of the loadings, we also standardize the returns so that they have unit variances (\citet[p.422]{stockwatson2016}). We shall use QML-all-res.

\subsection{Basic Results}

Table \ref{table_asia} displays the estimated factor loadings and idiosyncratic variances for Asia-Pacific. The indices for China's mainland, Hong Kong of China,  and Japan are all presented, while only the Mid-Value and Mid-Growth indices   are shown for other Asia-Pacific markets in the interest of space. Recall that $\boldsymbol{\beta}_{\text{as, as\_time}}^G$ represents the loading of the global factor during the Asian sub-period, and likewise for others. %Recall that $\boldsymbol{\beta}_{\text{as, as\_time}}^G$, $\boldsymbol{\beta}_{\text{as,am\_time}}^G$, and $\boldsymbol{\beta}_{\text{as,eu\_time}}^G$ represent the loading of the global factor during the Asian sub-period (from the American close on day $s-1$ to the Asian close on day $s$), the American sub-period (from the European close on day $s$ to the American close on day $s$), and the European sub-period (from the Asian close on day $s$ to the European close on day $s$), respectively. Also recall that $\boldsymbol{\beta}_{\text{as}}$ represents the loading of the Asia-Pacific continental factor.
We have made two main observations. First, for the Asia-Pacific markets, the global factor has the largest loadings during the Asian sub-period, with Australia having the largest ones at approximately 0.7. Second, China's mainland and Hong Kong of China have large loadings of the continental factor, around 0.7, whereas Japan has minimal loadings of the continental factor.

\begin{table}[ptb]
\centering
{\footnotesize 
\begin{tabular}
[c]{lcccccc}%
\toprule[0.3mm] \toprule[0.3mm] & $\boldsymbol{\beta}_{\text{as,as\_time},i}^G$ &
$\boldsymbol{\beta}_{\text{as,am\_time},i}^G$ & $\boldsymbol{\beta}_{\text{as,eu\_time},i}^G$ &
$\boldsymbol{\beta}_{\text{as},i}$ & $\sigma^{2}_{\text{as},i}$ & \\
\midrule[0.3mm]
 China's mainland LG&    0.383 (0.011)&    0.055 (0.011)&   -0.144 (0.011)&    0.720 (0.011)&    0.296 (0.009) \\  
China's mainland LV&    0.431 (0.009)&    0.135 (0.009)&   -0.034 (0.009)&    0.753 (0.010)&    0.211 (0.006) \\  
 China's mainland MG&    0.369 (0.011)&    0.067 (0.011)&   -0.142 (0.011)&    0.735 (0.011)&    0.284 (0.008) \\  
 China's mainland MV&    0.387 (0.008)&    0.109 (0.008)&   -0.087 (0.008)&    0.823 (0.008)&    0.141 (0.004) \\  
 China's mainland SG&    0.337 (0.009)&    0.088 (0.009)&   -0.159 (0.008)&    0.818 (0.009)&    0.172 (0.005) \\  
 China's mainland SV&    0.365 (0.008)&    0.115 (0.008)&   -0.091 (0.008)&    0.820 (0.008)&    0.162 (0.005) \\  
  Hong Kong of China LG&    0.403 (0.012)&    0.139 (0.012)&   -0.002 (0.012)&    0.672 (0.012)&    0.351 (0.010) \\  
 Hong Kong of China LV&    0.438 (0.014)&    0.140 (0.014)&    0.136 (0.014)&    0.533 (0.014)&    0.457 (0.013) \\  
 Hong Kong of China MG&    0.392 (0.015)&    0.081 (0.015)&   -0.073 (0.015)&    0.518 (0.016)&    0.554 (0.016) \\  
 Hong Kong of China MV&    0.436 (0.013)&    0.132 (0.013)&    0.075 (0.013)&    0.604 (0.013)&    0.400 (0.012) \\  
 Hong Kong of China SG&    0.378 (0.009)&    0.152 (0.009)&   -0.096 (0.009)&    0.779 (0.010)&    0.204 (0.006) \\  
 Hong Kong of China SV&    0.406 (0.011)&    0.171 (0.010)&    0.010 (0.010)&    0.725 (0.011)&    0.260 (0.008) \\  
Japan LG&    0.347 (0.014)&    0.604 (0.013)&   -0.159 (0.013)&    0.113 (0.014)&    0.430 (0.013) \\  
Japan LV&    0.333 (0.014)&    0.573 (0.014)&   -0.076 (0.014)&    0.058 (0.015)&    0.485 (0.014) \\  
Japan MG&    0.344 (0.014)&    0.598 (0.014)&   -0.166 (0.014)&    0.101 (0.014)&    0.444 (0.013) \\  
Japan MV&    0.345 (0.014)&    0.577 (0.014)&   -0.068 (0.014)&    0.054 (0.014)&    0.469 (0.014) \\  
Japan SG&    0.337 (0.014)&    0.593 (0.014)&   -0.166 (0.014)&    0.092 (0.014)&    0.456 (0.013) \\  
Japan SV&    0.331 (0.014)&    0.582 (0.014)&   -0.089 (0.014)&    0.039 (0.015)&    0.478 (0.014) \\  
\midrule[0.3mm]
 South Korea MV&    0.505 (0.014)&    0.248 (0.014)&    0.242 (0.013)&    0.340 (0.014)&    0.432 (0.013) \\  
 Taiwan of China MV&    0.436 (0.014)&    0.271 (0.014)&    0.177 (0.014)&    0.385 (0.015)&    0.489 (0.014) \\  
Australia MV&    0.731 (0.012)&    0.278 (0.012)&   -0.021 (0.012)&    0.033 (0.012)&    0.329 (0.010) \\  
India MV&    0.517 (0.016)&   -0.030 (0.016)&    0.255 (0.016)&    0.206 (0.017)&    0.625 (0.018) \\  
Indonesia MV&    0.442 (0.016)&    0.025 (0.016)&    0.336 (0.016)&    0.197 (0.017)&    0.632 (0.018) \\  
Malaysia MV&    0.482 (0.016)&    0.076 (0.015)&    0.326 (0.015)&    0.238 (0.016)&    0.564 (0.016) \\  
New Zealand MV&    0.521 (0.017)&    0.193 (0.017)&   -0.019 (0.017)&   -0.021 (0.017)&    0.662 (0.019) \\  
Philippines MV&    0.345 (0.018)&    0.088 (0.018)&    0.272 (0.017)&    0.207 (0.018)&    0.725 (0.021) \\  
Singapore MV&    0.645 (0.014)&    0.144 (0.014)&    0.164 (0.014)&    0.237 (0.014)&    0.441 (0.013) \\  
Thailand MV&    0.521 (0.016)&    0.007 (0.016)&    0.244 (0.016)&    0.218 (0.016)&    0.614 (0.018) \\  
\midrule[0.3mm]
 South Korea MG&    0.481 (0.014)&    0.259 (0.014)&    0.176 (0.014)&    0.359 (0.015)&    0.475 (0.014) \\  
 Taiwan of China MG&    0.414 (0.015)&    0.271 (0.015)&    0.101 (0.015)&    0.373 (0.016)&    0.555 (0.016) \\  
Australia MG&    0.733 (0.011)&    0.321 (0.011)&   -0.124 (0.011)&    0.059 (0.011)&    0.284 (0.008) \\  
India MG&    0.512 (0.016)&    0.004 (0.016)&    0.230 (0.016)&    0.221 (0.017)&    0.630 (0.018) \\  
Indonesia MG&    0.414 (0.017)&    0.010 (0.017)&    0.353 (0.016)&    0.186 (0.017)&    0.651 (0.019) \\  
Malaysia MG&    0.374 (0.017)&    0.102 (0.017)&    0.214 (0.017)&    0.268 (0.018)&    0.705 (0.021) \\  
New Zealand MG&    0.380 (0.018)&    0.215 (0.018)&   -0.094 (0.018)&    0.021 (0.019)&    0.783 (0.023) \\  
Philippines MG&    0.303 (0.018)&    0.087 (0.018)&    0.270 (0.018)&    0.161 (0.018)&    0.771 (0.022) \\  
Singapore MG&    0.640 (0.014)&    0.139 (0.014)&    0.156 (0.013)&    0.271 (0.014)&    0.437 (0.013) \\  
Thailand MG&    0.482 (0.017)&   -0.013 (0.017)&    0.216 (0.017)&    0.201 (0.017)&    0.678 (0.020) \\  
\bottomrule[0.3mm] \bottomrule[0.3mm] 
\end{tabular}} 
\caption{{\protect \footnotesize Selected QML-all-res estimates of the factor loadings and idiosyncratic variances. LG, LV, MG, MV, SG, SV, stand for Large-Growth, Large-Value, Mid-Growth, Mid-Value, Small-Growth and Small-Value, respectively. The standard errors are in parentheses.}}%
\label{table_asia}%
\end{table}

Table \ref{table_euro} displays the estimated factor loadings and idiosyncratic variances for the Europe. All the indices of the UK, France, and Germany are reported, while only the Mid-Value and Mid-Growth indices are reported for other European markets to save space. We have identified three key findings. First, all the European markets have substantial loadings of the global factor during the Asian sub-period, but small loadings of the global factor during the American sub-period. Second, most Value indices exhibit non-negative loadings of the global factor during the European sub-period, with the Large-Value ones having particularly large loadings. Conversely, most Growth indices have negative loadings of the global factor during the European sub-period. Third, most European markets have considerably large loadings of the continental factor, except for South Africa and Turkey. %, which have the largest loadings on the global factor during the Asian sub-period but small loadings on the other factors.

\begin{table}[ptb]
\centering
{\footnotesize 
\begin{tabular}
[c]{lcccccc}%
\toprule[0.3mm] \toprule[0.3mm] & $\boldsymbol{\beta}_{\text{eu,eu\_time},i}^G$ &
$\boldsymbol{\beta}_{\text{eu,as\_time},i}^G$ & $\boldsymbol{\beta}_{\text{eu,am\_time},i}^G$ &
$\boldsymbol{\beta}_{\text{eu},i}$ & $\sigma^{2}_{\text{eu},i}$ & \\
\midrule[0.3mm] 
United Kingdom LG&   -0.028 (0.016)&    0.600 (0.013)&    0.069 (0.012)&    0.520 (0.013)&    0.358 (0.010) \\  
United Kingdom LV&    0.396 (0.012)&    0.536 (0.010)&    0.059 (0.009)&    0.481 (0.010)&    0.209 (0.006) \\  
United Kingdom MG&   -0.085 (0.012)&    0.680 (0.010)&    0.064 (0.009)&    0.587 (0.009)&    0.201 (0.006) \\  
United Kingdom MV&    0.127 (0.013)&    0.602 (0.011)&    0.096 (0.009)&    0.559 (0.010)&    0.233 (0.007) \\  
United Kingdom SG&   -0.085 (0.013)&    0.665 (0.011)&    0.097 (0.010)&    0.563 (0.010)&    0.239 (0.007) \\  
United Kingdom SV&    0.153 (0.013)&    0.600 (0.011)&    0.115 (0.010)&    0.534 (0.010)&    0.242 (0.007) \\  
France LG&   -0.078 (0.011)&    0.619 (0.009)&    0.028 (0.008)&    0.675 (0.009)&    0.175 (0.005) \\  
France LV&    0.379 (0.008)&    0.466 (0.007)&    0.070 (0.006)&    0.671 (0.006)&    0.094 (0.003) \\  
France MG&   -0.087 (0.010)&    0.636 (0.008)&    0.067 (0.007)&    0.684 (0.008)&    0.134 (0.004) \\  
France MV&    0.260 (0.010)&    0.518 (0.009)&    0.084 (0.007)&    0.653 (0.008)&    0.145 (0.004) \\  
France SG&   -0.203 (0.010)&    0.662 (0.009)&    0.049 (0.008)&    0.668 (0.008)&    0.153 (0.004) \\  
France SV&    0.173 (0.010)&    0.576 (0.008)&    0.095 (0.007)&    0.642 (0.008)&    0.141 (0.004) \\  
Germany LG&   -0.150 (0.011)&    0.621 (0.009)&    0.025 (0.008)&    0.689 (0.009)&    0.169 (0.005) \\  
Germany LV&    0.209 (0.009)&    0.515 (0.008)&    0.058 (0.007)&    0.701 (0.007)&    0.124 (0.004) \\  
Germany MG&   -0.329 (0.010)&    0.699 (0.008)&    0.000 (0.007)&    0.653 (0.008)&    0.143 (0.004) \\  
Germany MV&    0.082 (0.012)&    0.558 (0.010)&    0.036 (0.009)&    0.650 (0.010)&    0.217 (0.006) \\  
Germany SG&   -0.300 (0.010)&    0.689 (0.008)&    0.031 (0.007)&    0.669 (0.008)&    0.129 (0.004) \\  
Germany SV&    0.020 (0.009)&    0.652 (0.008)&    0.053 (0.007)&    0.650 (0.007)&    0.123 (0.004) \\  
\midrule[0.3mm] 
Austria MV&    0.372 (0.016)&    0.423 (0.013)&    0.046 (0.011)&    0.519 (0.012)&    0.335 (0.010) \\  
Belgium MV&    0.211 (0.016)&    0.424 (0.013)&    0.060 (0.011)&    0.614 (0.012)&    0.339 (0.010) \\  
Denmark MV&    0.145 (0.022)&    0.354 (0.018)&    0.105 (0.016)&    0.375 (0.017)&    0.665 (0.019) \\  
Finland MV&    0.042 (0.017)&    0.539 (0.014)&    0.028 (0.012)&    0.533 (0.013)&    0.398 (0.012) \\  
Greece MV&    0.164 (0.023)&    0.283 (0.020)&    0.077 (0.017)&    0.318 (0.018)&    0.757 (0.022) \\  
Italy MV&    0.318 (0.013)&    0.372 (0.011)&    0.024 (0.009)&    0.684 (0.010)&    0.234 (0.007) \\  
Netherlands MV&    0.325 (0.012)&    0.461 (0.010)&    0.072 (0.009)&    0.619 (0.010)&    0.212 (0.006) \\  
Norway MV&    0.114 (0.019)&    0.580 (0.016)&    0.016 (0.014)&    0.338 (0.015)&    0.484 (0.014) \\  
Poland MV&    0.131 (0.021)&    0.452 (0.018)&    0.074 (0.015)&    0.337 (0.016)&    0.613 (0.018) \\  
South Africa MV&    0.160 (0.016)&    0.709 (0.014)&   -0.062 (0.012)&    0.143 (0.013)&    0.363 (0.011) \\  
Spain MV&    0.394 (0.014)&    0.413 (0.012)&    0.067 (0.011)&    0.553 (0.011)&    0.289 (0.008) \\  
Sweden MV&   -0.028 (0.013)&    0.646 (0.011)&    0.009 (0.010)&    0.574 (0.010)&    0.251 (0.007) \\  
Switzerland MV&    0.144 (0.013)&    0.580 (0.011)&    0.082 (0.009)&    0.590 (0.010)&    0.220 (0.006) \\  
Turkey MV&    0.118 (0.024)&    0.352 (0.020)&   -0.052 (0.018)&    0.129 (0.019)&    0.815 (0.024) \\  
\midrule[0.3mm] 
Belgium MG&   -0.205 (0.018)&    0.555 (0.015)&    0.027 (0.013)&    0.538 (0.014)&    0.434 (0.013) \\  
Denmark MG&   -0.415 (0.016)&    0.660 (0.013)&    0.037 (0.012)&    0.499 (0.012)&    0.350 (0.010) \\  
Finland MG&   -0.117 (0.019)&    0.490 (0.016)&    0.013 (0.014)&    0.516 (0.015)&    0.513 (0.015) \\  
Ireland MG&   -0.245 (0.021)&    0.460 (0.017)&    0.030 (0.015)&    0.462 (0.016)&    0.595 (0.017) \\  
Italy MG&    0.011 (0.014)&    0.470 (0.012)&    0.007 (0.010)&    0.701 (0.011)&    0.277 (0.008) \\  
Netherlands MG&   -0.353 (0.017)&    0.587 (0.014)&    0.030 (0.012)&    0.550 (0.013)&    0.383 (0.011) \\  
Norway MG&   -0.041 (0.020)&    0.579 (0.016)&    0.015 (0.014)&    0.363 (0.015)&    0.535 (0.016) \\  
Poland MG&   -0.024 (0.019)&    0.557 (0.016)&    0.030 (0.014)&    0.397 (0.015)&    0.527 (0.015) \\  
Portugal MG&    0.070 (0.022)&    0.362 (0.018)&    0.044 (0.016)&    0.417 (0.017)&    0.666 (0.019) \\  
South Africa MG&    0.031 (0.018)&    0.684 (0.015)&   -0.064 (0.014)&    0.164 (0.014)&    0.474 (0.014) \\  
Spain MG&   -0.116 (0.022)&    0.345 (0.019)&    0.032 (0.016)&    0.434 (0.018)&    0.702 (0.020) \\  
Sweden MG&   -0.262 (0.012)&    0.720 (0.010)&   -0.005 (0.009)&    0.574 (0.010)&    0.211 (0.006) \\  
Switzerland MG&   -0.418 (0.012)&    0.722 (0.010)&    0.040 (0.009)&    0.556 (0.010)&    0.217 (0.006) \\  
Turkey MG&    0.083 (0.024)&    0.373 (0.020)&   -0.052 (0.018)&    0.150 (0.019)&    0.806 (0.023) \\   
\bottomrule[0.3mm] \bottomrule[0.3mm] &  &  &  &  &  &
\end{tabular}
} \caption{{\protect \footnotesize Selected QML-all-res estimates of the factor loadings and idiosyncratic variances. LG, LV, MG, MV, SG, SV, stand for Large-Growth, Large-Value, Mid-Growth, Mid-Value, Small-Growth and Small-Value, respectively. The standard errors are in parentheses.}}%
\label{table_euro}%
\end{table}

Table \ref{table_US} presents the estimated factor loadings and idiosyncratic variances for Americas. We see that the United States and Canada have substantial loadings of the global factor during the American sub-period, whereas their loadings of the continental factor are negligible. On the other hand, the emerging American markets tend to have significantly large loadings of the global factor during the Asian sub-period. Brazil also has large loadings of the continental factor.

\begin{table}[ptb]
{\small \centering
\begin{tabular}
[c]{lcccccc}%
\toprule[0.3mm] \toprule[0.3mm] & $\boldsymbol{\beta}_{\text{am,am\_time},i}^G$ &
$\boldsymbol{\beta}_{\text{am,eu\_time},i}^G$ & $\boldsymbol{\beta}_{\text{am,as\_time},i}^G$ &
$\boldsymbol{\beta}_{\text{am},i}$ & $\sigma^{2}_{\text{am},i}$ & \\
\midrule[0.3mm]
United States LG&    1.018 (0.012)&   -0.342 (0.012)&    0.327 (0.010)&    0.093 (0.009)&    0.188 (0.005) \\  
United States LV&    0.999 (0.010)&    0.096 (0.010)&    0.277 (0.008)&    0.055 (0.008)&    0.133 (0.004) \\  
United States MG&    1.075 (0.008)&   -0.330 (0.008)&    0.348 (0.007)&    0.078 (0.007)&    0.096 (0.003) \\  
United States MV&    1.027 (0.007)&    0.092 (0.007)&    0.302 (0.006)&    0.059 (0.006)&    0.072 (0.002) \\  
United States SG&    1.096 (0.007)&   -0.201 (0.008)&    0.316 (0.006)&    0.069 (0.006)&    0.078 (0.002) \\  
United States SV&    1.021 (0.008)&    0.132 (0.008)&    0.268 (0.007)&    0.051 (0.006)&    0.089 (0.003) \\ 
Brazil LG&    0.286 (0.008)&    0.138 (0.008)&    0.414 (0.007)&    0.758 (0.007)&    0.097 (0.003) \\  
Brazil LV&    0.253 (0.011)&    0.274 (0.012)&    0.379 (0.010)&    0.672 (0.009)&    0.185 (0.005) \\  
Brazil MG&    0.279 (0.007)&    0.119 (0.007)&    0.423 (0.006)&    0.779 (0.006)&    0.069 (0.002) \\  
Brazil MV&    0.262 (0.009)&    0.193 (0.009)&    0.408 (0.007)&    0.741 (0.007)&    0.108 (0.003) \\  
Brazil SG&    0.284 (0.008)&    0.120 (0.008)&    0.413 (0.006)&    0.777 (0.006)&    0.080 (0.002) \\  
Brazil SV&    0.248 (0.006)&    0.165 (0.006)&    0.415 (0.005)&    0.788 (0.005)&    0.053 (0.002) \\  
Canada LG&    0.702 (0.015)&   -0.132 (0.015)&    0.558 (0.012)&    0.066 (0.011)&    0.298 (0.009) \\  
Canada LV&    0.602 (0.013)&    0.222 (0.013)&    0.524 (0.011)&    0.053 (0.010)&    0.233 (0.007) \\  
Canada MG&    0.656 (0.014)&   -0.110 (0.014)&    0.606 (0.012)&    0.071 (0.011)&    0.273 (0.008) \\  
Canada MV&    0.524 (0.015)&    0.115 (0.015)&    0.575 (0.013)&    0.057 (0.012)&    0.317 (0.009) \\  
Canada SG&    0.582 (0.016)&   -0.122 (0.016)&    0.600 (0.013)&    0.066 (0.013)&    0.357 (0.010) \\  
Canada SV&    0.531 (0.015)&    0.082 (0.015)&    0.598 (0.012)&    0.057 (0.011)&    0.300 (0.009) \\  
Chile LG&    0.300 (0.021)&    0.031 (0.021)&    0.504 (0.017)&    0.145 (0.016)&    0.599 (0.017) \\  
Chile LV&    0.210 (0.021)&    0.150 (0.021)&    0.479 (0.018)&    0.149 (0.016)&    0.608 (0.018) \\  
Chile MG&    0.167 (0.022)&    0.064 (0.023)&    0.448 (0.019)&    0.117 (0.018)&    0.713 (0.021) \\  
Chile MV&    0.168 (0.022)&    0.138 (0.022)&    0.443 (0.019)&    0.135 (0.017)&    0.677 (0.020) \\  
Chile SG&    0.165 (0.022)&    0.099 (0.022)&    0.487 (0.018)&    0.103 (0.017)&    0.658 (0.019) \\  
Chile SV&    0.181 (0.021)&    0.094 (0.021)&    0.502 (0.018)&    0.106 (0.017)&    0.637 (0.019) \\  
Colombia LG&    0.331 (0.021)&    0.190 (0.021)&    0.397 (0.018)&    0.159 (0.016)&    0.607 (0.018) \\  
Colombia LV&    0.335 (0.019)&    0.220 (0.020)&    0.441 (0.016)&    0.174 (0.015)&    0.530 (0.015) \\  
Colombia SV&    0.199 (0.023)&    0.150 (0.023)&    0.377 (0.019)&    0.130 (0.018)&    0.726 (0.021) \\  
Mexico LG&    0.317 (0.019)&    0.146 (0.020)&    0.491 (0.016)&    0.181 (0.015)&    0.529 (0.015) \\  
Mexico LV&    0.311 (0.019)&    0.118 (0.019)&    0.536 (0.016)&    0.169 (0.015)&    0.499 (0.015) \\  
Mexico MG&    0.265 (0.020)&    0.008 (0.020)&    0.554 (0.017)&    0.137 (0.016)&    0.570 (0.017) \\  
Mexico MV&    0.272 (0.019)&    0.113 (0.019)&    0.557 (0.016)&    0.153 (0.015)&    0.501 (0.015) \\  
Mexico SG&    0.312 (0.019)&    0.079 (0.019)&    0.557 (0.016)&    0.165 (0.015)&    0.497 (0.014) \\  
Mexico SV&    0.305 (0.019)&    0.085 (0.019)&    0.547 (0.016)&    0.168 (0.015)&    0.508 (0.015) \\  
Peru LG&    0.530 (0.020)&    0.106 (0.020)&    0.375 (0.017)&    0.113 (0.016)&    0.550 (0.016) \\  
\bottomrule[0.3mm] \bottomrule[0.3mm] 
\end{tabular}
} \caption{{\protect \footnotesize Selected QML-all-res estimates of the factor loadings and idiosyncratic variances. LG, LV, MG, MV, SG, SV, stand for Large-Growth, Large-Value, Mid-Growth, Mid-Value, Small-Growth and Small-Value, respectively. The standard errors are in parentheses.}}%
\label{table_US}%
\end{table}

Overall, a large number of developing markets have large loadings of the global factor during the Asian sub-period. This means that the international news happened during the period between the American close on day $s-1$ and the Asian close on day $s$ has a substantial impact on the stock returns of these developing markets. Table \ref{table_phi} reports $\hat{\phi}^*=0.2374$ which is significantly different from 0. This implies that the weak-form efficient market hypothesis does not hold and it takes more than a $1/3$ day for some international news to fully unfold or dissipate.\footnote{Figure \ref{figure factor} of the OS plots the estimated factors.} %We see that the global factor had a very large negative shock during the March of 2020, a period where the global stock markets crashed in the fear of Covid-19.}

\begin{table}[ptb]
\centering 
		\begin{tabular}
			[c]{lcc}%
			\toprule[0.3mm] 
\toprule[0.3mm]
& Estimated value & Standard error\\
\midrule[0.3mm]
			 $\phi$ & 0.2374  & 0.0129 \\
\midrule[0.3mm] 
\midrule[0.3mm]
\end{tabular}
\caption{{\protect \footnotesize The QML-all-res estimate of $\phi$, measuring the autocorrelation of the global factor. }}%
	\label{table_phi}%
\end{table}

\subsection{Integration}
%\label{sec application}

Financial market integration holds significant importance in the field of international finance. We adopt a methodology similar to \cite{2009Global}  to assess market integration.  \cite{2009Global} used the R-squared as a measure of integration, which is the proportion of variance explained by the explanatory variables. We define the proportion of variance explained by the global factor as the \textit{global integration index}, while that explained by the continental factor as the \textit{regional integration index}. The \textit{non-integration index} is calculated as: $1-\text{global integration index}-\text{regional integration index}$.

Table \ref{table_inte_asian} presents an analysis of integration of Asia-Pacific. We find that China's mainland and Hong Kong of China play a significant role in driving the Asia-Pacific regional integration. Specifically, most these indices demonstrate high regional integration, with approximately 50\% of their variances explained by the continental factor. On the contrary, Japan shows high global integration, with over 50\% of its variance explained by the global factor, but displays very low regional integration. In other words, Japanese indices are more globally integrated, while Chinese indices are more regionally integrated.

Table \ref{table_inte_euro} presents an analysis of integration of Europe. The developed European markets, such as the UK, France, and Germany, exhibit higher levels of both the global and regional integration. On the other hand, the less developed European markets, such as Greece, Portugal, and Turkey, show lower levels of integration. For instance, the non-integration index of Turkish Mid-Growth index is as high as 81.28\%, indicating that 81.28\% of the variance of the Turkish Mid-Growth index is explained by neither the global nor continental factor.

Table \ref{table_inte_us} presents an analysis of integration of Americas. The US market demonstrates the highest level of global integration, with the global factor accounting for approximately 90\% of its variance. However, the US market lacks regional integration, meaning it is not influenced by the American continental factor. Within Americas, Brazil holds the most prominent role in regional integration.

\begin{table}[ptb]
\centering
\begin{tabular}
[c]{lccc}%
\toprule[0.3mm] 
\toprule[0.3mm]
& Global integration & Regional integration & Non-integration  \\
\toprule[0.3mm] 
China's mainland LG&   0.1812&   0.5213&   0.2975 \\  
 China's mainland LV&   0.2376&   0.5558&   0.2066 \\  
 China's mainland MG&   0.1722&   0.5425&   0.2853 \\  
 China's mainland MV&   0.1898&   0.6709&   0.1393 \\  
 China's mainland SG&   0.1573&   0.6704&   0.1723 \\  
 China's mainland SV&   0.1745&   0.6654&   0.1601 \\  
Hong Kong of China LG&   0.2156&   0.4412&   0.3432 \\  
 Hong Kong of China LV&   0.2825&   0.2746&   0.4428 \\  
 Hong Kong of China MG&   0.1832&   0.2667&   0.5501 \\  
 Hong Kong of China MV&   0.2566&   0.3550&   0.3884 \\  
 Hong Kong of China SG&   0.2000&   0.5989&   0.2011 \\  
 Hong Kong of China SV&   0.2358&   0.5113&   0.2530 \\  
Japan LG&   0.5723&   0.0123&   0.4154 \\  
Japan LV&   0.5261&   0.0032&   0.4706 \\  
Japan MG&   0.5609&   0.0099&   0.4291 \\  
Japan MV&   0.5438&   0.0028&   0.4534 \\  
Japan SG&   0.5497&   0.0081&   0.4422 \\  
Japan SV&   0.5346&   0.0015&   0.4639 \\  
\midrule[0.3mm]
 South Korea MV&   0.4797&   0.1097&   0.4106 \\  
 Taiwan of China MV&   0.3885&   0.1420&   0.4695 \\  
Australia MV&   0.6930&   0.0010&   0.3060 \\  
India MV&   0.3491&   0.0414&   0.6095 \\  
Indonesia MV&   0.3460&   0.0377&   0.6163 \\  
Malaysia MV&   0.4001&   0.0548&   0.5452 \\  
New Zealand MV&   0.3618&   0.0004&   0.6378 \\  
Philippines MV&   0.2462&   0.0422&   0.7115 \\  
Singapore MV&   0.5310&   0.0530&   0.4160 \\  
Thailand MV&   0.3584&   0.0460&   0.5956 \\  
\midrule[0.3mm]
 South Korea MG&   0.4244&   0.1225&   0.4531 \\  
 Taiwan of China MG&   0.3320&   0.1338&   0.5342 \\  
Australia MG&   0.7318&   0.0032&   0.2649 \\  
India MG&   0.3398&   0.0474&   0.6128 \\  
Indonesia MG&   0.3283&   0.0337&   0.6380 \\  
Malaysia MG&   0.2414&   0.0700&   0.6886 \\  
New Zealand MG&   0.2329&   0.0005&   0.7666 \\  
Philippines MG&   0.2143&   0.0255&   0.7602 \\  
Singapore MG&   0.5179&   0.0693&   0.4128 \\  
Thailand MG&   0.2970&   0.0396&   0.6634 \\  
\bottomrule[0.3mm] \bottomrule[0.3mm] 
\end{tabular}
\caption{{\protect\footnotesize Integration of the Asia-Pacific markets: ``Global integration'' represents the global integration index, ``Regional integration'' represents the regional integration index, and ``Non-integration'' represents the non-integration index.}}
\label{table_inte_asian}%
\end{table}

\begin{table}[ptb]
\centering
\begin{tabular}
[c]{lccc}%
\toprule[0.3mm] 
\toprule[0.3mm] 
& Global integration & Regional integration & Non-integration  \\
\toprule[0.3mm] 
United Kingdom LG&   0.3885&   0.2634&   0.3481 \\  
United Kingdom LV&   0.5765&   0.2224&   0.2011 \\  
United Kingdom MG&   0.4758&   0.3312&   0.1930 \\  
United Kingdom MV&   0.4679&   0.3046&   0.2274 \\  
United Kingdom SG&   0.4677&   0.3031&   0.2292 \\  
United Kingdom SV&   0.4880&   0.2771&   0.2349 \\  
Germany LG&   0.3793&   0.4576&   0.1631 \\  
Germany LV&   0.3943&   0.4835&   0.1222 \\  
Germany MG&   0.4763&   0.3920&   0.1316 \\  
Germany MV&   0.3674&   0.4180&   0.2146 \\  
Germany SG&   0.4673&   0.4133&   0.1193 \\  
Germany SV&   0.4671&   0.4130&   0.1198 \\  
France LG&   0.3865&   0.4433&   0.1702 \\  
France LV&   0.4771&   0.4321&   0.0908 \\  
France MG&   0.4190&   0.4517&   0.1294 \\  
France MV&   0.4436&   0.4150&   0.1414 \\  
France SG&   0.4333&   0.4219&   0.1449 \\  
France SV&   0.4602&   0.4019&   0.1379 \\  
\midrule[0.3mm] 
Austria MV&   0.4153&   0.2603&   0.3244 \\  
Belgium MV&   0.2957&   0.3704&   0.3339 \\  
Denmark MV&   0.2093&   0.1378&   0.6529 \\  
Finland MV&   0.3258&   0.2809&   0.3933 \\  
Greece MV&   0.1532&   0.0996&   0.7472 \\  
Italy MV&   0.3126&   0.4582&   0.2292 \\  
Netherlands MV&   0.4239&   0.3706&   0.2055 \\  
Norway MV&   0.4058&   0.1134&   0.4808 \\  
Poland MV&   0.2844&   0.1118&   0.6038 \\  
SouthAfrica MV&   0.6095&   0.0207&   0.3698 \\  
Spain MV&   0.4300&   0.2931&   0.2770 \\  
Sweden MV&   0.4293&   0.3240&   0.2467 \\  
Switzerland MV&   0.4437&   0.3408&   0.2155 \\  
Turkey MV&   0.1609&   0.0168&   0.8223 \\  
\midrule[0.3mm] 
Belgium MG&   0.3076&   0.2774&   0.4150 \\  
Denmark MG&   0.4637&   0.2232&   0.3131 \\  
Finland MG&   0.2381&   0.2600&   0.5018 \\  
Ireland MG&   0.2277&   0.2037&   0.5686 \\  
Italy MG&   0.2372&   0.4877&   0.2751 \\  
Netherlands MG&   0.3698&   0.2779&   0.3522 \\  
Norway MG&   0.3441&   0.1298&   0.5260 \\  
Poland MG&   0.3267&   0.1553&   0.5180 \\  
Portugal MG&   0.1662&   0.1730&   0.6608 \\  
South Africa MG&   0.4942&   0.0271&   0.4787 \\  
Spain MG&   0.1247&   0.1853&   0.6900 \\  
Sweden MG&   0.4934&   0.3088&   0.1978 \\  
Switzerland MG&   0.5327&   0.2750&   0.1923 \\  
Turkey MG&   0.1645&   0.0227&   0.8128 \\  
\bottomrule[0.3mm] \bottomrule[0.3mm] 
\end{tabular}
\caption{{\protect\footnotesize Integration of the European markets: ``Global integration'' represents the global integration index, ``Regional integration'' represents the regional integration index, and ``Non-integration'' represents the non-integration index.}}
\label{table_inte_euro}%
\end{table}

\begin{table}[ptb]
\centering
\begin{tabular}
[c]{lccc}%
\toprule[0.3mm] 
\toprule[0.3mm]
& Global integration & Regional integration & Non-integration  \\
\midrule[0.3mm] 
United States LG&   0.8531&   0.0064&   0.1405 \\  
United States LV&   0.9013&   0.0022&   0.0965 \\  
United States MG&   0.9258&   0.0044&   0.0698 \\  
United States MV&   0.9463&   0.0025&   0.0512 \\  
United States SG&   0.9406&   0.0034&   0.0559 \\  
United States SV&   0.9349&   0.0018&   0.0632 \\  
Brazil LG&   0.3431&   0.5619&   0.0950 \\  
Brazil LV&   0.3853&   0.4361&   0.1787 \\  
Brazil MG&   0.3368&   0.5952&   0.0680 \\  
Brazil MV&   0.3577&   0.5364&   0.1060 \\  
Brazil SG&   0.3312&   0.5908&   0.0780 \\  
Brazil SV&   0.3378&   0.6099&   0.0523 \\  
Canada LG&   0.7340&   0.0038&   0.2621 \\  
Canada LV&   0.7903&   0.0025&   0.2072 \\  
Canada MG&   0.7503&   0.0046&   0.2451 \\  
Canada MV&   0.7023&   0.0030&   0.2947 \\  
Canada SG&   0.6675&   0.0040&   0.3284 \\  
Canada SV&   0.7172&   0.0030&   0.2798 \\  
Chile LG&   0.3899&   0.0206&   0.5895 \\  
Chile LV&   0.3748&   0.0221&   0.6031 \\  
Chile MG&   0.2745&   0.0136&   0.7120 \\  
Chile MV&   0.3078&   0.0180&   0.6741 \\  
Chile SG&   0.3324&   0.0107&   0.6569 \\  
Chile SV&   0.3531&   0.0113&   0.6356 \\  
Colombia LG&   0.3912&   0.0245&   0.5843 \\  
Colombia LV&   0.4618&   0.0291&   0.5092 \\  
Colombia SV&   0.2654&   0.0168&   0.7179 \\  
Mexico LG&   0.4519&   0.0319&   0.5162 \\  
Mexico LV&   0.4822&   0.0280&   0.4898 \\  
Mexico MG&   0.4170&   0.0185&   0.5644 \\  
Mexico MV&   0.4811&   0.0233&   0.4956 \\  
Mexico SG&   0.4847&   0.0268&   0.4885 \\  
Mexico SV&   0.4727&   0.0279&   0.4994 \\  
Peru LG&   0.4855&   0.0116&   0.5029 \\  

\bottomrule[0.3mm] \bottomrule[0.3mm] 
\end{tabular}
\caption{{\protect\footnotesize Integration of the American markets: ``Global integration'' represents the global integration index, ``Regional integration'' represents the regional integration index, and ``Non-integration'' represents the non-integration index.}}
\label{table_inte_us}%
\end{table}

\subsection{Market Connectivity under Different VIX Levels}
\label{sec VIX}

We now investigate whether market connectivity increases when the US VIX increases. We divide the sample into two subsamples according to the level of the VIX. Since our QML estimators ignore the cross-$t$-period correlation, we could assign non-consecutive time units to a sub-sample, something we cannot do when using the MLE. Nevertheless, in order to maintain some extent of temporal continuity, we assign time units from a month whose monthly VIX exceeds the median monthly VIX into the high-VIX sub-sample, while that is below the median monthly VIX into the low-VIX sub-sample.\footnote{If the first day of a time unit falls into the high-VIX month and the second day of the time unit falls into the low-VIX month, we will assign the time unit to the high-VIX sub-sample, and so forth.}

We conduct this analysis for two separate time intervals: one from May 5th, 2014 to Nov 29th, 2018 (the before sub-sample), and the other from Nov 30th, 2018 to June 23rd, 2023 (the after sub-sample). %  prior to 2018-11-29 (which marks the midpoint of the entire sample period spanning from May 5th, 2014 to June 23rd, 2023) and the other subsequent to 2018-11-29. 
Consequently, we have a total of four sub-samples: Before-HighVix, Before-LowVix, After-HighVix, and After-LowVix, where Before-HighVix is the intersection of the before and high-VIX sub-samples, and so forth. %and pertains to the period prior to 2018-11-29 characterized by high VIX values, and similarly for the other sub-samples.
We note that the HighVix sub-sample consists of the time units mainly after Nov 29th, 2018, while the LowVix sub-sample mainly consists of the time units mainly before Nov 29th, 2018.\footnote{Figure \ref{figure vix} in Section \ref{sec figure factor} of the OS depicts these sub-samples.} Thus, it is possible that the HighVix and LowVix sub-samples differ in aspects other than the VIX. For example, the Covid-19 happened in the after sub-sample. We argue that the Before-HighVix and Before-LowVix sub-samples are more alike in aspects other the VIX, similarly for the After-HighVix and After-LowVix sub-samples. That is the reason we partition the sample into the four sub-samples.

%We partition time intervals into before and after periods because if we do not perform this division,  These differences could include variations in policies across different years, such as trade conflicts in the initial period and the emergence of the pandemic in the later period. In comparison, the Before-HighVix and Before-LowVix sub-samples may exhibit relatively smaller differences in aspects beyond VIX. The same reasoning applies to the after-HighVix and after-LowVix subsamples.

We estimate the model for each of the four sub-samples. Table \ref{table_inte_did} presents the average value of the non-integration or global integration indices of all the markets within each continent for each sub-sample. The average value of the non-integration indices is bigger in the Before-LowVix sub-sample than in the Before-HighVix sub-sample for each continent; the pattern also holds for the After-LowVix and After-HighVix sub-samples. Moreover, the average value of global integration indices is bigger in the Before-HighVix sub-sample than in the Before-LowVix sub-sample for each continent, similarly for the After-HighVix and After-LowVix sub-samples (except for Europe). This implies that most markets become more integrated during periods of high VIX.

Table \ref{table_loading_did} displays the average absolute value of loadings of a continent for a sub-period of a day for a sub-sample. %The first row stands for the average value of the absolute loading on the global factor during the Asia-subperiod ($\boldsymbol{\beta}_{\text{as,as\_time}}^G$), and others likewise.  We present two main findings. 
For all three continents, the global factor has larger average absolute values of the loadings during the Asian sub-period in the high VIX sub-samples than in the low VIX sub-samples. This means that the international news happened during the Asian sub-period becomes more influential when the VIX increases. %On the other hand, the loadings of the continental factors remain roughly the same in terms of either the average absolute values reported in Table \ref{table_loading_did} or the individual absolute values (available upon request).

\begin{table}[ptb]
\centering
\begin{tabular}[c]{lcccc}%
\toprule[0.3mm]
\toprule[0.3mm]
Non-integration & Before-LowVix  & Before-HighVix  & After-LowVix & After-HighVix  \\
\midrule[0.3mm] 
Asia-Pacific &    0.6221 &   0.3557  & 0.5399  &  0.3790\\  
Europe &   0.4597 &   0.2647   & 0.4103  &  0.2932 \\  
Americas &    0.4003  &  0.2804  &  0.4353  &  0.2209 \\  
\toprule[0.3mm]
\toprule[0.3mm]
Global integration & Before-LowVix  & Before-HighVix  & After-LowVix & After-HighVix  \\
\midrule[0.3mm] 
Asia-Pacific &     0.2741 &   0.5760  &  0.4118  &  0.5286\\  
Europe &      0.3328  &  0.5753   & 0.5157   & 0.4303\\  
Americas &    0.5028   & 0.6321   & 0.3923  &  0.6873\\  
\bottomrule[0.3mm]  
\bottomrule[0.3mm] 
\end{tabular}
\caption{{\protect\small The average value of the non-integration or global integration indices of all the markets within each continent for each sub-sample.}}
\label{table_inte_did}%
\end{table}

\begin{table}[ptb]
\centering
\begin{tabular}[c]{lcccc}%
\toprule[0.3mm]
\toprule[0.3mm]
 & Before-LowVix  & Before-HighVix  & After-LowVix & After-HighVix  \\
\midrule[0.3mm] 
Ave. of $|\boldsymbol{\beta}_{\text{as,as\_time}}^G|$ &   0.2446&   0.4993&   0.4030&   0.6840 \\  
Ave. of $|\boldsymbol{\beta}_{\text{eu,as\_time}}^G|$&    0.3404&   0.6136&   0.4947&   0.5139  \\  
Ave. of $|\boldsymbol{\beta}_{\text{am,as\_time}}^G|$&   0.3397&   0.6019&   0.3347&   0.7195\\  
Ave. of $|\boldsymbol{\beta}_{\text{as,eu\_time}}^G|$&    0.2032&   0.2037&   0.4510&   0.2057\\  
Ave. of $|\boldsymbol{\beta}_{\text{eu,eu\_time}}^G|$&  0.2494&   0.3977&   0.4018&   0.2050 \\  
Ave. of $|\boldsymbol{\beta}_{\text{am,eu\_time}}^G|$&  0.4180&   0.1783&   0.1915&   0.2581 \\  
Ave. of $|\boldsymbol{\beta}_{\text{as,am\_time}}^G|$&  0.1910&   0.4526&   0.1573&   0.1130 \\  
Ave. of $|\boldsymbol{\beta}_{\text{eu,am\_time}}^G|$&    0.0638&   0.0538&   0.2394&   0.1411 \\  
Ave. of $|\boldsymbol{\beta}_{\text{am,am\_time}}^G|$&  0.1369&   0.4371&   0.2961&   0.5072 \\  
Ave. of $|\boldsymbol{\beta}_{\text{as}}|$&  0.2334&   0.2185&   0.1667&   0.2436\\  
Ave. of $|\boldsymbol{\beta}_{\text{eu}}|$&   0.4387&   0.4456&   0.2260&   0.5223  \\  
Ave. of $|\boldsymbol{\beta}_{\text{am}}|$&  0.2250&   0.2288&   0.3257&   0.2984 \\ 
\bottomrule[0.3mm]  
\bottomrule[0.3mm]  
\end{tabular}
\caption{{\protect\small The average absolute value of loadings of a continent for a sub-period of a day for a sub-sample.}}
\label{table_loading_did}%
\end{table}

%Firstly, all continents, espcially Asian and European, exhibit larger loadings on the Asia-subperiod global factor in both 'Before-Highvix' and 'After-Highvix', indicating an increased significance of global information from the Asia-subperiod during high VIX periods. Secondly, the loadings on regional factors remain consistent between high VIX and low VIX periods, not only in terms of the average values but also for individual values(reported in OS).

%In summary, our findings indicate that the global integration increases when the VIX increases. This is primarily driven by the increase of the magnitudes of loadings of the global factor during the Asian sub-period.

In Section \ref{sec out of sample performance} of the OS, we compare our model with various competing models (the penalized maximum likelihood (PML) method of \cite{bailiao2016}, the sparse orthogonal factor regression (SOFAR) method of \cite{uematsu2019}, the principal component analysis (PCA) etc.) in terms of the forecasting performance of the next-day returns. Our model has the best forecasting performance. We also compare our model with various competing models (PML, linear shrinkage estimator of \cite{ledoitwolf2004}, quadratic shrinkage estimator  of \cite{ledoitwolf2022} etc.) in terms of the prediction accuracy of the covariance matrix. Although our model is not specifically designed for estimation or prediction of covariance matrices, our model performs reasonably well.\footnote{Section \ref{sec esti cova matrix} of the OS proves the average Frobenius norm consistency of the matrix estimator formed using QML-all-res.}

%\section{Conclusion}
%\label{sec conclusion}

%We introduce a confirmatory dynamic factor model to model a large number of daily stock returns across different time zones. This model accounts for the impacts of both the global and continental factors, and provides means to assess the efficient market hypothesis and levels of international and regional integration. We extend \cite{baili2016}'s results to this large confirmatory dynamic factor model. This has practical importance because economic or financial theories often imply confirmatory factor models whose identification schemes do not comply with those of \cite{baili2016}. %We apply our model to the 42 markets and find that the US is the most integrated market, and the global market exhibits higher integration when the VIX is large.

\bibliographystyle{chicago}
\bibliography{over}

\end{document}